\font \sevenrm=cmr7

\newenvironment{disarray}%
 {\everymath{\displaystyle\everymath{}}\array}%
 {\endarray}
\documentclass[10pt,a4paper]{article}
\usepackage{amssymb}
\newtheorem{thm}{Theorem}
\newtheorem{defn}{Definition}
\newtheorem{cor}{Corollary}
\newtheorem{lem}{Lemma}

\newtheorem{prop}{Proposition} 
\newtheorem{ex}{Example}
\newtheorem{rk}{Remark}

\usepackage{amssymb}
\usepackage{color}
\usepackage{colordvi}
\usepackage{amsmath}
\newcommand{\ignore}[1]{}

\newcommand{\field}[1]{{\mathbb #1}}
\newcommand{\R}{\field{R}}
\newcommand{\C}{\field{C}}

\newcommand{\Z}{\field{Z}}
\newcommand{\Q}{\field{Q}}
\newcommand{\N}{\field{N}}
\newcommand\cutoffsum{\mathop{-\hskip -4mm\sum}\limits}

\newcommand\altcutoffsum{\mathop{-\hskip -3mm\sum}\limits}
\newcommand\altscutoffsum{\mathop{=\hskip -3.5mm\sum}\limits}

\def\shu{\joinrel{\!\scriptstyle\amalg\hskip -3.1pt\amalg}\,}
\def\sshu{\joinrel{\hskip 2pt\scriptscriptstyle\amalg\hskip -2.5pt\amalg}\,}
\def \endsquare{ $\sqcup \!\!\!\! \sqcap$ }
\def \mop#1{\mathop{\hbox{\rm #1}}\nolimits}
\def \smop#1{\mathop{\hbox{\sevenrm #1}}\nolimits}

\def \mopl#1{\mathop{\hbox{\rm #1}}\limits}
\def \smop#1{\mathop{\hbox{\sevenrm #1}}\nolimits}

\def\sshu{\joinrel{\hskip 2pt\scriptscriptstyle\amalg\hskip -2.5pt\amalg}\,}
\def\cutoffint{-\hskip -10pt\int}

\def\otherterm#1{{\it#1}}
\def \e {{\epsilon}}
\def \restr#1{\mathstrut_{\textstyle |}\raise-6pt\hbox{$\scriptstyle #1$}}
\def \srestr#1{\mathstrut_{\scriptstyle |}\hbox to
-1.5pt{}\raise-4pt\hbox{$\scriptscriptstyle #1$}}

\def \Ci {{C^\infty}}

\begin{document}

\title{\bf  Nested  sums of symbols  and
  renormalised multiple zeta values}

\author{  Dominique MANCHON, Sylvie
PAYCHA }
\maketitle
\section*{Abstract}
We define
  renormalised nested sums of symbols on $\R$ 
which  satisfy stuffle relations. For appropriate symbols these give rise 
to  renormalised Euler-Zagier-Hoffman multiple zeta (and Hurwitz zeta)
functions which satisfy stuffle relations at all arguments. We show
the rationality of renormalised multiple zeta values at nonpositive integer arguments. These results generalise to radial symbols
on $\R^n$ giving rise to a higher-dimensional analog of multiple zeta functions.
\section*{Acknowledgements}
The second author would like to thank Pierre Cartier for interesting
discussions around multiple zeta functions as well as Steven Rosenberg for
pointing out reference \cite{GSW}. Thanks to both of them for giving
the second author the opportunity to present part of the results at a
preliminary stage of their development. The first author thanks Kurusch
Ebrahimi-Fard and Michael E. Hoffman for crucial references and details
about quasi-shuffle products, and Herbert Gangl
for checking the table at the end of section \ref{sect:entiersnegatifs} with a PARI
routine. Both authors are indebted to Li Guo
for various interesting discussions and for bringing to our attention and
commenting for us  his own
work with Bin Zhang \cite{GZ} which  also uses a Birkhoff factorisation to renormalise
multiple zeta functions. We also greatly appreciate comments sent  by   Shlomo
Sternberg and   Hoang Ngoc  Minh.\\

The second author also thanks the Max Planck Institute in Bonn  where parts of  this
 paper were completed; she is specially grateful to Bin Zhang for discussions  during her stay at the
  Max Planck Institute in the summer 2006, 
 which were essential to clarify some  fundamental points in the paper.\\

We
 thank the referees for their numerous and very helpful comments and
 suggestions which greatly contributed in improving the final presentation.
\eject
\section*{Introduction}
Multiple zeta functions\footnote{also known under the purely greek-rooted expression  ``polyzeta
   functions'' (P. Cartier). We refer the reader to e.g.  \cite{H},
   \cite{Z},\cite{W}, \cite{ENR}, \cite{ CEMP}, \cite{Mi}, \cite{Zu} among a
   long list of articles on algebraic relations obeyed by multiple zeta functions.}
\begin{equation}
  \zeta(s_1, \ldots, s_k):= \sum_{0<n_k< \cdots< n_1}
\frac{1}{n_1^{s_1}} \cdots \frac{1}{n_k^{s_k}}
\end{equation}
converge when $\sum_{j=1}^m\mop{Re}s_j>m$ for all $m\in\{1,\ldots ,k\}$ (see
\cite{G} Theorem 2.25, \cite{Zh})\footnote{In our opinion, there is a misprint in both
  references concerning this convergence condition.} and
can be meromorphically extended to $\C^k$, with singularities located at:
\begin{eqnarray*}
s_1&=&1,\\
s_1+s_2&=&2,1,0,-2,-4,-6,\ldots,\\
s_1+\cdots+ s_j&\in&\Z\cap]-\infty, j], \ j=3,4,\ldots,k
\end{eqnarray*}
(\cite{AET}, see also \cite{Zh} and \cite{G} Theorem 2.25). They satisfy quasi-shuffle (or ``stuffle'') relations on the convergence domain\footnote{We shall not
  consider shuffle relations and regularisation relations for multiple zeta
  values, which, at least to our knowledge, only make
  sense at positive integer arguments.}, an example of which reads:
\begin{equation}
\zeta(s_1)\, \zeta(s_2)= \zeta(s_1, s_2) +\zeta(s_2,s_1) +\zeta(s_1+s_2).
\end{equation}
The aim of this article is threefold:
\begin{itemize}
\item to extend  multiple zeta
functions to all  arguments $s_i$
using a renormalisation procedure \`a
la Connes and Kreimer,  in such a way that the renormalised values satisfy
stuffle relations,
\item to prove the rationality of multiple zeta values at nonpositive integer
  arguments $(s_1,\ldots ,s_k)$, and the very much related holomorphicity of the maps $z\mapsto
\zeta(s_1-c_1z, \ldots, s_k-c_k z)$ for
complex numbers $c_1, \ldots, c_k$ with positive real part.
\item to provide  natural   higher
 dimensional analogs of multiple zeta functions which we renormalise in a similar
 manner.
\end{itemize}
 In our approach,  the  requirement that   multiple zeta
functions should also obey    stuffle  relations at  all 
arguments
   boils down  to requiring that
 certain maps define characters on  Hopf algebras given by the tensor product
 algebras of {\sl classical symbols with constant coefficients\/} (see paragraph \ref{par:logpoly} for
 the definition) equipped  with the
  stuffle product\footnote{In  \cite{GZ},  the authors have a similar
  approach   using a Birkhoff factorisation to renormalise
multiple zeta functions. Our construction should relate to theirs via a Mellin
transform.}. Let $\chi$ be a smooth
cut-off function which vanishes around the origin and such that $\vert \xi\vert \geq
1 \Rightarrow \chi(\xi)=1$, and consider the functions $\sigma_s(\xi):= \chi(\xi) \, \vert
\xi\vert^{-s}$, which are elementary examples of classical symbols. It was
shown in \cite{MP} that, as $N\to+\infty$, for any $(s_1,\ldots,s_k)\in\C^k$ the sum:
$$\sum_{0<n_k<\cdots <n_1\le N}\sigma_{s_1}(n_1)\cdots\sigma_{s_k}(n_k)$$
admits an asymptotic expansion in powers of $N$ and $\log
N$, as does the iterated integral:
$$\int_{r\leq  x_k\leq \cdots\leq 
  x_1\leq N}
\sigma_{s_1}(x_1) \cdots \sigma_{s_k} (x_k)\, d\, x_1\cdots d\, x_k$$
(for some $r>0$) to which the above sum is related via an iterated use of the
Euler-MacLaurin formula. We define the {\sl cut-off nested sum\/} and
the {\sl cut-off nested integral\/} of the tensor product
$\sigma_{s_1}\otimes\cdots\otimes\sigma_{s_k}$ as the respective constant
terms in $N$, which we denote by
$$\cutoffsum_<^{\smop{Chen}}\sigma_{s_1}\otimes\cdots\otimes\sigma_{s_k}
\hskip 8mm\hbox{ and
}\hskip 8mm\cutoffint_r^{\smop{Chen}}\sigma_{s_1}\otimes\cdots\otimes\sigma_{s_k}$$
respectively, in honour of K.T. Chen \cite{Ch}. When $\mop{Re }(s_1+\cdots
+s_m)>m$ for all $m$ in $\{1,\ldots ,k\}$, then
the above cut-off nested integral and sum both converge and coincide with the
ordinary nested integral and sum respectively, which obey shuffle,
resp. stuffle relations. The nested sum is nothing but $\zeta(s_1,\ldots,s_k)$
in this case. However these shuffle and stuffle properties do not extend in a
 straightforward manner
 to cut-off integrals and sums; nevertheless, a holomorphic 
 regularisation of the symbols (which
 perturbs their order)  provides shuffle, resp. stuffle relations
 by analytic
 continuation as identities of meromorphic maps. These constructions lead to
 algebra morphisms $\Phi_r^{\cal R}$
 (see Proposition \ref{prop:Chencharacter}) and ${\Psi}^{\cal R},$ (see Proposition
 \ref{prop:meroChensums2}) defined on Hopf algebras ${\cal H}_0$ and ${\cal H}$
 respectively, with values in
meromorphic functions, where ${\cal H}_0$ (resp. ${\cal H}$) is the tensor
algebra built on the space of classical symbols on
$\R$  (or a suitable subalgebra of it) equipped 
with the shuffle product (resp. a suitable stuffle product) and
deconcatenation coproduct.\\

Let us make this statement more precise: on the one hand, a holomorphic regularisation 
${\cal R}: \sigma\mapsto \sigma(z)$ (e.g. Riesz or  dimensional
regularisation, see section \ref{sect:regul})
extends multiplicatively to a regularisation $\tilde{\cal R}$ on the tensor algebra of  classical
symbols defined by  $$\tilde {\cal R}(\sigma_1\otimes
\cdots\otimes\sigma_k):={\cal R}(\sigma_1)\otimes \cdots \otimes{\cal
  R}(\sigma_k),$$    leading  to a
meromorphic map
$$\Phi_r^{\cal R}(\sigma_1\otimes\cdots\otimes\sigma_k): z
\mapsto\cutoffint_{  r}^{{\rm Chen}}\tilde{\cal R}\left(\sigma_1\otimes\cdots\otimes \sigma_k\right)(z),$$
with poles of order $\leq k$. The regularisation  ${\tilde{\cal R}}$ is
compatible with the shuffle product:
$${\tilde {\cal R}}(\sigma\shu\tau)=\tilde{\cal
   R}(\sigma)\shu\tilde{\cal R}(\tau).$$
On the other hand, one can extend  a regularisation
${\cal R}$ on classical symbols to  a regularisation ${\tilde{\cal R}}^\star$ on the tensor algebra of  classical
symbols which is compatible with the above mentioned stuffle product $\star$,
i.e. such that
$${\tilde {\cal R}^\star}(\sigma\star\tau)=\tilde{\cal
   R}^\star(\sigma)\star\tilde{\cal R}^\star(\tau).$$
 by ``twisting'' the
 regularisation $\tilde {\cal R}$ by the {\sl Hoffman Hopf algebra
   isomorphism\/} from
 ${\cal H}_0$ to ${\cal H}$ (see section \ref{sect:9} and \cite{H2}). This leads to a meromorphic map
\begin{equation}
{\Psi}^{\cal R}(\sigma_1\otimes \cdots\otimes \sigma_k): z
\longmapsto\cutoffsum^{{\rm Chen}}_{<}  {\tilde{\cal
    R}}^\star\left(\sigma_1\otimes\cdots\otimes \sigma_k\right)(z).
\end{equation}
The meromorphicity property for cut-off iterated
integrals and the upper bound on the order of the poles follow from properties of cut-off integrals of holomorphic families of 
{\sl log-polyhomogeneous symbols\/}, a natural generalisation of classical
symbols which includes positive integer powers of logarithms (see paragraph
\ref{par:logpoly} and \cite{L}). Similar meromorphicity properties and a similar pole structure  hold for cut-off Chen
sums (see section \ref{sect:cutoffchen}), which can be derived from the
corresponding properties of cut-off Chen integrals 
via  the Euler-MacLaurin formula. We prove this way  (Proposition
\ref{thm:meroChensums})  that 
the maps \begin{equation}\label{eq:Chensumzi} (z_1, \ldots, z_k)\mapsto \cutoffsum^{{\rm Chen}}_{<}{\cal R}
\left(\sigma_1\right)(z_1)\otimes\cdots\otimes {\cal R}
\left(\sigma_k\right)(z_k)
\end{equation}
are meromorphic with poles on a discrete set of  hyperplanes
$$q\sum_{i=1}^j z_i\in -\sum_{i=1}^j \alpha_i(0)+[-j,
\infty[\cap \Z $$
in the case when $\alpha_i(z)=\alpha_i(0)+qz$ is the holomorphic order of the
symbol $\sigma_i(z)$. When $\sigma_j(z)(\xi)=\chi(\xi)\vert\xi\vert^{-z}$ (and hence
$\alpha_i(z)=-z$), we recover  (Corollary \ref{cor:meroChensums})
 the
pole structure $\sum_{i=1}^jz_i\in ]-\infty ,j ]\cap \Z$ of multiple zeta
functions recalled at the beginning of the introduction, up to the
specificities in depths one and  two which are discussed later in the paper
 (Theorem
\ref{thm:recurrence}).\\

The {\sl Birkhoff factorisation\/} of these
characters (see Theorems \ref{thm:renChenint}, \ref{thm:renChensums}) then
provides renormalised values  at $z=0$ for any tensor product
$\sigma_1\otimes\cdots\otimes \sigma_k$ of classical symbols, 
   which  obey the same
shuffle, resp. stuffle relations.
Applying this to 
 $\sigma_i(\xi)= \vert \xi\vert^{-s_i} \chi(\xi)$ and a regularisation
  procedure
 ${\cal R}_\mu(\sigma_i )(z)= \chi(\xi)\, \vert \xi\vert^{-s_i-z-\mu z^2} $
 \quad ($\mu \in \R$), leads on the one hand to renormalised continuous analogs  of 
multiple zeta functions
(see Corollary  \ref{cor:rencontmultiplezeta}) 
 which  obey the expected
shuffle relations and on the other hand to renormalised multiple zeta values  
 $\zeta^{{\cal R}}(s_1, \ldots, s_k)$
 (depending
on the regularisation procedure ${\cal R}$) at all arguments
$s_i$,
  which indeed verify the
stuffle relations. Multiple zeta functions are sometimes considered in their ``weak inequality version'':
\begin{equation}
\overline\zeta(s_1, \ldots, s_k):= \sum_{0<n_k\le \cdots\le n_1}
\frac{1}{n_1^{s_1}} \cdots \frac{1}{n_k^{s_k}},
\end{equation}
which are  simultaneously renormalised along the same lines.\\

We show that renormalised multiple zeta functions $\overline\zeta(s_1, \ldots, s_k)$ and
  $\zeta(s_1, \ldots, s_k)$ at non-positive  integer arguments are
  rational numbers which do not depend on the regularisation chosen, so that
  we can drop the superscript ${\cal R}$ in this case (Theorem \ref {thm:fpRieszsumhom}). This is due
  to the holomorphicity of the maps $z\mapsto\zeta(s_1+c_1z,\ldots, s_k+c_kz)$
  at these arguments which holds for any complex numbers $c_1,\ldots,c_k$. Li Guo and Bin Zhang
  \cite{GZ} have proposed another construction for
  stuffle-compatible multiple zeta
  values. Their construction however makes sense only when all the arguments have the same sign (positive or negative). They also use Birkhoff decomposition in an
  essential way, but starting from a different Hopf algebra, and the values are
  rational numbers at non-positive arguments as well. \\

Proposals for multiple zeta values at nonpositive arguments
  have already been formulated (see e.g. \cite{G}, \cite{AET}) but, ref. \cite{GZ}
  excepted, without taking quasi-shuffle relations into account. Let us also
  mention the work of Jean Ecalle \cite{E} on multiple zeta functions, in
  which a brief discussion of values at
  nonpositive arguments is outlined, in the language of mould calculus.\\

We give an explicit formula (in terms of
  Bernoulli numbers) for double
  zeta values (equation (\ref{multiple
    zeta2arg})). Our double
  zeta values at nonpositive arguments turn out to coincide with those proposed by
  S. Akiyama, S. Egami and Y. Tanigawa at the end of their paper as:
\begin{equation}
\zeta(-a_1,-a_2)=\mopl{lim}_{z\to 0}\zeta(-a_1+z,-a_2+z)
\end{equation}
(\cite{AET}, Remark 2
  therein), but this phenomenon does not survive in depth $k\ge 3$. Our table does not
  coincide with the table given in \cite{GZ}, except at the diagonal arguments
  $(-a,-a)$ (as a consequence of the stuffle relations) and at points where
  the double zeta function is holomorphic in the two variables, i.e. $(-a,-b)$
  with odd $a+b$.\\

It turns out that the method outlined above can be applied to more general
situations: the Hurwitz multiple zeta functions\footnote{The case $k=1$ gives
  the Hurwitz
  zeta function with parameter $q=v+1$ since $\zeta(s; \, v)= \zeta(s,
  v+1)$. This non-standard convention somewhat simplifies the notations when
  dealing with Hurwitz and ordinary multiple zeta functions simultaneously.}:
\begin{equation}
\zeta(s_1, \ldots, s_k;\, v_1, \ldots, v_k)
:= 
\sum_{0<n_k< \cdots< n_1}
\frac{1}{(n_1+v_1)^{s_1}} \cdots \frac{1}{(n_k+v_k)^{s_k}},
 \quad v_i\geq 0,
\end{equation}
are renormalised exactly along the same lines, by considering symbols $\sigma_{s,v}(\xi):= \chi(\xi) \, (\vert
\xi\vert+v)^{-s}$ for some non negative real number $v$. Here again,
at depth $2$  we have \begin{equation}\label{eq:AET}
\zeta(-a_1,-a_2;v)=\mopl{lim}_{z\to 0}\zeta(-a_1+z,-a_2+z;v).
\end{equation} Furthermore, 
the values of renormalised Hurwitz multiple zeta functions at non-positive
integer arguments:
$$\zeta(-a_1, \ldots, -a_k;\, v, \ldots, v)$$
are rational numbers when $v\in\Q$, but the rationality issue for
non-identical rational
parameters $v_1,\ldots v_k$ is still an open question. Our proof for $v_1=\cdots=v_k=v$ relies
on the algebra structure on the set of functions $\{\sigma_{-a, v}(x)=
(x+v)^{a},x>0, 
\quad a\in \N\}$. We show that our renormalized Hurwitz multiple zeta values
at nonpositive integer arguments satisfy the two following identities:
\begin{equation}\label{eq:hurwitz1}
\zeta(-a_1,\ldots ,-a_k;v+1)=\zeta(-a_1,\ldots
,-a_k;v)-(v+1)^{a_k}\zeta(-a_1,\ldots ,-a_{k-1};v+1)
\end{equation}
and
\begin{equation}\label{eq:hurwitz2}
\frac{d}{dv}\zeta(-a_1,\ldots ,-a_k;v)=\sum_{j=1}^k a_j\zeta(-a_1,\ldots
,-a_{j-1},-a_j+1,-a_{j+1},\ldots ,-a_k; v),
\end{equation}
which are well-known for arguments in the domain of convergence. Whereas  the multiple zeta values $\zeta^{\rm alt}(-a_1,-a_2; v):=\mopl{lim}_{z\to 0}\zeta(-a_1+z,-a_2+z, \ldots,
-a_k+z; v)$ considered in  \cite{AET} for non positive integers $a_i$ do not
satisfy stuffle relations, they do  satisfy the relations \eqref{eq:hurwitz1}
and \eqref{eq:hurwitz2} as a result of the corrseponding identities of meromorphic functions (with respect to the variables $s_j$):
\begin{equation}
\zeta(s_1,\ldots ,s_k;v+1)=\zeta(s_1,\ldots ,s_k;v)-(v+1)^{-s_k}\zeta(s_1,\ldots ,s_{k-1};v+1)
\end{equation}
and:
\begin{equation}
\frac{d}{dv}\zeta(s_1,\ldots ,s_k;v)=-\sum_{j=1}^k s_j\zeta(s_1,\ldots
,s_{j-1},s_j+1,s_{j+1},\ldots ,s_k; v).
\end{equation}
In the last part of the paper, we extend the renormalisation procedure
described above to 
higher-dimensional iterated sums on $\R^n$ of the following type:
$$\sum_{1\le|n_k|<\cdots<|n_1|}\sigma_1\otimes\cdots\otimes \sigma_k.$$
Here $|.|$ stands for the supremum norm on $\R^n$, and the $\sigma_j$'s are radial functions $f_j\circ |.|$ where $f_j$ is an
ordinary symbol on $\R$ with support in $]0,+\infty[$. Such a
nested sum can be explicitly expressed as an ordinary
(one-dimensional) iterated sum, namely:
\begin{equation*}
\sum_{1\le|n_k|<\cdots<|n_1|}\sigma_1\otimes\cdots\otimes \sigma_k=\sum_{1\le n_k<\cdots<n_1}A_nf_1\otimes\cdots\otimes A_nf_k,
\end{equation*}
where $A_n(t)=(2t+1)^n-(2t-1)^n$ is the natural interpolation of the number of points in the intersection of the sphere of
radius $t$ with $\Z^n$.\\

This yields a natural way to define higher-dimensional renormalized
Hurwitz multiple zeta  values 
$\zeta_n(s_1,\ldots,
s_k;v_1,\ldots,v_k)$ in terms of renormalized nested sums of one-dimensional symbols. When $v_1=\cdots= v_k=v$, these turn out to be linear combinations with coefficients given by
  polynomial expressions in $v$  of  renormalized one-dimensional  Hurwitz multiple zeta
 values, which converge for $\mop{Re}(s_1+\cdots+s_m)>nm$ for any $m\in\{1,\ldots,k\}$.
  An interesting but probably  difficult problem is how to build similar
renormalised   higher dimensional multiple zeta values with the supremum norm
replaced by the Euclidean norm, the main difficulty arising from the lack of an exact
formula for  the number of integer points in the euclidean sphere  of given
radius, centered at
zero. 
\tableofcontents
\vfill\eject
\section{NESTED INTEGRALS OF SYMBOLS }
Throughout this part as well as parts 2 and 3, $\vert\cdot \vert$ stands for a
continuous  norm on $\R^n$, smooth outside the origin.
\subsection{Cut-off integrals of  log-polyhomogeneous symbols}
\subsubsection{Log-polyhomogeneous symbols}\label{par:logpoly}
We partly follow \cite{L}, Sections 2 and 3. For any complex
number  $\alpha$ and any non-negative integer $k$ let us  denote by
${\cal S}^{\alpha, k}(\R^n)$ the set of complex valued smooth functions on
$\R^n$ which can be written $f(\xi)=\sum_{l=0}^k f_l(\xi)\log ^l|\xi|$ where the $f_l$'s
have the following asymptotic
behaviour as $\vert\xi \vert \to \infty$:
\begin{equation}\label{eq:asymptf}
f_l(\xi)\sim\sum_{j=0}^\infty f_{\alpha-j, l}(\xi).
\end{equation}
Here $f_{\alpha-j, l}(\xi)$ is positively homogeneous of degree
$\alpha-j$. The notation $\sim$ stands for the existence, for any positive
integer $N$ and any multi-index $\gamma$, of  a positive constant $C_{\gamma,N}$ such that 
$$\left\vert \partial^\gamma\left(f_l(\xi)-\sum_{j=0}^{N}\,  f_{\alpha-j,
      l}(\xi)\right)\right\vert \le C_{\gamma,N}(1+\vert \xi\vert)^{\smop{Re }\alpha-N-1-\vert\gamma\vert},$$
with $\partial^\gamma=(\frac{\partial}{\partial \xi_1})^{\gamma_1}\cdots(\frac{\partial}{\partial \xi_k})^{\gamma_k}$, with $|\gamma|=\gamma_1+\cdots +\gamma_k$ and where $\mop{Re }\alpha$ stands for the real part of $\alpha$. We simplify  notations setting
${\cal S}^{\alpha, k}:={\cal S}^{\alpha, k}(\R^n)$ when there is no ambiguity
on the dimension, which is one in most of the paper, except for Section 5
where it is $n>1$. We set:
$$ {\cal S}^{\alpha, *}:= \bigcup_{k=0}^\infty {\cal S}^{\alpha, k},$$
and ${\cal S}^{*, k}$ stands for the linear span over $\C$ of all ${\cal S}^{\alpha,k}$ for
${\alpha\in \C}$. We also define:
$${\cal S}^{*, *}:= 
\bigcup_{k=0}^\infty {\cal S}^{*, k}$$
which is a filtered algebra (with respect to the logarithmic power $k$) for
the ordinary product of functions. If $f\in{\cal S}^{\alpha,k}$ and if there
is an $l\in\N$ such that $f_{\alpha,l}\not =0$ in (\ref{eq:asymptf}) we call $\alpha$
the {\sl order\/} of $f$. Note that the union
$\bigcup_{\alpha\in -\N}\bigcup_{k=0}^\infty {\cal S}^{\alpha, k}$ is a
  subalgebra of ${\cal S}^{*,*}$, and that $\bigcup_{\alpha\in -\N} {\cal S}^{\alpha, 0}$ is a
  subalgebra of ${\cal S}^{*,0}$. The {\sl real order\/} of $f\in{\cal
  S}^{*,*}$ will be defined as:
$$o(f)=\mop{inf}\{\lambda\in\R,\, \mopl{lim}_{|\xi|\to +\infty}|\xi|^{-\lambda}f(\xi)=0\},$$
hence $o(f)=\mop{Re }\alpha$ for $f\in{\cal S}^{\alpha,*}$ of order
$\alpha$.\\

The
algebra ${\cal S}^{*,*}$ (resp. the subalgebra ${\cal S}^{*,0}$) coincides with the algebra of
{\sl log-polyhomogeneous\/} (resp. {\sl classical\/}) constant coefficient symbols on
$\R^n$, which is often denoted by $CS^{*,*}$ (resp. $CS$) in the
pseudodifferential literature. Typical examples of classical symbols are the
polynomials in $n$ variables, and more generally complex powers of
polynomials, which can be defined for any polynomial provided it nowhere
vanishes on $\R^n$ and takes values outside an angular sector pointed at the origin of the complex plane. Log-polyhomogeneous symbols naturally appear by taking
iterated primitives of suitable classical ones (\cite{L}, Section 3).\\

We recall  results which belong to folklore knowledge, namely that the
definition of ${\cal
  S}^{*, *}$  is independent of the choice of the norm. This is due to the
fact that all norms are equivalent on a finite dimensional space: Let
$\vert\cdot \vert$ and  $\vert\cdot \vert_1$ be two norms smooth outside the origin. Let  $f \in
{\cal S}^{\alpha, k}$ which we write
\begin{equation}\label{sigmadev}
f(\xi)= \sum_{l=0}^k\sum_{j=0}^{K_N} f_{\alpha -j,l}(\xi)\log^l\vert
\xi \vert+ f_{(N)}(\xi)
\end{equation}
with $f_{\alpha-j,l}$ positively homogeneous of order $\alpha-j$ and $ f_{(N)}$ of order $\leq -N$. Then 
\begin{eqnarray*}
f(\xi)&=&\sum_{l=0}^k\sum_{j=0}^{K_N} f_{\alpha-j,l}(\xi)\log^l\vert \xi \vert+ f_{(N)}(\xi)\\
&=& \sum_{l=0}^k\sum_{j=0}^{K_N} f_{\alpha-j,l}(\xi)\left(\log\vert 
\xi \vert_1+ \log \frac{\vert \xi \vert}{{\vert \xi\vert}_1}\right)^l+ f_{(N)}(\xi)\\
&=& \sum_{l=0}^k\sum_{j=0}^{K_N} f^1_{\alpha-j,l}(\xi)\,  \log^l\vert 
\xi \vert_1 + f_{(N)}(\xi)
\end{eqnarray*}
where we have set 
$$f^1_{\alpha-j,i}(\xi):= \sum_{l=0}^kf_{\alpha-j,l}(\xi){l\choose i}\, 
\log^{l-i} \frac{\vert \xi \vert}{{\vert \xi\vert}_1}.$$
It follows immediately that $f^1_{\alpha-j,i}$ is positively homogeneous of
order $\alpha-j$, which proves the claim.
\subsubsection{Cut-off integrals revisited}
Cut-off integrals can be defined on the class ${\cal S}^{*, *}$ of
(constant coefficient) log-polyhomo\-geneous symbols \cite{L}. The
result depends on the choice of the norm. We recall however the asymptotic
behaviour of integrals on balls, which can be derived independently of the
specific norm chosen.
\begin{prop}\label{prop:extRotaBaxter} Let $d_S\xi$ be the volume measure on the unit sphere
  induced by the canonical Lebesgue measure on $\R^n$. Let $B(0, R)$
  be the  closed ball in $\R^n$ centered at $0$ with radius
  $R$ and for any $0\leq r \leq R$, let us set $B(r, R)=
 B(0, R)- B(0, r)$. Then  for any fixed $r$ and for any $f \in {\cal S}^{*,
   k}$ displayed as in \eqref{sigmadev}, the integral 
$\int_{B(r, R)}f(\xi) \, d\xi$ has an asymptotic expansion as $R\to +\infty$ of the type
\begin{eqnarray}\label{eq:cutoffasympt}
\int_{B(r, R)} f(\xi) d\xi&\sim_{R\to\infty} &C_{r}(f)+ \sum_{j=0,\alpha-j+n\neq 0}^\infty \sum_{l=0}^k P_{l}(f_{\alpha-j, l})(\log R)  \, R^{\smop{Re}\alpha-j+n}\nonumber\\
&+&
\sum_{l=0}^k r_{l}(f)\log^{l+1}  R
\end{eqnarray}
where the $r_{l}(f)$ are positive constants depending on  $f_{l,-n}$,
$P_{l}(f_{\alpha-j, l})(X)$ is a polynomial of degree
$l$ with coefficients depending on $f_{\alpha-j, l}$ and where the constant term
$C_{r}(f)$ corresponds to the finite part:
\begin{eqnarray*} 
C_{r}(f) &
:=&\int_{\R^n} f_{(N)}(\xi)\, d\xi+  \int_{B(r, 1)}
 f(\xi)\, d\xi\nonumber\\
 &+&\sum_{j=0, \alpha-j+n\neq 0}^{N}\sum_{l=0}^k \frac{(-1)^{l+1}l!}{(\alpha-j+n)^{l+1}}
 \int_{|\xi|=1} f_{\alpha-j,l} (\xi) d_S\xi\\
\end{eqnarray*}
which is independent of $N$ for  $N\geq \alpha +n-1$.
\end{prop}
 {\bf Proof:}
A proof can be
found e.g. in \cite {L} (see also \cite{MP}) and references therein for the case $r=0$. The case $r>0$ then easily follows subtracting $\int_{B(0, r)}f(\xi) \, d\xi$ from the constant term so that $C_{r}(f)= C_{0}(f)- \int_{B(0, r) }f(\xi) \, d\xi.$
\endsquare\\ \\
The {\sl cut-off integral\/} of $f$ outside the ball of radius $r$ is defined
as the finite part $C_r(f)$, and is  denoted by:
\begin{equation}\label{eq:finitepart}
\cutoffint_{r\leq \vert \xi \vert }f(\xi) \, d\xi.
\end{equation}
It coincides with the ordinary integral when the latter converges, but may
depend on the choice of the norm in the general case.
\subsubsection{Holomorphic families} 
 Let us first recall what we mean by a holomorphic family of elements of some
 topological vector space ${\cal A}$. Let $W\subset\C$ be a complex domain. A
 family $\{f(z)\}_{z\in W}\subset {\cal A}$   is holomorphic at
$z_0\in W$ if
the corresponding function $f:W\to{\cal A}$ admits a  Taylor expansion in a neighbourhood $N_{z_0}$ of $z_0$
\begin{equation}\label{e:Texpansion}
f(z) = \sum_{k=0}^{\infty}f^{(k)}(z_0)\,\frac{(z-z_0)^k}{k!}
\end{equation}
which is convergent, uniformly on compact subsets of $N_{z_0}$,
with respect to the  topology on ${\cal A}$. The vector spaces of functions we consider here are $C(\R^n,\C)$ and
$C^\infty(\R^n,\C)$ equipped with their usual topologies, namely uniform
convergence on compact subsets, and uniform convergence of all derivatives on
compact subsets respectively. 
\begin{defn}\label{defn:holfamilies} Let $k$ be a non-negative integer, and
  let $W$ be a domain in $\C$. A {\rm simple holomorphic
family of log-polyhomogeneous symbols\/} $f(z)\in {\cal S}^{*, k}$ parametrised by $W$ of order
$\alpha:W\to\C$   means a holomorphic family 
$f(z)(\xi) := f(z,\xi)$  of smooth functions on $\R^n$ for which: 
\begin{enumerate}
\item $f(z)(\xi)=\sum_{l=0}^k f_l(\xi)\, \log^l \vert \xi\vert$
with $$
 f_l(z)(\xi) \sim \sum_{j\geq 0}\, 
 f(z)_{\alpha(z)-j,l}(\xi).$$
Here $\alpha: W\to \C$ is a holomorphic map and  $f(z)_{\alpha(z)-j,l}$ is positively homogeneous of degree $\alpha(z)-j$.
\item For any positive integer $N$ there is some positive integer $K_N$ such
  that the remainder term 
$$f_{(N)}(z)(\xi):= f(z)(\xi) -\sum_{l=0}^k \sum_{j= 0}^{K_N}
 f(z)_{\alpha(z)-j, l}(\xi)\log^l\vert\xi\vert=o(\vert \xi\vert^{-N})$$
is holomorphic in $z\in W$ as a function of $\xi$ and verifies for any
$\epsilon>0$ the following estimates:
\begin{equation}\label{e:kthderivlogclassical}
\partial_\xi^\beta\partial_z^kf_{(N)}(z)(\xi)=o(|\xi|^{-N-|\beta|+\e})
\end{equation}
for $k\in \N$ and $\beta\in\N^n$.
\end{enumerate}
A {\rm holomorphic
family of log-polyhomogeneous symbols\/} is a finite linear combination (over
$\C$) of simple holomorphic families.
\end{defn}
Note that the notion of holomorphic order $\alpha(z)$ only makes sense for
simple holomorphic families. Cut-off integrals of log-polyhomogeneous
simple holomorphic families behave as follows:
\cite{L}:
\begin{prop} \label{prop:cutoffinthol} Given a  non-negative integer $k$, let $\{f(z)\}$ be a  family
  in ${\cal S}^{*,k}$  holomorphic on $W\subset \C$ and simple. Let us assume that the
  order $\alpha(z)$ of $f(z)$ is of the form $\alpha(z)=-q\, z+\alpha(0)$ with
  $q\neq 0$. Then the map $$z\mapsto \cutoffint_{\R^n} f(z)(\xi) $$ is meromorphic with
  poles of order $\leq k+1$ in the discrete set:
$$\alpha^{-1}\left([-n, +\infty[\cap \, \Z\right)=\{\frac{\alpha(0)+n-j}{q},\ j\in\N\}.$$ \end{prop}
\subsection{Nested integrals on classical symbols}
Let us introduce  the following nested
 integrals of symbols on compact sets:
\begin{defn}
Given $R>r>0$, 
\begin{eqnarray*}
&{}&\int_{B(r, R)}^{{\rm Chen}}\sigma_1\otimes \cdots \otimes \sigma_k\\
&:= &\int_{r\leq \vert \xi_k\vert \leq \vert
  \xi_{k-1}\vert \leq \cdots\leq \vert \xi_1\vert\leq R }d\xi_{1}\cdots d\xi_k\, 
\sigma_{1}(\xi_1)\, \sigma_{2}(\xi_2)\cdots \sigma_{k} (\xi_k)\\
&=& \int_{r\leq  r_k \leq 
  r_{k-1} \leq \cdots \leq r_1\leq R } \, dr_{1}\cdots  \,
  dr_k\, f_1(r_1)\cdots f_k(r_k)  \\
\end{eqnarray*}
where we have set $f_i(r):= r^{n-1} \int_{|\xi|=1}\sigma_i(r\xi) \, d_S\xi$.  
\end{defn}
 These nested integrals correspond to ordinary nested (or iterated) integrals
 $$\int_{r\leq r_k\leq  \cdots\leq r_1\leq R} \omega_1\wedge\cdots\wedge\omega_k,$$ 
with  $\omega_i(t)=
f_i(t) dt$. As such they enjoy the usual properties of one-dimensional nested
integrals (see e.g  \cite{Ch}, or Appendix XIX.11 in \cite{Ka}): 
\begin{defn}
The tensor algebra   ${\cal T}\left ({\cal S}^{*,*}\right):=
\bigoplus_{k=0}^\infty ({\cal S}^{*,*})^{\otimes k}$ (over the base field
$\C$) is equipped with a shuffle product: 
$$\left(\sigma_1\otimes \cdots \otimes \sigma_k\right)\shu\left( \sigma_{k+1}\otimes \cdots
\otimes \sigma_{k+l}\right):= \sum_{\tau\in \Sigma_{k; l}} \sigma_{\tau^{-1}(1)}\otimes\cdots
\otimes \sigma_{\tau^{-1}(k+l)} $$
where $\tau $ runs over the set $\Sigma_{k; l}$ of $(k,l)$-shuffles, i.e. permutations $\tau$ on
$\{1, \ldots, k+l\}$ satisfying $\tau(1)< \cdots < \tau(k)$  and
$\tau(k+1)< \cdots < \tau(k+l)$. 
\end{defn}
The shuffle product and the deconcatenation coproduct:
$$\Delta \left(\sigma_1\otimes \cdots \otimes \sigma_k\right):= \sum_{j=0}^k
\left(\sigma_1\otimes\cdots \otimes \sigma_j\right) \bigotimes\left( 
\sigma_{j+1}\otimes \cdots \otimes \sigma_k\right)$$
endow ${\cal T}\left ({\cal S}^{*,*}\right)$ with a structure of connected
graded commutative Hopf algebra. This algebraic construction holds if we
replace ${\cal S}^{*,*}$ by any vector space \cite{H2}. As in the case of ordinary nested integrals, we have the following shuffle relations:
\begin{prop}\label{prop:starChen}
 For any  $\sigma, \tau \in {\cal T}\left(
  {\cal S}^{*,*}\right)$ 
\begin{equation}\label{eq:shuffle}
\int^{{\rm Chen}}_{B(r, R)}  \sigma\shu \tau
 =\int^{{\rm Chen}}_{B(r, R)}\sigma \, \int^{{\rm Chen}}_{B(r,
  R)}\tau.
\end{equation}
 \end{prop}
{\bf Proof:} We may assume that $\sigma$ and $\tau$ are indecomposable
elements. Setting $\sigma= \sigma_1\otimes \cdots \otimes\sigma_k$ and
$\tau=\tau_1\otimes \cdots \otimes \tau_l$, $f_i(r)= r^{n-1} \int_{|\xi|=1}
\sigma_i (r \xi) \, d_S\xi$ and $g_j(r)= r^{n-1} \int_{|\xi|=1}
\tau_j (r \xi) \, d_S\xi$ we have:
\begin{eqnarray*}
\int^{{\rm Chen}}_{B(r, R)}  \sigma\shu \tau
&=&
\int^{{\rm Chen}}_{[ r, R]} \left( f_1\otimes \cdots\otimes f_k\right)\, \shu  \left(\, g_1\otimes \cdots
\otimes g_l\right)\\
&=&\left(\int^{{\rm Chen}}_{[r, R]}f_1\otimes \cdots \otimes f_k \right) \,\left( \int^{{\rm Chen}}_{[r, R]}
  g_1\otimes \cdots\otimes g_l\right)\\
&=& \int^{{\rm Chen}}_{B(r, R)}\sigma \, \int^{{\rm Chen}}_{B(r,
  R)}\tau.
\end{eqnarray*}
where we have used the shuffle property for usual nested Chen integrals
between the bounds $r$ and $R$.
\endsquare\\ \\
Finally, let us recall Chen's lemma:
\begin{prop} Given $\sigma=\otimes_{i=1}^k\sigma_i \in
\otimes_{i=1}^k {\cal S}^{*,*}$ we have:
\begin{eqnarray}\label{chenlemma}
\int_{B(r, R)}^{{\rm Chen}}\sigma &=& \int_{B(r, \Lambda)}^{{\rm Chen}}\sigma+
\int_{B( \Lambda, R)}^{{\rm Chen}}\sigma\nonumber \\
&+&\sum_{j=1, \ldots, k-1} \int_{B(r,
  \Lambda)}^{{\rm Chen}}\sigma_{i_1}\otimes \cdots\otimes\sigma_{i_j}\,
\int_{B(\Lambda, R)}^{{\rm Chen}}\sigma_{i_{j+1}}\otimes \cdots
\otimes\sigma_{i_k}.
\end{eqnarray}
\end{prop}
To summarise these two facts, the three linear maps $\int_{B(r, \Lambda)}^{{\rm Chen}}$,
$\int_{B(\Lambda, R)}^{{\rm Chen}}$ and  $\int_{B(r, R)}^{{\rm Chen}}$ are
characters of the Hopf algebra ${\cal T}({\cal S}^{*,*})$, and moreover
satisfy the relation:
\begin{equation}
\int_{B(r, R)}^{{\rm Chen}}=\int_{B(r, \Lambda)}^{{\rm Chen}}\star\int_{B(\Lambda, R)}^{{\rm Chen}},
\end{equation}
where $\star$ stands for the convolution product.
\subsection{Cut-off  nested  integrals of symbols}\label{sect:cutoffchenint}
In this paragraph we let $R$ tend to infinity, which requires a
regularisation since the integrals do not a priori converge when $R\to
\infty$.  Let us first recall how nested integrations can turn classical symbols into
log-polyhomogeneous symbols.
\begin{lem}\label{lem:contRotaBaxter} \cite{MP} Let $r\ge 0$. The following operator on
  $C^\infty(\R^n)$:
\begin{equation}\label{eq:contRotaBaxter}
\tilde P_{r}(f)(\eta):= 
\int_{r\leq \vert \xi \vert\leq \vert  \eta\vert} f(\xi) \, d\xi
\end{equation}
 maps ${\cal S}^{*, k-1}$ to ${\cal S}^{*, k}$ for any positive integer $k$.
\end{lem} 
{\bf Proof:}
It follows from setting $R= \vert \eta\vert$ in equation  (\ref{eq:cutoffasympt}).
The smoothness of the norm outside the origin  ensures that whenever $f$ is a
symbol, then  the resulting
map $\tilde P_r(f)$  is also a symbol.
\endsquare\\ \\
As a result of Proposition \ref{prop:starChen} the map 
$\tilde P_{r}$  satisfies the following Rota-Baxter relation  (see \cite{MP} Theorem 3):
\begin{equation} \label{eq:contRotaBaxterequation}
\tilde P_{r}(f)\tilde P_{r}(g)=
\tilde P_{r}\left(f\,\tilde  P_{r}(g)\right)+
 \tilde P_{r}\left(\tilde
  P_{ r}(f)\, g\right).
\end{equation}
Iterating Lemma \ref{lem:contRotaBaxter} we get:
\begin{lem}\label{lem:CStoCSlog} 
Given $ \sigma_i\in {\cal S}^{\alpha_i,0}$ and 
 $\sigma:=\otimes_{i=1}^k\sigma_i  $, for any $r\geq 0$ and $k\ge 2$ the function
\begin{equation}
\sigma_r^{\smop{Chen}}:\xi \mapsto  \sigma_1(\xi )\int^{{\rm
    Chen}}_{B(r,\vert\xi \vert)}
 \sigma_2\otimes\cdots\otimes \sigma_{k}
\end{equation}
lies  in ${\cal S}^{*, k-1}$ as a linear
combination of log-polyhomogeneous  functions of order
$\alpha_{1}+\alpha_{2}+\cdots +\alpha_{j}+(j-1)n$, $j=1, \ldots, k$. It has
real order $\omega_1$, where we successively define $\omega_k=\mop{Re }
  \alpha_k$ and $\omega_{j}=\mop{Re }\alpha_{j}+\mop{max}(0,\,\omega_{j+1}+n)$.
\end{lem}
  We can define the cut-off nested integral of $\sigma$ as the cut-off
  integral of the symbol $\sigma_r^{\smop{Chen}}$, namely:  
\begin{defn} Given $\sigma_i \in {\cal S}^{\alpha_i,0}, i=1, \ldots, k$, setting 
$\sigma=\otimes_{i=1}^k\sigma_i$, for any $r>0$  we set 
\begin{equation}
\cutoffint_{r}^{{\rm Chen}}\sigma:= 
{\rm fp}_{R\to \infty} \int^{{\rm Chen}}_{B(r,R)} \sigma.
 \end{equation}
 \end{defn}
\begin{rk}\label{rk:convChen}
 It follows from the above lemma that  if
  $\sum_{j=1}^m\mop{Re }\alpha_j< -nm$ for any $m\in\{1,\ldots ,k\}$, then the symbol
  $\sigma^{{\rm Chen}}_r$ also has real order $<-n$ so that the nested integral
  converges and the cut-off integrals become ordinary integrals.
\end{rk}
We now want to extend  the shuffle relation (\ref{eq:shuffle}) in the $R\to \infty$ limit. 
\begin{thm}\label{thm:cutoffshuffle}
Given
$\sigma\in\bigotimes_{i=1}^k{\cal S}^{\alpha_i,0}$ and 
$\tau\in\bigotimes_{j=1}^l{\cal
  S}^{\beta_j,0}$
then,
\begin{enumerate}
\item  if all the partial sums
$\alpha_1+\alpha_{2}+\cdots +\alpha_{m}+ \beta_1+ \beta_{2}+\cdots
+\beta_{p}$ with  $m\in\{1,
  \dots, k\}$, $p\in\{1, \ldots, l\}$ are non-integer valued, the following shuffle
relation holds:
$$ \cutoffint_{r}^{{\rm Chen}}
\sigma\, \shu  \, \tau
  =\left(\cutoffint_{r}^{{\rm Chen}}\sigma\right) \, \left(\cutoffint_{r}^{{\rm Chen}}\tau\right).$$
\item If  one of the two conditions: $\mop{Re }(\alpha_1+\cdots +\alpha_m)<
  -nm$ for any $m\in\{1,\ldots ,k\}$, or $\mop{Re }(\beta_1+\cdots +\beta_m)<
  -nm$ for any $m\in\{1,\ldots ,l\}$   holds,  the shuffle relation
  also holds.
\item If $\mop{Re }(\alpha_1+\cdots +\alpha_m)<
  -nm$ for any $m\in\{1,\ldots ,k\}$ {\bf and} $\mop{Re }(\beta_1+\cdots +\beta_m)<
  -nm$ for any $m\in\{1,\ldots ,l\}$  both hold,  the shuffle relation
  moreover becomes an
  equality of ordinary nested integrals.
 \end{enumerate}
\end{thm}
 {\bf Proof:} This follows from taking  finite parts when $R\to\infty$ in   Proposition
 \ref{prop:starChen}.  
\begin{enumerate}
\item  Let us now assume that  all the partial sums
$\alpha_1+\alpha_{2}+\cdots +\alpha_{m}+ \beta_1+ \beta_{2}+\cdots
+\beta_{p}$ with  $m\in{1,
  \dots, k}$, $p\in\{1, \ldots, l\}$ are non-integer valued. Under this assumption,  the finite part of the product of the integrals coincides with
the product of the finite parts of the integrals. Let us make this statement
more precise: On the right hand
side of 
$$\int_{B(r, R)}^{{\rm Chen}} \sigma\shu \tau
  =\left(\int_{B(r, R)}^{{\rm Chen}}\sigma \right)\, \left(\int_{B(r,
      R)}^{{\rm Chen}}\tau\right),$$ one might expect cancellations of divergences when $R\to \infty$
  which could lead to new contributions to the finite part.
Since the asymptotic expansion of 
$$\int_{B(r, R)}^{{\rm Chen}} \sigma=
\int_{r\leq \vert \xi\vert \leq  R}   \sigma_r^{\smop{Chen}}(\xi)
\,d\xi$$ when $R\to \infty$ involves
powers $R^{\alpha_1+\alpha_{2}+\cdots +\alpha_{m} +m\,n  -j}$ for non-negative integers $j$  combined with  powers
of $\log R$ (and similarly for $\int_{B(r, R)}^{{\rm Chen}} \tau$
), such
cancellations  of
divergences do not occur  if $\alpha_1+\alpha_{2}+\cdots +\alpha_{m}+ \beta_1+ \beta_{2}+\cdots
+\beta_{p}$  do not lie in $[-n, \infty[\cap \Z$.

\item If one of the conditions of item 2 holds, then one of the  expressions  $\int_{B(r, R)}^{{\rm Chen}}\sigma $ or $ \int_{B(r, R)}^{{\rm Chen}}\tau$  
converges in the limit $R\to \infty$ so that   such cancellations do not occur either.  As a result, 
 the shuffle property also holds in the limit.
\item If both conditions of item 2 hold, all the integrals converge and the finite parts are ordinary limits. Taking the limits when $R\to \infty$ in   Proposition
 \ref{prop:starChen} yields the result. 
\end{enumerate}
\endsquare
\subsection{A continuous analog of multiple zeta functions}\label{sect:contmz}
We define a continuous  analog of multiple zeta functions via nested integrals of
symbols. For $\e>0$, let $\chi_\e$ be a smooth cut-off function on $\R^n$ that
vanishes in  the ball $B(0, \frac{\e}{4})$ and is constant equal to $1$
outside the ball $B(0, \frac{\e}{2})$. The symbols
$$\sigma_i(\xi) =\chi_\e(\xi)|\xi|^{-s_i}$$
lie in  ${\cal S}^{*,0}$. 
For $\underline s:= (s_1, \ldots, s_k)\in \C^k$ we set $\sigma_{\underline
  s}:= \sigma_{s_1}\otimes \cdots \otimes \sigma_{s_k}$ and define for $r\geq
\frac{\e}{2}$,  
\begin{eqnarray*}
\tilde \zeta^{r }_{n}(\underline s)&:=& \cutoffint_{r\leq \vert \xi_k\vert \leq
  \cdots \leq \vert \xi_1\vert}\sigma_{\underline s}\\
&=& \cutoffint_{r\leq \vert \xi_k\vert \leq
  \cdots \leq \vert \xi_1\vert}
d\xi_1\cdots d\xi_k\, \prod_{i=1}^k |\xi_i|^{-s_i}.
\end{eqnarray*}
If $s_1> n$ and $s_1+\cdots +s_k> nk$ then, by  Remark
  \ref{rk:convChen}, $\tilde \zeta_{n}^{r}(\underline s)$
  converges and the regularised integrals reduce to ordinary ones. In dimension one we write:
$$\tilde \zeta^r (\underline s)
:=\cutoffint_{r\leq  t_{k}\leq t_{k-1} \leq \cdots\leq t_1} dt_1\cdots  dt_k\, \prod_{i=1}^k  \, t_i^{-s_i}.
$$
This corresponds to $2^{-k}\tilde \zeta_{1}^{r}(\underline s)$ the norm chosen being the usual absolute value. This clearly does not extend to $r=0$. 
\begin{rk}
Note the analogy with the usual multiple zeta functions, an analogy  which was made precise in \cite{MP}  via an Euler-MacLaurin formula.
\end{rk}
The following proposition gives an  analog of the ``second shuffle relations''
\cite{ENR}, also called stuffle relations for multiple zeta functions (but without ``diagonal
  terms'', which arise only in the discrete setting. As such they appear as
  ``ordinary'' shuffle relations). 
\begin{prop}Let $\underline s=(s_1,\ldots,s_k)\in \R^k, \underline s^\prime=(s_{k+1},\ldots ,s_{k+l})\in \R^{l}$
such that $s_1+\cdots +s_m>nm$ for any $m\in\{1,\ldots ,k\}$ and $s_{k+1}+\cdots +s_{k+m}>nm$ for any $m\in\{1,\ldots ,k\}$. For any $r\geq 0$
$$ \tilde \zeta_n^{r }(\underline s)\, \tilde  \zeta_n^r
  (\underline s^\prime)=\tilde  \zeta_n^r(\underline s\shu \underline s^\prime).$$
where $\underline s\shu \underline s^\prime$ is defined as usual by:
 \begin{equation}
(s_1,\ldots,s_k)\shu(s_{k+1},\ldots ,s_{k+l})=\sum_{\tau\in\Sigma_{k;l}}(s_{\tau^{-1}(1)},\ldots,s_{\tau^{-1}(k+l)}).
\end{equation}
\end{prop}
{\bf Proof:} It follows from applying the first item of   Theorem \ref{thm:cutoffshuffle} to 
$\sigma=\sigma_{\underline s}$ and $\tau=\sigma_{\underline s^\prime}$
 since 
$$\sigma_{\underline s}\shu \sigma_{\underline s^\prime}= \sigma_{\underline s\sshu \underline s^ \prime}.$$
\endsquare\\ \\
When $n=1$, $\e<r$ and $\vert\cdot\vert$  is  the usual absolute value, this yields:
\begin{equation}\label{eq:tildezetashuffle1}
 \tilde \zeta^{r}(\underline s)\, \tilde  \zeta^{r}(\underline s^\prime)=\tilde
 \zeta^{r}(\underline s\shu \underline s^\prime).
\end{equation}
We also describe a continuous analog of
polylogarithms in the Appendix.
\subsection{Nested  integrals of  holomorphic families of symbols}\label{sect:regul}
Let us briefly recall the notion of holomorphic regularisation 
inspired by \cite{KV}.
\begin{defn}
A  holomorphic
regularisation procedure on ${\cal S}^{*, *}$ 
is a  map 
\begin{eqnarray*}
{\cal R}: {\cal S}^{*, *}&\to & {\rm Hol}_W\,\left( {\cal S}^{*, *}\right)\\
f  &\mapsto & \{f(z)\}_{z\in W}
\end{eqnarray*}
where $W$ is an open subset of $\C$ containing $0$, and ${\rm Hol}_W\left( {\cal S}^{*, *}\right)$ is the algebra of
holomorphic families in ${\cal S}^{*, *}$ as defined in Section 1, 
such that for any $f\in {\cal S}^{*, *}$,
\begin{enumerate}
\item  $f(0)=f$, 
\item  the holomorphic family $f(z)$ can be written as a linear combination of
  simple ones:
$$
f(z)=\sum_{j=1}^kf_j(z),
$$
the holomorphic  order  $\alpha_j(z)$ of which verifying $\mop{Re
  }\alpha_j^\prime(z)< 0$ for any $z\in W$ and any $j\in\{1,\ldots ,k\}$.\\
\end{enumerate}
The holomorphic regularisation ${\cal R}$ is {\rm simple} if for any
  log-polyhomogeneous symbol $f\in{\cal S}^{\alpha,k}$ the holomorphic family
  ${\cal R}(f)$
  is simple.
\end{defn}
A similar definition holds with suitable subalgebras of ${\cal S}^{*,*}$,
e.g. classical symbols ${\cal S}^{*,0}$ instead of
log-polyhomogeneous ones. Simple holomorphic regularisation procedures naturally arise in physics:
\begin{ex}\label{ex:regphys}
 Let $z\mapsto \tau(z)  \in {\cal S}^{*,0}$ be a holomorphic family of
  classical symbols such that $\tau(0)=1$ and $\tau(z)$ has holomorphic  order
 $t(z)$ with $\mop{Re }t^\prime(z)< 0$. Then 
$${\cal R}: \sigma \mapsto \sigma(z):= \sigma\tau(z)$$
yields a holomorphic regularisation on ${\cal S}^{*, *}$ as well as on ${\cal S}^{*, 0}$. Choosing $\tau(z)(\xi):= \chi(\xi)+\big(1-\chi(\xi)\big)\big(H(z)\, \vert\xi \vert^{-z}\big)$
  where $H$ is a scalar valued  holomorphic map such that  $H(0)=1$, and where
  $\chi$ is a cut-off function as
 defined in the introduction, we get
\begin{equation}\label{eq:dimreg}
{\cal R}(\sigma)(z)(\xi)= \chi(\xi)\sigma(\xi)+\big(1-\chi(\xi)\big)\big(H(z)\, \sigma\, \vert \xi \vert^{-z}\big).
\end{equation}
Dimensional regularisation commonly used in physics is of this type, where $H$
is expressed in terms of Gamma functions which account for a
``complexified'' volume of the unit sphere. When $H\equiv 1$, such a
regularisation   ${\cal
  R}$ is called   Riesz regularisation.
  \end{ex}
From a given simple holomorphic regularisation procedure ${\cal R}:\sigma \mapsto
\sigma(z)$ on ${\cal S}^{*,0}$ which yields ordinary regularised integrals
$$\cutoffint_{r}^{{\cal R}}\sigma:= {\rm
  fp}_{z=0}\cutoffint_{r\leq \vert \xi \vert}  \sigma(z)(\xi)\, d\xi, \quad \forall
r\geq 0,$$ 
we  wish
to build regularised nested integrals
$$\cutoffint_{r }^{
{\cal R},{\rm Chen}} \sigma_1\otimes \cdots \otimes \sigma_k:=\cutoffint^{
{\cal R}}_{r\leq \vert \xi_k \vert\leq \cdots \leq \vert
  \xi_1\vert}\sigma_1(\xi_1)\cdots  \sigma_k(\xi_k)\, d\xi_1\cdots d\xi_k $$ with the following

\begin{enumerate}
\item For $k=1$, the nested integral coincides with the ordinary regularised integral $\cutoffint^{{\cal
 R}}_{r} .$
\item It coincides with the ordinary   Chen
  nested integral 
$$\int_{r }^{{\rm Chen}} \sigma_1\otimes \cdots \otimes \sigma_k:=\int_{r\leq \vert \xi_k \vert\leq \cdots \leq \vert
  \xi_1\vert}\sigma_1(\xi_1)\cdots  \sigma_k(\xi_k)\, d\xi_1\cdots d\xi_k$$
on symbols $\sigma_1, \ldots, \sigma_k$ with respective
    orders $ \alpha_1,\ldots,\alpha_k$ such that $ \mop{Re }(\alpha_1+\cdots+\alpha_m)<-nm$ for any $m\in\{1,\ldots ,k\}$.
\end{enumerate}
Given a simple holomorphic regularisation ${\cal R}:\sigma \mapsto
\sigma(z)$ on ${\cal S}^{*,0}$, we can  assign to  $\sigma=\sigma_1\otimes
\cdots\otimes\sigma_k$ holomorphic  maps
\begin{equation}\label{eq:sigmatilde} z\mapsto \tilde{\cal R}( \sigma)(z):= \sigma_1(z)\otimes
 \cdots \otimes\sigma_k(z).\end{equation}
Clearly, for any tensor products $\sigma=\sigma_1\otimes
\cdots\otimes\sigma_k$ and $\tau=\tau_1\otimes
\cdots\otimes\tau_l$ we have
$\tilde{\cal R}( \sigma\otimes \tau)= \tilde{\cal R}( \sigma)\otimes \tilde{\cal R}(
\tau)$,
 so that the regularisation is compatible with the shuffle product: 
\begin{equation}\label{eq:regshu}\tilde{\cal R}( \sigma\shu \tau)= \tilde{\cal R}( \sigma)\shu\tilde{\cal R}(
\tau).
\end{equation}
\begin{prop}\label{prop:Chenmero}
Let ${\cal R}: \sigma\mapsto \{\sigma(z)\}_{z\in W}$ be a holomorphic regularisation
procedure defined on a neighbourhood $W$ of $0$ in the complex plane. Given any  non
negative real number $r$, for any $\sigma_i\in {\cal S}^{\alpha_i,0}$, $i=1,
\ldots, k$, the  map $$z\mapsto  \cutoffint_{r}^{{\rm Chen}}
  \sigma_1(z)\otimes \cdots\otimes\sigma_k(z) $$
  is  meromorphic with poles 
 of order at most $k$. If moreover the holomorphic regularisation is simple
 and if for any $i\in\{1,\ldots ,k\}$ the symbol $\sigma_i(z)$ has order 
$\alpha_i(z)=-qz+\alpha_i(0)$ affine in $z$ with
$q>0$, the poles lie in  the discrete set of points $W\cap D_k$, with:
$$ D_k:= \bigcup_{j=1}^k \left\{\frac{\sum_{i=1}^j
    \alpha_i(0)+jn-\gamma}{jq},\ \gamma\in\N\right\}.$$
\end{prop}
{\bf Proof:}  This result extends 
 Theorem 5 in \cite{MP}  in the case $r=0$. The first assertion follows from    $$\cutoffint_{r}^{{\rm Chen}}
  \sigma(z)= \cutoffint_{|\xi|\ge r}\left( 
  \sigma(z)\right)^{{\rm Chen}}(\xi)\, d\xi.$$
The second assertion comes from the
    fact that by Lemma \ref{lem:CStoCSlog}, 
     $\left( \sigma(z)\right)^{{\rm Chen}}$ lies in
${\cal S}^{*,k-1}$ as a linear combination of functions in ${\cal S}^{*, j-1}$
of order 
 $\alpha_{1}(z)+\alpha_{2}(z)+\cdots +\alpha_{j}(z)+ (j-1)n$, $j=1, \ldots,
 k$. Applying Proposition \ref{prop:cutoffinthol}, we infer that the map
 $z\mapsto \cutoffint_{|\xi|\ge r}\left( 
  \sigma(z)\right)^{{\rm Chen}}(\xi)\, d\xi$ has  poles at $z_{j,\gamma}:=
 \frac{\sum_{i=1}^j \alpha_i(0)+jn-\gamma}{jq},  \gamma\in \N$. A closer look shows that the symbols of order
  $\alpha_1(z)+\cdots
 +\alpha_j(z)+(j-1)n$ arise from the $j$-th outermost integration so that poles at $z_{j,\gamma}$ are of order $j$. 
\endsquare
\subsection{Shuffle relations for renormalised nested integrals of symbols}
\begin{prop}\label{prop:Chencharacter}
Given an open neighbourhood $W$ of $0$ in $\C$ and a holomorphic regularisation ${\cal R}: \tau \mapsto \{\tau(z)\}_{z\in W}$ on
${\cal S}^{*,0}$, for any non-negative real number $r$ the map 
\begin{eqnarray*}
\Phi_r^{{\cal R}}: {\cal T}\left({\cal S}^{*,0}\right)&\to& {\cal 
  M}(W)\\
\sigma&\mapsto &
\cutoffint_{r}^{\smop{Chen}}\tilde {\cal R}(\sigma)
\end{eqnarray*}
defines
 a character
from the commutative algebra $\left( {\cal T}\left({\cal S}^{*,0}\right),
  \shu\right)$  to the algebra  ${\cal
  M}(W)$ of meromorphic functions on $W$.
\end{prop}
{\bf Proof:}
By Proposition \ref{prop:Chenmero}, $\Phi_r^{{\cal R}}$ is ${\cal 
  M}(W)$-valued.
Applying  Theorem \ref{thm:cutoffshuffle} to $\tilde{\cal R}( \sigma)(z)$ and
$\tilde {\cal R}(\tau)(z)$ for any two $\sigma, \tau \in  {\cal T}\left({\cal S}^{*,0}\right)$,  
which makes sense outside a discrete set of
complex numbers $z$, we have by (\ref{eq:regshu}): 
\begin{eqnarray*}
\Phi_r^{{\cal R}}(\sigma\shu \tau)&=&\cutoffint^{{\rm Chen}}_r\left(
  \tilde{\cal R}( \sigma \,\shu \, \tau)(z)\right)\\
&=&\cutoffint^{{\rm Chen}}_r\left( \tilde{\cal R}( \sigma)(z) \,\shu \, \tilde{\cal R}( \tau)(z)\right)\\
&= &\left(\cutoffint^{{\rm Chen}}_r\tilde{\cal R}( \sigma)(z)\right) \, \left(\cutoffint^{{\rm Chen}}_{r}
\tilde{\cal R}( \tau)(z)\right)\\
&=& \Phi_r^{{\cal R}}(\sigma)\, \Phi_r^{{\cal R}}( \tau)
\end{eqnarray*}
as an equality of meromorphic functions, from which the above proposition
follows.\endsquare
\begin{thm} \label{thm:renChenint}
Let $W\subset\C$ be an open neighbourhood of $0$. Given a holomorphic regularisation ${\cal R}: \sigma \mapsto \{\sigma(z)\}_
{z\in W}$ on
classical symbols, for any non-negative real number $r$ the regularised integral $\cutoffint_{r}^{{\cal R}}$  extends  to a character
$$
\phi_r^{{\cal R}}: \left({\cal
  T}\left({\cal S}^{*,0}\right),
  \shu\right)\to \C$$ on the tensor algebra equipped with the shuffle
product, given by 
$\phi_r^{{\cal R}}:=\Phi_{r,+}^{{\cal R}}(0)$ from  the Birkhoff decomposition of the ${\cal
  M}(W)$-valued character $\Phi_r^{{\cal R}}$. It coincides with $\cutoffint_r^{{\rm Chen}}$ on tensor products of symbols with non-integer partial sums of orders, i.e.:
$$\phi_r^{{\cal R}}(\sigma_1\otimes \cdots \otimes \sigma_k)=
\cutoffint_r^{{\rm Chen}}\sigma_1\otimes \cdots \otimes \sigma_k$$
 for any  $\sigma_i\in {\cal S}^{\alpha_i,0}, i=1, \ldots, k$ whenever the (left) partial sums of the orders  $\alpha_1+ \alpha_{2}+\cdots +\alpha_m$, $m=1,
  \ldots, k$ are  non-integer valued. 
\end{thm}
{\bf Proof:} For convenience, we temporarily drop the subscript $r$. The Birkhoff decomposition \cite{CK}, \cite{M} for the minimal subtraction scheme reads:
$$\Phi^{{\cal R}}= \left(\Phi_-^{{\cal R}}\right)^{\star -1}\, \star\, \Phi_+^{{\cal
    R}},$$
where  $\star$ is the convolution defined from the deconcatenation coproduct  by
$\phi \star \psi(\sigma)= \phi(\sigma)+\psi(\sigma)+ \sum_{(\sigma)}\phi(\sigma')\cdot \psi(\sigma'')$
with Sweedler's notations and  
where $ \Phi_+^{{\cal    R}}(\sigma) $ and  $\Phi_-^{{\cal R}} (\sigma)$ are
defined inductively on the degree of $\sigma\in {\cal T}\left({\cal S}^{*,0}\right)$.
If $\pi$ denotes the  projection   onto the ``pole part'' $z^{-1}\C[z^{-1}]$
we have:
\begin{equation}\label{eq:birkhoff1} 
\Phi_-^{{\cal
    R}} = -\pi\left[ \Phi^{{\cal R}}+\sum_{(\sigma)} \Phi_-^{{\cal
  R}}(\sigma')\,\Phi^{\cal
  R}(\sigma'') \right], \quad \Phi_+^{{\cal
    R}} =(1 -\pi)\left[ \Phi^{{\cal R}}+\sum_{(\sigma)} \Phi^{{\cal
  R}}_-(\sigma')\,\Phi^{{\cal
  R}}(\sigma'') \right].
\end{equation}
Since $\Phi_+^{{\cal
    R}}$ is also a character on $\left( {\cal T}\left({\cal S}^{*,0}\right), \shu\right)$, it  obeys the
shuffle relation:
\begin{equation}\label{eq:shufflehol}\Phi_+^{{\cal
    R}}\left(\sigma \shu \tau\right)(z)= \Phi_+^{{\cal
    R}}(\sigma)(z) \, \Phi_+^{{\cal
    R}} (\tau)(z)
\end{equation}
which holds as an equality of meromorphic functions holomorphic at $z=0$.
 Setting
$$\phi^{{\cal R}}:= \Phi_+^{{\cal
    R}}(0),$$
and applying the shuffle relation (\ref{eq:shufflehol}) at $z=0$ yields 
$$\phi^{{\cal R}}\left(\sigma\, \shu\, \tau\right)= \phi^{{\cal
    R}}(\sigma)\, \phi^{{\cal R}}(\tau).$$
Now suppose that $\sigma= \sigma_1\otimes \cdots \otimes \sigma_k$ is such
that none of the (left) partial sums of the orders $\alpha_1+ \alpha_{2}+\cdots +\alpha_{m}$, $m=1,
  \ldots, k$ is an integer. Due to the particular form of the
deconcatenation coproduct, the component $\sigma'$ shares the same property
(with $k$ replaced by the degree of $\sigma'$). In other words, the linear
span ${\cal K}$ of such $\sigma$'s is a right co-ideal of the coalgebra ${\cal
  T}({\cal S}^{*,0})$. For any $\sigma\in{\cal K}$, the map $z\mapsto \Phi^{{\cal
  R}}(\sigma)(z)$ is holomorphic at $0$ since poles only occur for values $z_0$ such
that the partial orders lie in $ \Z$. In particular $\pi(\Phi^{{\cal
    R}}(\sigma))=0$, and it moreover shows for $k=1$
that:
\begin{equation}\label{pluss}
\Phi_-^{{\cal
  R}}(\sigma)=0 \hbox{ and }\Phi_+^{{\cal
  R}}(\sigma)= \Phi^{{\cal
  R}}(\sigma).
\end{equation}
 On the grounds of   the above recursive formulas for $\Phi_-^{{\cal R}}$ and
$\Phi_+^{{\cal R}}$ we can then prove inductively on $k$ that (\ref{pluss}) holds
for any $k$, from which we infer that  $$ \phi^{{\cal
  R}}(\sigma)= \Phi^{{\cal
  R}}(\sigma)(0)= \cutoffint^{{\rm Chen}} \sigma_1(0)\otimes \cdots \otimes
\sigma_k(0)=\cutoffint^{{\rm Chen}} \sigma_1\otimes \cdots \otimes \sigma_k$$ since
$\sigma_i(0)=\sigma_i$. \endsquare\\\\
Motivated by this result we set the following definition. 
\begin{defn}\label{def:rci}  Given a regularisation ${\cal R}: \tau \mapsto
  \{\tau(z)\}_{z\in W}$ on
  ${\cal S}^{*,0}$, we define for any non-negative real number $r$ and any $\sigma_1, \ldots,\sigma_k$ the renormalised Chen
  integral of $\sigma_1\otimes \cdots \otimes \sigma_k$ by
$$\cutoffint_r^{{\rm Chen}, {\cal R}} \sigma:= \phi_r^{{\cal R}}(\sigma).$$
\end{defn}
We illustrate this on an example which offers a continuous $n$-dimensional
analog of 
regularised multiple zeta functions familiar to number theorists. For any
$s\in \R$ we set for some non negative real number $v$ $$\varphi_{s,v}(\xi)=
\chi(\xi)( \vert \xi\vert+v)^{-s}$$
where $\chi$ is some smooth cut-off function around $0$. For $\mu \in \R$  let $\gamma_\mu: \C\to \C$ be a holomorphic  function such that
  $\gamma_\mu(z)= z+\mu z^2+ o(z^2)$ for $z\to 0$ (e.g. $\gamma_\mu(z)= \frac{z}{1-\mu
    z}$),
\begin{equation}\label{eq:regmu}{\cal R}_\mu:  \sigma\mapsto (1-\chi)\sigma+
  \chi \, \sigma \, \vert \xi \vert^{-\gamma_\mu(z)}
\end{equation} in which case
$${\cal R}_\mu\left(\varphi_s\right)\sim
    \varphi_{s+\gamma_\mu(z)}.$$
Applying the above construction with $r=0$ yields for non negative real
numbers $v_1, \ldots,v_k$ renormalised quantities: 
$$\tilde \zeta^{\mu}_n(s_1, \ldots, s_k; v_1, \ldots, v_k):= 
\cutoffint^{\smop{Chen}, {\cal R}_\mu} \varphi_{s_1, v_1}\otimes \cdots
\otimes \varphi_{s_k, v_k}.$$ When $v_1=\cdots =v_k=v$ we set
$$\tilde \zeta^{\mu}_n(s_1, \ldots, s_k; v)=\tilde \zeta^{\mu}_n(s_1, \ldots, s_k; v_1, \ldots, v_k).$$
\begin{cor}\label{cor:rencontmultiplezeta} For any $\mu \in \R$, $n\in \N$ and
  for any two words $u,v$ with letters in $\C$ we have:
$$  \tilde \zeta^{\mu}_{n} (u\, \shu\,u^\prime;\, v)=\tilde
\zeta^{\mu}_{n}(u;\ v)\cdot
\tilde \zeta^{\mu}_{n}
(u^\prime; \, v).$$ 
\end{cor}
{\bf Proof:} This is a straightforward consequence of   Theorem
\ref{thm:renChenint} applied to
${\cal R}_\mu$ and $\sigma_i= \varphi_{s_i,\, v}.$
\endsquare
\vfill \eject \noindent
\section{DISCRETE  SUMS OF SYMBOLS}
We denote by $|.|$ the usual absolute value on $\R$. The notation ${\cal
  P}^{\alpha,k}$ stands for {\sl positively supported\/} symbols,
i.e. symbols in ${\cal S}^{\alpha,k}(\R)$ with support included in
$]0,+\infty[$. We keep {\sl mutatis mutandis\/} the notations of subsection
1.1: in particular ${\cal P}^{*,0}$ is a subalgebra of the filtered algebra  ${\cal P}^{*,*}$.
\\ \\

We investigate in this second part a discrete counterpart of the results of
part I. We will mostly stick to dimension $n=1$: the reason for this lies in
the crucial use of the Euler-MacLaurin formula. To the best of our knowledge a higher dimensional
analogue, the Khovanskii-Pukhlikov formula (\cite{KSW2}, \cite{GSW}) exists
only on polytopes, i.e. on balls with respect to norms which are not smooth outside
the origin. The suitable framework for dealing with nested sums is then not
clear. We return however to a higher dimensional situation at the end of the paper
by considering radial functions with respect to the supremum norm $|.|_\infty$
in $\R^n$, namely functions $\sigma=f\circ|.|_\infty$ with $f\in{\cal
  P}^{*,*}$.
\subsection{Discrete  sums versus integrals of symbols}\label{sect:sumint}
A Rota-Baxter operator on an algebra $A$ (defined, say over a field $k$) is a
linear operator $P:A\to A$ such that the relation:
\begin{equation}\label{defrb}
P(f)P(g)=P\big(P(f)g+fP(g)\big)-\theta P(f)P(g)
\end{equation}
holds for any $f,g\in A$. Here $\theta$ is a scalar in the field $k$ called the {\sl weight\/}\footnote{Some authors use the opposite sign convention for the weight.}.
The operator $\tilde { P} :{\cal
  P}^{*,k}\to {\cal P }^{*, k+1}$:
$$\tilde { P}(f) (\eta):= \int_0^\eta f(\xi)\, d\xi$$
satisfies the weight zero Rota-Baxter relation.
On the other hand the operator $P$ defined on sequences $f:\N\to\C$ by:
\begin{equation}\label{operatorP}
P(f)(N)=\sum_{m=0}^N f(k)
\end{equation}
satisfies the Rota-Baxter relation of weight $1$, whereas
the operator $Q=P-\mop{Id}$ satisfies the Rota-Baxter relation of
weight $-1$. We now want to compare the behaviour in the large of the continuous integrals
$\tilde P (f)$  with that of discrete sums $P(f)$ when $f$ is a positively
supported symbol (note that the sum actually begins with $m=1$ in this case).
In order to compare ${ P}(f) (N)$ and $\tilde { P}(f)
(N)$ for large $N$  we make use of the Euler-MacLaurin formula
(\cite{Ha},\cite{KSW1}, see also \cite{MP}). Formula (6) in \cite{KSW1} for
$f\in{\cal P}^{*,*}$ yields:
\begin{eqnarray}\label{eq:eml}
P(f)(N)-\tilde
P(f)(N)&=&\frac 12 f(N)+\sum_{k=2}^{2K}\frac{B_{k}}{k!}f^{(k-1)}(N)\nonumber\\
&+&
\frac{1}{(2K+1)!}\int_0^N\overline{B_{2K+1}}(x)f^{(2K+1)}(x)\,dx.
\end{eqnarray}
with $\overline{B_k}(x)= B_k\left(x-[x] \right)$. Here $
B_k(x)= \sum_{i=0}^k {k\choose i} \, B_{k-i} \, x^k$ are the Bernoulli
polynomials of degree $k$,
the $B_i$ being the Bernoulli numbers, defined by the generating series:
$$\frac {t}{e^t-1}=\sum_i\frac{B_i}{i!}t^i,$$
so that
$$B_0=1,\,B_1=-\frac 12,\,B_2=\frac 16,\,B_4=-\frac 1{30},\ldots,\hskip 12mm
B_{2j+1}=0\hbox{ for }j\ge 1.$$
\begin{prop}\label{prop:interpolation}
For any $f\in {\cal P}^{\alpha, k}$ and any $j\in \N$ chosen large enough the
function $\overline P(f)$ defined by:
\begin{eqnarray}\label{eq:interpolation}
\overline P(f)(\eta)&:=&\tilde P(f)(\eta)+\frac 12 f(\eta)+\sum_{j=2}^{2K}\frac{B_{j}}{j!}f^{(j-1)}(\eta)\nonumber\\
&+&
\frac{1}{(2K+1)!}\int_0^\eta\overline{B_{2K+1}}(x)f^{(2K+1)}(x)\,dx.
\end{eqnarray}
lies   in ${\cal P}^{\alpha+1, k+1}+{\cal
  P}^{0, k+1} $ and the difference $\overline
  P(f)- \tilde P(f)$ lies
  in ${\cal P}^{\alpha, k}$.
\end{prop}
\begin{rk}
Note that when $\eta=N\in \N$ then 
$\overline P(f) (N)= 
  P(f) (N)$. Hence $\overline P$ is an interpolation of $P$. Of course
  $\overline Q:=\overline P-\mop{Id}$ is an interpolation of $Q$ in the same sense.
\end{rk}
{\bf Proof:}
The sum $\frac 12 f(\eta)+\sum_{j=2}^{2K}\frac{B_{j}}{j!}f^{(j-1)}(\eta)$  lies in ${\cal
  P}^{\alpha, k}$, whereas the integral  $\tilde P(f)$ lies in ${\cal
  P}^{\alpha+1, k+1}+{\cal
  P}^{0, k+1}$. The result then follows from  splitting  the
integral remainder term into $\int_0^{+\infty}(...)-\int_\eta^{+\infty}(...)$: the
first term  in the sum is a constant for large enough $K$, and the second term
is a symbol (with respect to the variable $\eta$) with arbitrarily small order
$\alpha-(2K+1)$ as $K$ grows, which lies in ${\cal P}^{\alpha, k}$.\endsquare
\\ \\ 
We are now ready to introduce cut-off sums. 
\begin{cor} 
For any $f\in {\cal P}^{*, *}$ and any positive integer $N$,  the discrete sum  
$\sum_{m=0}^N f(m) $ has an asymptotic expansion when $N\to \infty$ of
 the type (\ref{eq:cutoffasympt}) so that we can define the following cut-off sums as finite parts:
 $$\cutoffsum  f:= \mopl{fp }_{N\to \infty}
 \sum_{m=0}^Nf(m)$$
which coincides with $ \altscutoffsum f:= \mopl{fp }_{N\to \infty}
 \sum_{m=0}^{N-1}f(m)$ whenever the order of $f$ does not belong to $\{
 0,1, 2, 
 \ldots\}$.
 \end{cor}
{\bf Proof:} When $f\in  {\cal P}^{\alpha, k}$ the maps $\eta\mapsto\overline
P(f) (\eta)$ and $\eta\mapsto\overline
P(f) (\eta)-f(\eta)$
lie in ${\cal P}^{*, k+1}$ and therefore have the expected asymptotic behaviour as
$\eta\to \infty$. Setting $\eta=N$ and picking the constant
term in the $N\to \infty$ expansion yields the cut-off sums
$\altcutoffsum f$ and $\altscutoffsum f$. They a priori differ by $ \mopl{fp
}_{N\to \infty}f(N)$ which vanishes whenever $f$ has  order $\alpha  \notin\{
 0,1, 2, 
 \ldots\}$.
\endsquare 
\begin{prop}\label{prop:sumholsymb}
Given a holomorphic regularisation ${\cal R}: h \mapsto h(z)$ on
${\cal P}^{*,k}$, for any $f\in {\cal P}^{*,k}$,  the map
$$z\mapsto \cutoffint_{\R}f(z)-
\cutoffsum  f(z)$$
is holomorphic  for any $f\in  {\cal P}^{*,k}$.\\
Consequently, the map
 $z\mapsto\altcutoffsum f(z)$ is meromorphic with the same 
 poles (of order $\leq k+1$)
as the map  $z\mapsto \cutoffint_{\R}f(z)$, which lie  in the discrete
set $\alpha^{-1}\left(\{ -1,0, 1, 2,\ldots \}\right)$  whenever $f(z)$ is a
simple  holomorphic family of order $\alpha(z)$. Moreover,  $\altcutoffsum
f(z)$  and  $\altscutoffsum
f(z)$  coincide as
meromorphic functions.
\end{prop}
{\bf Proof:}
Let $f\mapsto f(z)$ be a holomorphic
perturbation in ${\cal P}^{*,k}$. By Definition (\ref{eq:interpolation}) in Proposition \ref{prop:interpolation}, the difference ${\rm fp}_{N\to \infty} 
 P(f(z))(N)-{\rm fp}_{N\to \infty} 
 \tilde P(f(z))(N)$
 is holomorphic.  On the other hand,  since the  map $z\mapsto \cutoffint 
 f(z)= {\rm fp}_{N\to \infty} 
\tilde P(f(z))(N)$ is meromorphic with  poles 
of order $\leq k+1$, the same
property therefore holds for  $z\mapsto {\rm fp}_{N\to \infty} 
 P(f(z))(N)= \altcutoffsum 
f(z)$ and hence for  $\altscutoffsum f(z)$. When the regularisation is simple these poles lie in the discrete
set $\alpha^{-1}\left(\{ -1,0, 1, 2,\ldots \}\right)$ where $\alpha(z)$ stands for
the order of $f(z)$. Outside these poles, the  sums  $\altcutoffsum 
f(z)$ and $\altscutoffsum f(z)$ coincide so that they coincide as meromorphic functions.
\endsquare\\ \\
On the grounds of this result  we can set the following definition:
\begin{defn}Given a holomorphic regularisation ${\cal R}: f\mapsto
  \{f(z)\}_{z\in W}$ on
 ${\cal P}^{*,k}$ we define 
 for any $f\in
{\cal P}^{*, k}$ 
 the regularised sum
$$\cutoffsum^{{\cal R}}_{ }f:= \mopl{fp }_{z=0}\cutoffsum f(z)$$
where by finite part we mean the constant term in the Laurent series expansion.
\end{defn}
\subsection{The Riemann and Hurwitz zeta functions as regularised sums of a symbol}
Let $v$ be a non-negative real number. Riesz regularisation ${\cal R}: f(x)\mapsto f(x) \, \vert x\vert^{-z}$
applied to the symbol $f_{s;v}: x\mapsto (x+v)^{-s} $ for any $s\in ]0,+\infty[$ gives rise to the Hurwitz
$\zeta$-function\footnote{see Introduction for the notation convention}:
$$\zeta(s;v):= \cutoffsum^{\cal R} f_{s;v}=\mopl{fp
}_{z=0}\cutoffsum_{n=1}^{\infty} (n+v)^{-s-z}.$$
We recall the following elementary property of the Hurwitz zeta
  function on its domain of convergence:
$$\zeta(s; v+1)=\zeta(s; v)-(v+1)^{-s}.$$
In the following, we  investigate its values together with the
values of its derivatives
$$\frac{d}{dv} \zeta(s;v)= -s\, \zeta(s+1;v)$$ at non-positive
integer points.\\ \\
Let us first recall the Euler-MacLaurin formula  (\cite{Ha}, see also
\cite{C2}, \cite{KSW1}).  For any   $ f\in \Ci(\R^+)$ and for any integer $K\ge 0$,
\begin{eqnarray}\label{eq:EulerMacLaurin} \sum_{n=1}^{N} f(n)&=& \int_1^N f(x)\, dx+\frac{f(1)+f(N)}{2}+
\sum_{k=2}^{2K}\frac{B_k}{k!}\left( f^{(k-1)}(N)-f^{(k-1)}(1)\right)\nonumber\\
&+&
\frac{1}{(2K+1)!} \int_1^N \overline{B_{2K+1}} (x)\, f^{(2K+1)}(x) \, dx,
\end{eqnarray} 
with $\overline{B_k}(x)= B_k\left(x-[x] \right)$
and where $
B_n(x)= \sum_{k=0}^n {n\choose k} \, B_{n-k} \, x^k$ are the Bernoulli polynomials,
the $B_i$ being the Bernoulli numbers. For any $a\in\R$ and $m\in\N-\{0\}$ we
use the notation:
\begin{equation}\label{arrangement}
[a]_m:=a(a-1)\cdots(a-m+1)
\end{equation}
which should not lead to any confusion with the integer part $[a]$ of $a$. The following proposition yields back the known Hurwitz zeta values
  at non positive integers (see e.g. \cite{C2}) as well as their derivatives w.r. to the parameter.
\begin{prop}\label{prop:sumlogn}Let $z\mapsto \gamma(z)$ be a holomorphic
  function with $\gamma'(0)=1$. Let $\lambda\in\C-\{0\}$. 
The map $z\mapsto  \cutoffsum_{n=1}^\infty
(n+v)^{a-\lambda\gamma(z)}$ is holomorphic at zero for any $a\neq -1$. For any
$a\in\N$  and any  non negative real $v$ the Hurwitz zeta value $\zeta(-a;v)=-\frac{B_{a+1}(1+v)}{a+1}$ can be obtained as:
\begin{equation}\label{eq:zetaneg}
\zeta(-a;v)=\mopl{fp}_{z=0} \cutoffsum_{n=1}^\infty
(n+v)^{a-\lambda\gamma(z)}=\lim_{z\to 0}\cutoffsum_{n=1}^\infty
(n+v)^{a-\lambda\gamma(z)}.\end{equation}
In particular it is a rational number when $v$ is rational. 
When $a=-1$ the residue at $0$ which reads $\frac{1}{\lambda}$, is
  independent of $v$. Similar formulae
hold for the derivatives at non positive integers:
\begin{eqnarray*}\frac{d}{dv} \zeta(-a;v)&=&a\, \zeta(-a;v)\quad {\rm if}\quad
a\neq 0\\
&=&-1\quad {\rm if}\quad
a= 0.
\end{eqnarray*}
\end{prop}
{\bf Proof:} Applied to $f(x)=(x+v)^a$ with $a\in \C$,   the Euler-MacLaurin formula  
gives:
\begin{eqnarray}\label{eq:cutoffsumN}
 &\hskip -8mm&\sum_{0<n\le N} (n+v)^a=\frac{(N+v)^a +(1+v)^a}{2}+ \int_1^N
(x+v)^a \, dx\nonumber \\
&\hskip -15mm+&\hskip -8mm
\sum_{k=2}^{2K}  B_k\,\frac{ [a]_{k-1}}{k!}\left(
  (N+v)^{a-k+1}- (1+v)^{a-k+1}\right)+
\frac{[a]_{2K+1}}{(2K+1)!} \int_1^N \overline{B_{2K+1}}
(x)\,(x+v)^{a-2K-1} \, dx\nonumber\\
&\hskip -15mm=&\hskip -8mm  (1-\delta_{a+1}) \frac{(N+v)^{a+1}}{a+1}- 
(1-\delta_{a+1}) \frac{(1+v)^{a+1}}{a+1}+\frac{(N+v)^a +(1+v)^{a}}{2}  +\delta_{a+1} \,( \log (N+v) -\log(1+v))
\nonumber\\
&+&
\sum_{k=2}^{2K}  B_k\,\frac{[a]_{k-1}}{k!}\left(
  (N+v)^{a-k+1}-(1+v)^{a-k+1}\right)+
\frac{ [a]_{2K+1}}{(2K+1)!} \int_1^N \overline{B_{2K+1}}
(x)\,(x+v)^{a-2K-1} \, dx.
\end{eqnarray}
Let us set 
$$R_K(a):=
\frac{ [a]_{2K+1}}{(2K+1)!} \int_1^N \overline{B_{2K+1}}
(x)\,(x+v)^{a-2K-1} \, dx; \quad S_K(a):=\sum_{k=2}^{2K}
B_k\,\frac{[a]_{k-1}}{k!}\, (1+v)^{a-k+1}.$$
\\ \\
 Replacing  $a$  in 
  (\ref{eq:cutoffsumN})  by $a-\lambda\gamma(z)$ and taking finite parts as
  $N\to\infty$ we have:
\begin{equation}\label{eq:proof2}
\cutoffsum_1^\infty (n+v)^{a-\lambda\gamma(z)}
=-\frac{(1+v)^{a-\lambda\gamma(z)+1}}{a-\lambda\gamma(z)+1}+\frac{(1+v)^{a-\lambda\gamma(z)}}{2}
-S_K\big(a-\lambda\gamma(z)\big)
+R_{K}\big(a-\lambda\gamma(z)\big).
\end{equation}
Hence ${\rm Res}_{z=0} \cutoffsum_1^\infty (n+v)^{a-\lambda\gamma(z)}=
\delta_{a+1}\,\frac{1}{\lambda}$. 
Taking the finite part at  $ 0$ then yields:
\begin{equation*}
 {\rm fp}_{z= 0} \cutoffsum_{1}^\infty (n+v)^{a-\lambda\gamma(z)} 
=-(1-\delta_{a+1}) \frac{(1+v)^{a+1}}{a+1}-\delta_{a+1}\,\log (1+v)+\frac{(1+v)^a}{2}
-S_K(a)
+R_{K}(a),
\end{equation*}
which for a non negative integer $a$ gives:
\begin{equation*}
\lim_{z\to 0}\cutoffsum_{0}^\infty (n+v)^{a-\lambda\gamma(z)} 
=-\sum_{k=0}^{a+1} {a+1\choose k}\,B_k\,
(1+v)^{a+1}= -\frac{B_{a+1}(1+v)}{a+1}.
\end{equation*}
We now turn to derivatives. Since $\frac{d}{dv} \zeta(s;v)= -s\,
  \zeta(s+1;v)$, at a non negative integer $a$ we have:
\begin{eqnarray*}
\frac{d}{dv} \zeta(-a;v)&=&\mopl{fp}_{z= 0}\big((a-\lambda\gamma(z)\big)\,
  \cutoffsum_{1}^\infty (n+v)^{a-\lambda\gamma(z)}\\
&=&a\,\mopl {fp}_{z= 0}
  \cutoffsum_{1}^\infty (n+v)^{a-\lambda\gamma(z)}\quad {\rm if}\quad
a\neq 0\\
&=&-\lambda \,  \mopl{Res}_{z= 0}
  \cutoffsum_{1}^\infty (n+v)^{a-\lambda\gamma(z)}\quad {\rm if}\quad
a= 0\\
\end{eqnarray*}
so that 
\begin{eqnarray*}\frac{d}{dv} \zeta(-a;v)&=&a\, \zeta(-a;v)\quad {\rm if}\quad
a\neq 0\\
&=& -\lambda \, \frac{1}{\lambda}= -1\quad {\rm if}\quad
a= 0.
\end{eqnarray*}

\endsquare
\vfill \eject \noindent
\section{NESTED  SUMS OF SYMBOLS}
\subsection{The Hoffman isomorphism and stuffle relations for discrete nested sums}\label{sect:9}
We recall here the explicit construction by M. Hoffman (\cite{H2}) of an
isomorphism between the stuffle Hopf algebra and the shuffle Hopf algebra.
\begin{defn}\label{def:qsh}
Let $k,l,r\in\N$ with $k+l-r>0$. A {\rm $(k,l)$-quasi-shuffle\/} of type $r$
is a surjective map $\pi$ from $\{1,\ldots ,k+l\}$ onto
$\{1,\ldots, k+l-r\}$ such that $\pi(1)<\cdots
<\pi(k)$ and $\pi(k+1)<\cdots <\pi(k+l)$. We shall denote by $\mop{mix
  sh}(k,l;r)$ the set of $(k,l)$-quasi-shuffles of type $r$. The elements of $\mop{mix sh}(k,l;0)$ are the ordinary
$(k,l)$-shuffles. Quasi-shuffles are also called {\rm mixable shuffles\/} or
{\rm stuffles\/}. We denote by $\mop{mix sh}(k,l)$ the set of
$(k,l)$-quasi-shuffles (of any type).
\end{defn}
Let $({\cal A},\bullet)$ be a commutative (not necessarily unital) algebra. Let $\Delta$
be the deconcatenation coproduct on ${\cal T}({\cal A})=\bigoplus_{k\ge 0}{\cal A}^{\otimes
  k}$, let $\shu$ be the shuffle product, and let $\star_\bullet$ the product
on ${\cal T}({\cal A})$ defined by:
$$(v_1\otimes\cdots\otimes v_k)\star_\bullet (v_{k+1}\otimes\cdots\otimes v_{k+l})
=\sum_{\pi\in\smop{mix sh}(k,l)}w^\pi_1\otimes\cdots\otimes w^\pi_{k+l-r},$$
with~:
$$w^\pi_j=\prod_{i\in\{1,\ldots,k+l\},\pi(i)=j}v_i.$$
(the product above is the product $\bullet$ of ${\cal A}$, and contains only one or
two terms).
\begin{thm}\label{th:hoffman}{\rm (M. Hoffman, \cite{H2} theorems 3.1 and 3.3)}\\
\begin{itemize}
\item $\big({\cal T}({\cal A}),\star_\bullet,\Delta\big)$ is a commutative connected filtered
  Hopf algebra.
\item  There 
is an isomorphism of Hopf algebras~:
$$\exp:\big(T({\cal A}),\shu,\Delta\big) \mathop{\longrightarrow}\limits^{\sim}\big(T({\cal A}),\star_\bullet,\Delta\big).$$
\end{itemize}
\end{thm}
M. Hoffman in \cite{H2} gives a detailed proof in a slightly more restricted
context, which can be easily adapted in full generality (see also \cite{EG}). Hoffman's isomorphism is built explicitly as follows:
let ${\cal P}(n)$ be the set of compositions of the integer $n$, i.e. the set of
sequences  $I=(i_1,\ldots,i_k)$ of positive integers such that
$i_1+\cdots +i_k=n$. For any  $u=v_1\otimes\cdots\otimes v_n\in T({\cal A})$
and any composition  $I=(i_1,\ldots,i_k)$ of $n$ we set:
$$I[u]:=(v_1\bullet\cdots\bullet v_{i_1})\otimes
(v_{i_1+1}\bullet\cdots\bullet v_{i_1+i_2})
\otimes\cdots\otimes
(v_{i_1+\cdots+i_{k-1}+1}\bullet\cdots\bullet v_n).$$
We then further define:
$$\exp u=\sum_{I=(i_1,\ldots ,i_k)\in{\cal P}(n)}\frac{1}{i_1!\cdots
  i_k!}I[u].$$
Moreover (\cite{H2}, lemma 2.4), the inverse $\log$ of $\exp$ is given by~:
$$\log u=\sum_{I=(i_1,\ldots ,i_k)\in{\cal P}(n)}\frac{(-1)^{n-k}}{i_1\cdots
  i_k}I[u].$$
For example for  $v_1,v_2,v_3\in {\cal A}$ we have~:
\begin{eqnarray*}
\exp v_1 =v_1 &,& \log v_1 =v_1,\\
\exp (v_1\otimes v_2) =v_1\otimes v_2+\frac{1}{2}v_1\bullet v_2&,&
\log (v_1\otimes v_2) =v_1\otimes v_2-\frac{1}{2}v_1\bullet v_2,\\
\exp (v_1\otimes v_2\otimes v_3) =v_1\otimes v_2\otimes v_3&+&\frac{1}{2}(v_1\bullet v_2\otimes
v_3+v_1\otimes v_2\bullet v_3)+\frac{1}{6}v_1\bullet v_2\bullet v_3,\\
\log (v_1\otimes v_2\otimes v_3) =v_1\otimes v_2\otimes v_3&-&\frac{1}{2}(v_1\bullet v_2\otimes
v_3+v_1\otimes v_2\bullet v_3)+\frac{1}{3}v_1\bullet v_2\bullet v_3.
\end{eqnarray*}
Let $V$ be the space of sequences $f:\N\to\C$ such that $f(0)=0$. It is a commutative algebra for
the ordinary product.
\begin{thm}\label{thm:Hopfstuffle}
\begin{enumerate}
\item  Consider the commutative algebra $(V,\bullet)$ where $\bullet$ is the ordinary product. For
any $N\in \N$, the  truncated discrete nested sums on ${\cal T}(V)$ defined by:
 $$\sum_<^{ N, \smop{Chen}} f_1\otimes \cdots \otimes f_k:= \sum_{0<
 n_k < \cdots <
 n_1 < N} f_1(n_1)\cdots f_k(n_k)$$
for any $f=f_1\otimes
\cdots\otimes f_k\in V^{\otimes k}$ fulfill the following relations:
\begin{equation}\label{eq:Nstuffle}
\sum_<^{ N, \rm Chen} f\, 
 \star_\bullet \, 
 g = 
 \left(\sum_<^{N, \rm Chen} f \right)\,
 \left(\sum_<^{ N, \rm Chen} g\right).
\end{equation}
\item Whenever the nested sums converge as $N\to \infty$, we have in the
limit:
\begin{equation}\label{eq:convstuffle}\sum_<^{\rm Chen} f\, 
 \star_\bullet \, 
 g = 
 \left(\sum_<^{{\rm Chen}}  f\right)\,
 \left(\sum_<^{\rm Chen}
 g\right).
 \end{equation}
\item The same statements hold with the weak inequality version:
$$\sum_\le^{ N, \smop{Chen}} f_1\otimes \cdots \otimes f_k:= \sum_{0<
 n_k \le \cdots \le 
 n_1 \le N} f_1(n_1)\cdots f_k(n_k)$$
provided $\bullet$ is now the
  {\rm opposite\/} of the ordinary product, i.e. $v_1\bullet v_2=-v_1v_2$.
\end{enumerate}  
\end{thm}\goodbreak
 {\bf Proof:} It is enough to prove the theorem for $f=f_1\otimes\cdots\otimes
 f_k\in V^{\otimes k}$ and $g=f_{k+1}\otimes\cdots\otimes
 f_{k+l}\in V^{\otimes l}$.
 Let us prove the theorem for the weak inequality case. The domain:
$$P_{k,l}:=\{n_1>\cdots >n_k>0\} \times \{n_{k+1}>\cdots >n_{k+l}>0\}
\subset( \N-\{0\})^{k+l}$$
is partitioned into:
$$P_{k,l}=\coprod_{\pi\in\smop{mix sh}(k,l)}P_\pi
,$$
where the domain $P_\pi$ is defined by:
$$P_\pi=\{(n_1,\ldots ,n_{k+l})\,/\,
n_{\pi_m}>n_{\pi_{p}}\hbox{ if } m>p \hbox{ and }\pi_m\not =\pi_p,
\hbox{ and } n_m=n_{p} \hbox{ if } \pi_m=\pi_{p}\}.$$
As we must replace strict inequalities by large ones, let us consider the
``closures'' 
$$\overline{P_\pi}:=\{(n_1,\ldots ,n_{k+l})\,/\,
n_{\pi_m}\ge n_{\pi_{p}}\hbox{ if } m\ge p 
\hbox{ and } n_m=n_{p} \hbox{ if } \pi_m=\pi_{p}\}.$$
which then overlap. By the inclusion-exclusion
principle we have:
\begin{equation}\label{eq:overlaps}
\overline{P_{k,l}}=\coprod_{0\le r\le \smop{min}(k,l)}(-1)^r\coprod_{\pi\in\smop{mix sh}(k,l;r)}\overline{P_\pi},
\end{equation}
where we have set:
$$\overline{P_{k,l}}:=\{n_1\ge
\cdots \ge n_k>0\} \times \{n_{k+1}\ge\cdots \ge n_{k+l}>0\}
\subset (\N-\{0\})^{k+l}$$
Each term in equation (\ref{eq:overlaps}) must be added if $r$ is even, and
removed if $r$ is odd. Considering the summation of $f_1\otimes\cdots\otimes
f_{k+l}$ over each $\overline{P_\pi}$, this decomposition immediately yields the equality:
\begin{eqnarray}\label{eq:decomp}
&{}&\left(\sum_{0\le n_{k}\le\cdots \le
n_{1}\le N}      f_1(n_1)\cdots f_{i}(n_{k})\right)\,\left(\sum_{0\le n_{k+l} \le\cdots \le
n_{k+1}\le N}f_{k+1}(n_{k+1})\cdots f_{k+l}(n_{k+l})\right)\nonumber\\
 &=& \sum_{\leq}^{N, \smop{Chen}} \sum_{\pi\in \smop{mix sh}(k,l)}f^
 \pi, 
 \end{eqnarray} 
where $f^\pi= f^\pi_1\otimes \cdots \otimes f^\pi_{k+l-r}$ is
the tensor product defined by:
$$f^\pi_j=\mopl{$\bullet$}_{i\in \{1, \ldots, k+l\}, \, \pi(i)=j} 
f_i.$$
The stuffle relations (\ref{eq:Nstuffle}) are then a
re-writing of equality (\ref{eq:decomp}) using the commutative algebra
$(V,\bullet)$. Taking the limit as $N\to \infty$ provides the second statement of the
  theorem. The proof is similar for the strict inequality case, using the domains $P_\pi$ rather than the
``closures'' $\overline{P_\pi}$. As there are no overlaps the signs
disappear in the formula (\ref{eq:overlaps}).
\endsquare
\subsection{Cut-off nested sums of symbols and multiple (Hurwitz)  zeta functions}\label{sect:cutoffchen}
In this
section we iterate the interpolated summation operators
$\overline P$ and $\overline Q$ defined on ${\cal P}^{*, *}$ by
\eqref{eq:interpolation} in a  way to be
made precise. As a consequence of Proposition \ref{prop:interpolation}
 we derive the following result. 
\begin{prop} \label{prop:nestedsums}
Given $\sigma_i\in{\cal P}^{\alpha_i,0}$, $i=1,\ldots ,k$, setting $\sigma:=\sigma_1\otimes\cdots\otimes\sigma_k$, the function $\widetilde \sigma$ defined by:
\begin{equation}\label{sigmatilde}
\widetilde\sigma:=\sigma_1\overline P\Big(\cdots\sigma_{k-2}\overline P\big(\sigma_{k-1}\overline P(\sigma_k)\big)...\Big)
\end{equation}
lies in ${\cal P}^{*, k-1}$ as a linear combination of
(positively supported) symbols in  ${\cal P}^{\alpha_1+\cdots +\alpha_m+m-1, m-1}$,
$m\in\{1,\ldots ,k\}$. It has real order $\omega_1$ as defined in Lemma
\ref{lem:CStoCSlog}. The same holds for the function $(\widetilde\sigma)'$
defined the same way with $\overline Q=\overline P-\mop{Id}$ instead of $\overline P$.
\end{prop}
{\bf Proof:} Let us first observe that by
  Proposition \ref{prop:interpolation}, given two symbols $\tau_1\in {\cal P}^{\beta_1, j_1}$ and  $\tau_2\in {\cal P}^{\beta_2, j_2}$,  the expression $\tau_1\, \overline
  P(\tau_2)$ lies in ${\cal P}^{*, j_1+j_2+1}$ as 
 a linear combination of symbols  in ${\cal P}^{\beta_1,j_1}$ and  ${\cal P}^{\beta_1+\beta_2+1, j_1+j_2+1}$. In particular, setting
  $\beta_i=\alpha_i$ and  $j_1=j_2=0$ shows  the proposition
  holds when $k=2$. \\The statement for  $k>2$ then follows from an induction
procedure on $k$. Assuming that $\tau_2:=\sigma_2\overline
P\Big(\cdots\sigma_{k-2}\overline P\big(\sigma_{k-1}\overline
P(\sigma_k)\big)...\Big)$ lies in ${\cal P}^{*, k-2}$ as  a linear combination 
 of  log-polyhomogeneous symbols $\tau_{2,m}$  in   ${\cal P}^{\alpha_2+\cdots
   +\alpha_m+m-2, m-2}$,
$m$ varying in $\{2,\ldots ,k\}$,  we infer from our preliminary observation applied to
$\tau_1=\sigma_1$ and each $\tau_{2,m}$ that $\sigma_1\,
\overline P(\tau_2)$  lies in  ${\cal P}^{*, k-1}$
as  a linear combination 
 of  log-polyhomogeneous symbols $\sigma_1\,
\overline P(\tau_{2,m})$ in   ${\cal P}^{\alpha_1+\alpha_2+\cdots
+\alpha_m+m-1, m-1}$, $m\in\{1,\ldots, k\}$. 

\endsquare\\ \\
We are now ready to define discrete nested sums of symbols. Combining Propositions
\ref{prop:nestedsums}  and \ref{prop:interpolation}
shows that the cut-off sums of the symbols $\widetilde\sigma$ and $(\widetilde\sigma)'$ are well defined
  so that we can set the following definitions. 
\begin{defn} For $\sigma_1, \ldots, \sigma_k \in {\cal P}^{*,0}$ and $\sigma:=\sigma_1\otimes\cdots\otimes\sigma_k$ 
we call 
\begin{equation*}
\cutoffsum^{\smop{Chen}}_{\leq}\sigma
:=\cutoffsum_{n\in\N} \widetilde\sigma(n)
= \cutoffsum_{0<  n_k\le\cdots \le  n_1} \sigma_1(n_1) \cdots  \sigma_k(n_k),
\end{equation*}
the cut-off nested sum of $\sigma=\sigma_1\otimes \cdots \otimes\sigma_k$. The strict inequality version is defined by:
\begin{equation*}
\cutoffsum^{\smop{Chen}}_{<}\sigma
:=\cutoffsum_{n\in\N} (\widetilde\sigma)'(n)
= \cutoffsum_{0<  n_k<\cdots < n_1} \sigma_1(n_1) \cdots  \sigma_k(n_k).
\end{equation*}
\end{defn}
\begin{rk}
This definition of cut-off nested sums (in its weak inequality version) matches with the one given in \cite{MP}
paragraph 6.3. The restriction to positively supported symbols allows us to
drop the absolute values here.
\end{rk}
\begin{cor}\label{cor:convChensums}
The discrete nested sums above converge  whenever $\mop{Re
}(\alpha_1+\cdots+\alpha_m)<-m$  for any $m\in\{1,\ldots ,k\}$ in which case
they are ordinary discrete nested sums:
\begin{eqnarray*}\cutoffsum^{\smop{Chen}}_{\leq}\sigma &=& \sum_{0<  n_k\le\cdots \le n_1}
\sigma_1(n_1)\otimes  \cdots \otimes \sigma_k(n_k)\\
\hbox{\it and }\hskip 4mm\cutoffsum^{\smop{Chen}}_{<}\sigma
&=& \sum_{0<  n_k <\cdots < n_1}
\sigma_1(n_1)\otimes  \cdots \otimes \sigma_k(n_k)
\end{eqnarray*}
\end{cor}
  {\bf Proof:}
   By Proposition \ref{prop:interpolation}, the orders of
   $\widetilde\sigma=\sigma_1\overline P\Big(\cdots\sigma_{k-2}\overline
   P\big(\sigma_{k-1}\overline P(\sigma_k)\big)...\Big)$ and ${\widetilde\sigma}'=\sigma_1\overline Q\Big(\cdots\sigma_{k-2}\overline
   Q\big(\sigma_{k-1}\overline Q(\sigma_k)\big)...\Big)$ coincide with
     that of $\sigma^{\smop{Chen}}=\sigma_1\tilde P\Big(\cdots\sigma_{k-2}\tilde
   P\big(\sigma_{k-1}\tilde P(\sigma_k)\big)...\Big)$. Since the latter
    is smaller than $-1$   under the assumptions of the corollary, the
    convergence  follows from the definition of the nested sums. 
   \endsquare\\ \\
In the latter case, we drop the bar across the summation sign and write
$\sum_\le^{\smop{Chen}}\sigma$ and $\sum_<^{\smop{Chen}}\sigma$. Let us  now apply the above results to $\sigma_{i}(\xi):=\sigma_{s_i,v_i}(\xi):= 
 (\xi+v_i)^{-s_i}\, \chi(x),$
where $s_1, \ldots, s_k$ are complex numbers and  $v_1, \ldots, v_k$   non negative real numbers. Here $\chi$ is a cut-off function 
which vanishes on $]-\infty,\varepsilon[$ with $\varepsilon>0$ and such that
$\chi(\xi)=1$ for $|\xi|\ge 1$. This will give back the
multiple  zeta functions familiar to number theorists as well as  multiple Hurwitz
zeta functions which provide a multiple analog of ordinary Hurwitz zeta
functions (see e.g.\cite{C2})\footnote{The case $k=1$ gives back the  Hurwitz
  zeta function $\zeta(s,v+1)$. We choose these conventions  in order to
  deal with  multiple zeta and multiple Hurwitz zeta  functions simultaneously.}.
\begin{thm}\label{thm:mzext}
If  $\mop{Re }(s_1+\cdots +s_m)>m$ for any $m\in\{1,\ldots ,k\}$ the
discrete nested sums:
\begin{equation*}
\zeta(s_1,\ldots, s_k;\, v_1, \ldots, v_k):=\sum_<^{\smop{Chen}}\sigma_{s_1,v_1}\otimes \cdots\otimes\sigma_{s_k,v_k}
= \sum_{1\leq   n_{k} <  n_{k-1}<
\cdots< n_1} (n_{k}+v_k)^{-s_k} \cdots 
 (n_1+v_1)^{-s_1}
\end{equation*}
for non negative $v_1, \ldots, v_k$ converge and extend meromorphically to all $s_i\in \C$ by the cut-off nested sum:
\begin{equation*}
\zeta(s_1,\ldots, s_k;\,v_1, \ldots,v_k):=\cutoffsum^{{\rm Chen}}_{<}\sigma_{s_1,v_1}\otimes \cdots\otimes\sigma_{s_k,v_k}.
\end{equation*}
Setting $v_1=\cdots= v_k=0$ gives similar statements for multiple zeta functions
\begin{equation*}
\zeta(s_1,\ldots, s_k):=\cutoffsum^{{\rm Chen}}_{<}\sigma_{s_1}\otimes \cdots
\otimes\sigma_{s_k}.
\end{equation*}
 A  similar statement holds for the weak inequality version
 \begin{equation*}
\overline\zeta(s_1,\ldots, s_k;\, v_1, \ldots, v_k):=\sum_\le^{\smop{Chen}}\sigma_{s_1,v_1}\otimes \cdots\otimes\sigma_{s_k,v_k}
= \sum_{1\leq   n_{k} \le  n_{k-1}\le
\cdots\le n_1} (n_{k}+v_k)^{-s_k} \cdots 
 (n_1+v_1)^{-s_1}.
\end{equation*}
 \end{thm}

{\bf Proof:} The existence of the extension follows immediately from applying
the results of Corollary \ref{cor:convChensums} to  $\sigma_i=\sigma_{s_i,v_i}$  of
order $-s_i$. We postpone the proof of
meromorphicity to  the next section (Corollary \ref{cor:meroChensums}).\endsquare\\ \\ 
We have the following relations between both versions (see \cite{H}):
\begin{eqnarray}\label{eq:conversion}
\overline\zeta(a_1,\ldots ,a_k; \,v_1,\ldots, v_k)&=&\sum_{I=(i_1,\ldots ,i_r)\in{\cal
    P}(k)}\zeta(b^I_1,\ldots b^I_r;\, v_1,\ldots, v_r), \nonumber\\
\zeta(a_1,\ldots ,a_k; \, v_1,\ldots, v_k)&=&\sum_{I=(i_1,\ldots ,i_r)\in{\cal
    P}(k)}(-1)^{k-r}\overline\zeta(b^I_1,\ldots b^I_r;\, v_1,\ldots, v_r)
\end{eqnarray}
with $b^I_s:=a_{i_1+\cdots +i_{s-1}+1}+\cdots+a_{i_1+\cdots
  +i_{s}}$. For any multi-index $\alpha=(\alpha_1,\ldots ,\alpha_k)$ we have the following identity between meromorphic functions of $(s_1,\ldots ,s_k)$:
\begin{equation}
\partial^{\alpha_1}_{v_1}\cdots\partial^{\alpha_k}_{v_k}\zeta(s_1,\ldots ,s_k;v_1,\ldots ,v_k)
=(-1)^{|\alpha|}s_1^{\alpha_1}\cdots s_k^{\alpha_k}\zeta(s_1+\alpha_1,\ldots ,s_k+\alpha_k;v_1,\ldots ,v_k).
\end{equation}
Higher-dimensional analogues of multiple zeta functions are
discussed in Section \ref{sect:hdim}. Let us also mention the following vanishing result:
\begin{lem}\label{lem:fpcutoffsumhom} For any non-negative integers  $a_1,
  \ldots, a_k$ and any rational numbers $v_1,\ldots v_k$, the  expression 
$$ \sum_{0<n_k<\cdots < n_1\le N} (n_1+v_1)^{a_1}\cdots
  (n_k+v_k)^{a_k}$$
is a polynomial expression in  $N$ with rational coefficients. The following corresponding cut-off sum vanishes:
$${\rm fp}_{N\to \infty}\sum_{0<n_k<\cdots <
  n_1\le N}
(n_1+v_1)^{a_1}\cdots (n_k+v_k)^{a_k}=0.$$ 
\end{lem}
{\bf Proof:} First observe that once we know that the expression:

\begin{equation}f(N):= \sum_{0<n_k<\cdots <
  n_1\le N} (n_1+v_1)^{a_1}\cdots
  (n_k+v_k)^{a_k}
\end{equation}
 is a polynomial expression in  $N$, then its finite part as $N\to+\infty$
  corresponds to its value at $N=0$. Next, for any rational $v$ the
  polynomials $D_m^v:N\mapsto (N+v)^m-(N+v-1)^m, m\ge 1$ form a basis of the vector space $\Q[N]$. This means that for any polynomial $P\in\Q[N]$ there is a polynomial $Q\in\Q[N]$ (defined up to an additive constant) such that $P(N)=Q(N+v)-Q(N+v-1)$. We obviously get for $N\ge 1$:
  $$\sum_{n=1}^{N}P(n)=Q(N+v)-Q(v),$$
  hence a polynomial expression in $N$ with rational coefficients, which
  vanishes at $N=0$. The lemma follows then by induction on the depth $k$, by writing:
  $$f(N)= \sum_{n_1=1}^{N} (n_1+v)^{a_1}\left(\sum_{0<n_k<\cdots <
  n_2< n_1} (n_2+v)^{a_2}\cdots
  (n_k+v)^{a_k}\right),$$
  since the expression between parentheses is in $\Q[n_1]$ by the induction assumption. 
  \endsquare
\subsection{Renormalised nested sums of symbols}\label{sect:renormchen}
In this section we mimic the construction of renormalised nested integrals of
symbols carried out previously under the requirement that they obeyed shuffle
relations. Via the Euler-MacLaurin formula, which  provides a bridge to
nested sums, we can similarly build meromorphic families of nested sums  of
symbols on $\R$.
\begin{thm}\label{thm:meroChensums}  Let $W\subset \C$ be an open
  neighbourhood of $0$. Let $\sigma_1, \ldots ,\sigma_l \in {\cal
      P}^{*,0}$, and simple holomorphic perturbations ${\cal R}_1(\sigma_1)(z), \ldots ,{\cal R}_l(\sigma_l)(z) \in {\cal
      P}^{*,0}$ with holomorphic orders $\alpha_1(z), \ldots, \alpha_l(z)$,
    $z\in W$,
    such that $\mop{Re }\big(\alpha_1^\prime(z)+\cdots +\alpha_m^\prime(z)\big)< 0$ for any $m\in\{1,\ldots ,m\}$ and for any $z\in
    W$. Let $\sigma:=\sigma_1\otimes\cdots\otimes\sigma_l$ and
    $\sigma(z):=\sigma_1(z)\otimes\cdots\otimes\sigma_l(z)$ where we have set
    $\sigma_i(z):= {\cal R}_i(\sigma_i)(z).$
 Then
\begin{enumerate}
\item The maps 
 $$(z_1, \ldots, z_k)\mapsto \altcutoffsum^{\smop{Chen}}_{<}\sigma_1(z_1)\otimes \cdots \otimes
 \sigma_l(z_k)$$
are meromorphic with poles on a countable number of hypersurfaces
$$\sum_{i=1}^j \alpha_i(z_i)\in 
-j +\N,$$ 
of multiplicity $j$ varying in $\{1, \ldots, k\}$. 
\item
the maps $z\mapsto 
\altcutoffsum^{\smop{Chen}}_{\leq}\sigma(z)$ and $z\mapsto 
\altcutoffsum^{\smop{Chen}}_<\sigma(z)$
are meromorphic on $W$ with poles of order $\leq k$.
\item If $\mop{Re }(\alpha_1+\cdots +\alpha_j)<-j$ for any $j\in\{1,\ldots
  ,k\}$, the nested sums converge and boil down to ordinary nested  sums (independently
  of the perturbation):
\begin{eqnarray*}\cutoffsum^{\smop{Chen}, {\cal R}}_{<}
 \sigma&=&
 \lim_{z\to 0} \cutoffsum^{\smop{Chen}}_{<}
 \sigma(z)=
 \sum_<^{\smop{Chen}}\sigma,\\
\cutoffsum^{\smop{Chen}, {\cal R}}_{\leq}
 \sigma&=&
 \lim_{z\to 0} \cutoffsum^{\smop{Chen}}_{\leq}
 \sigma(z)=
 \sum^{\smop{Chen}}_{\leq}\sigma.
\end{eqnarray*}
These convergent nested sums moreover obey  stuffle relations:
\begin{eqnarray}\label{eq:stuffle}
\sum_<^{\smop{Chen}}\sigma\, \star_\bullet \, \tau&=& 
\sum_<^{\smop{Chen}}\sigma\,\cdot  \,\sum_<^{\smop{Chen}}\tau,\\
\sum^{\smop{Chen}}_{\leq}\sigma\, \star_\bullet \, \tau&=& 
\sum^{\smop{Chen}}_{\leq}\sigma\,\cdot  \,\sum^{\smop{Chen}}_{\leq }\tau,
\end{eqnarray}
where $\bullet$ stands for the ordinary product in the first
case, and stands for the opposite of the ordinary product in the second case.
\end{enumerate}
\end{thm} 
{\bf Proof:} We give  the proof in  the weak inequality case since the strict inequality case
 is  similar.
\begin{enumerate}
\item With the notations of Proposition \ref{prop:nestedsums} and setting  
$\sigma(\underline z_k):= \sigma_1(z_1)\otimes \cdots \otimes 
\sigma_k(z_k)$  we have
\begin{equation*}
 \cutoffsum^{\smop{Chen}}_{\leq}\sigma(\underline z_k) = \cutoffsum \widetilde{\sigma(\underline z_k)} =
 \mopl{fp}_{N\to \infty}\overline P\left(\widetilde{\sigma(\underline z_k)}  \right)(N).
\end{equation*}
By Proposition \ref{prop:nestedsums}, the map $ \widetilde{\sigma(\underline
  z_k)}$ lies  in
${\cal P}^{*,k-1}$  as a linear combination of symbols $\tau_j(z_1, \ldots, z_j)$ in 
${\cal P}^{\sum_{i=1}^j \alpha_i(z_i)+j-1,j-1}$, $j$ varying from $1$ to
$k$. \\  By Proposition \ref{prop:sumholsymb}, 
the  maps $(z_1, \ldots, z_j) \mapsto \cutoffsum \tau_j(z_1, \ldots, z_j)  $ are
 meromorphic with poles  of order $\leq j$ on a countable set of  hypersurfaces
 $\sum_{i=1}^j \alpha_i(z_i)\in  [-j, \infty[\cap \Z $. Consequently, the  map $(z_1, \ldots, z_l)\mapsto 
\cutoffsum^{\smop{Chen}}_{\leq}\sigma(\underline z_k)$ is  meromorphic with
poles of order
$\leq j$ on the countable set of hypersurfaces  $\sum_{i=1}^j \alpha_i(z_i) \in 
 [-j, \infty[\cap \Z $,   $j$ varying in $\{1, \ldots, k\}$. 
\item Setting $z_j=z$ yields the second item of the proposition.
\item
By Corollary \ref{cor:convChensums} we know that under the
assumptions of the theorem the nested sums converge. Consequently, the meromorphic
map  $\cutoffsum^{\smop {Chen}}_{\leq }
 \sigma(z)$  is holomorphic at zero and the
 value at $z=0$  coincides with the ordinary nested sum. The stuffle
relations then hold as a result of (\ref{eq:convstuffle}). 
\end{enumerate}
\endsquare\\ Applying  this theorem to $\sigma_i(z_i)(\xi)=
\chi(\xi)\, \xi^{-z_i} $ (resp. $\sigma_i(z_i)(\xi)=
 (\xi+v_i)^{-z_i}, v_i> 0$)  leads to the following properties of
multiple zeta functions.
\begin{cor}\label{cor:meroChensums}
The extensions of multiple (resp. Hurwitz) zeta functions $\zeta(z_1,
\ldots,z_k)$  (resp. $\zeta(z_1, \ldots, z_k;\, v_1, \ldots , v_k)$) given by
Theorem \ref{thm:mzext} satisfy stuffle relations and are
meromorphic in all variables with poles on a countable family of hyperplanes $z_1+\cdots +z_j\in
]-\infty, j]\cap \Z$, $j$ varying from $1$ to $k$.
\end{cor}
 \begin{rk} More can be said about the pole structure for multiple zeta
   functions of depth $2$ as recalled in the introduction and as we shall see
   in Theorem \ref{thm:recurrence}. 
 
\end{rk}
When the nested sums do not converge, one does not expect the stuffle relations
to hold  in general neither for cut-off nested sums, nor for finite parts of meromorphic perturbations of these sums, because of extra terms which occur when taking finite parts that might spoil
the stuffle relations. However, just as for nested integrals, one can use instead a renormalisation procedure
which takes care of these extra terms and thereby  build  renormalised nested sums  which do obey the required stuffle
relations. For this purpose,  we first  extend a holomorphic regularisation
${\cal R}$ on ${\cal P}^{*,0}$ to one on the tensor algebra  ${\cal T}({\cal
  P}^{*,0})$ which is
compatible with the stuffle product (with respect to a product $\bullet$ which
will be $\mp$ the ordinary product of functions). This can be carried out by twisting by
the Hoffman isomorphism the
previously extended regularisation $\tilde {\cal R}$ compatible with the
shuffle product. Applying this isomorphism to any subalgebra ${\cal A}$ of
${\cal P}^{*,0}$ provides a
regularisation $\tilde{\cal R}^\star$  on $ {\cal T}({\cal A})$
compatible with the stuffle product $\star_\bullet$, just as $\tilde{\cal R}$ was with the 
  shuffle product $\shu$ (see (\ref{eq:regshu})):
\begin{lem}\label{lem:rstar}
A regularisation ${\cal
      R}$ on any subalgebra ${\cal A}$ of ${\cal P}^{*,0}$ extends to one  on $ {\cal T}({\cal A})$ defined by
$$\tilde {\cal R}^\star:=\exp\circ \tilde {\cal R}\circ\log$$
  which is compatible with $\star_\bullet$:
\begin{equation}\label{eq:regstu}\tilde {\cal R}^\star(\sigma\star_\bullet\tau)= \tilde {\cal R}^\star(\sigma)\star_\bullet
\tilde {\cal R}^\star(\tau)\quad \forall \sigma, \tau \in  {\cal T}\left({\cal P}^{*,0}\right). \end{equation}
\end{lem}
{\bf Proof:}
This is an immediate consequence of Theorem \ref{th:hoffman}.
\endsquare
\begin{prop}:\label{prop:meroChensums2} Let $W$ be an open neighbourhood of
  $0$ in $\C$, let ${\cal A}$ be a subalgebra
    of ${\cal P}^{*,0}$, and let  ${\cal R}: \tau \mapsto
  \{\tau(z)\}_{z\in W}$ be a holomorphic regularisation on  ${\cal A}$  such that  the perturbed symbol  $\sigma(z)$ of any  symbol $\sigma$ has
 order $\alpha(z) $  with $\mop{Re }\alpha^\prime(z)<0$ for any $z\in W$. We
 suppose that for any real $t$ there is a $z\in W$ such that for any symbol
 $\sigma$ of real order $t$, the symbol $\sigma(z)$ has real order $<-1$. Then the map
\begin{eqnarray*}
{\cal A}&\to & {\cal M}(W)\\
\sigma&\mapsto& \cutoffsum\sigma(z)
\end{eqnarray*}
extends to  multiplicative maps
\begin{eqnarray*}
\Psi^{{\cal R}}(\hbox{ resp. }{\Psi'}^{{\cal R}}):\left( {\cal T}\left({\cal A}\right), \star_\bullet\right)&\to& {\cal M}(W)\\
 \sigma &\mapsto & 
\cutoffsum^{\smop{Chen}}_<\tilde{\cal R}^\star(\sigma)(z)\hskip 6mm\big(\hbox{
   resp. } \cutoffsum^{\smop{Chen}}_{\leq}\tilde{\cal R}^\star(\sigma)(z)\big) \\
\end{eqnarray*}
where $\bullet$ stands for the ordinary product $\cdot$ in the first case, and
stands for the opposite of the ordinary product in the second case. In other
words, $\Psi^{\cal R}$ and ${\Psi'}^{\cal R}$ satisfy the stuffle relations: 
\begin{eqnarray*}
\Psi^{\cal R}(\sigma\star_{+\displaystyle\cdot} \tau)&=&\Psi^{\cal R} (\sigma)\,
 \cdot\, \Psi^{\cal R} (\tau),\\
{\Psi'}^{\cal R}(\sigma\star_{-\displaystyle\cdot} \tau)&=&{\Psi'}^{\cal R} (\sigma)\,
 \cdot\, {\Psi'}^{\cal R} (\tau)
\end{eqnarray*}
which hold as equalities of meromorphic functions. 
\end{prop}
{\bf Proof:} The proof is carried out along the same lines as for continuous
integrals. 
By construction, $\tilde{\cal R}^\star\left(\sigma_1\otimes \cdots \otimes
  \sigma_{k}\right)(z) $ is a finite linear combination of bullet-products and tensor
products of the ${\cal
  R}(\sigma_i)(z)$'s. For example:
\begin{equation}\label{eq:Rtildestar}
\tilde{\cal R}^\star(\sigma_1\otimes\sigma_2)(z)=\sigma_1(z)\otimes\sigma_2(z)
-\frac 12(\sigma_1\bullet\sigma_2)(z)+\frac 12\sigma_1(z)\bullet\sigma_2(z).
\end{equation}
It therefore follows from  Theorem \ref{thm:meroChensums}  that the   cut-off nested sums
$\Psi^{\cal R}(\sigma)(z)$ and ${\Psi'}^{\cal R}(\sigma)(z)$ define meromorphic functions. Moreover, the compatibility of $\tilde{\cal
  R}^\star$ with the stuffle product $\star_\bullet$ implies that
 \begin{eqnarray*}
\Psi^{{\cal R}}(\sigma\star_\bullet\tau)(z)&=& \cutoffsum^{\smop{Chen}}\tilde{\cal R}^\star(\sigma\star_\bullet\tau)(z)\\
&=&  \cutoffsum^{\smop{Chen}}\tilde{\cal R}^\star(\sigma)(z)\star_\bullet\tilde{\cal
  R}^\star(\tau)(z)\\
&=&  \cutoffsum^{\smop{Chen}}\tilde{\cal R}^\star(\sigma)(z)\cdot \cutoffsum^{\smop{Chen}}\tilde{\cal R}^\star(\tau)(z),\\
&=&\Psi^{{\cal R}}(\sigma)(z)\star_\bullet\Psi^{{\cal R}}(\tau)(z),
\end{eqnarray*}
 holds on a common domain  of convergence of these sums as a result of the
 usual stuffle relations for convergent nested sums (see  (\ref{eq:stuffle})). This
identity of holomorphic functions on some common domain (which is non-empty
due to the hypothesis on $W$ and the regularisation procedure) then  extends by
analytic continuation to an
identity of meromorphic functions. 
  \endsquare\\ \\
Birkhoff factorisation provides a way to extract finite parts which still
obey stuffle relations:
\begin{thm} \label{thm:renChensums}
Given a subalgebra ${\cal A}\subset  {\cal P}^{*,0}$ and a  holomorphic regularisation ${\cal R}: \sigma
\mapsto \{\sigma(z)\}_{z\in W}$ on ${\cal A}$  which satisfies the same
assumption as in Proposition \ref{prop:meroChensums2},  
the map
\begin{eqnarray*} {\cal A}&\to &\C\\
\sigma&\mapsto &\mopl{fp }_{z=0} \cutoffsum\sigma(z)
\end{eqnarray*}
extends to a character:
\begin{eqnarray*}\psi^{{\cal R}}(\hbox{ resp. }{\psi'}^{{\cal R}}):\left({\cal
    T}\left({\cal A}\right), \star_\bullet\right),
&\longrightarrow &\C\\ 
\sigma&\mapsto & \Psi_{+}^{{\cal R}}(\sigma)(0)\hskip 6mm\big(\hbox{ resp. }
{\Psi'}_{+}^{{\cal R}}(\sigma)(0)\big)
\end{eqnarray*}
 defined from  the Birkhoff decomposition of $\Psi^{{\cal
     R}}$
(resp. ${\Psi'}^{{\cal R}}$). It coincides with the ordinary nested sums
$\sum^{\smop{Chen}}_{<}$  (resp. $\sum_{\leq}^{\smop{Chen}}$) on elements
$\sigma=\sigma_1\otimes\cdots\otimes\sigma_k$ with $\sigma_j\in{\cal
  P}^{\alpha_j,0}$ when
 $\mop{Re }(\alpha_1+\cdots +\alpha_m)<-m$ for any $m\in\{1,\ldots ,k\}$. Here $\bullet$ stands for the product $\pm\cdot$ as in Proposition \ref{prop:meroChensums2}.
\end{thm}
{\bf Proof:} 
Let us prove  the strict inequality case: the Birkhoff decomposition \cite{CK}, \cite{M} for the minimal subtraction scheme reads:
$$\Psi^{{\cal R}}= \left(\Psi_-^{{\cal R}}\right)^{* -1}\, *\, \Psi_+^{{\cal
    R}}$$
with $\bullet=\cdot$ (see Theorem \ref{thm:renChenint}). Since $\Psi_+^{{\cal
    R}}$ is multiplicative  on $\left({\cal T}\left({\cal A}\right), \star_\bullet\right)$, it  obeys the
stuffle relation:
\begin{equation}\label{eq:stufflehol}\Psi_+^{{\cal
    R}}\left(\sigma \star_\bullet \tau\right)(z)= \Psi_+^{{\cal
    R}}(\sigma)(z) \, \Psi_+^{{\cal
    R}} (\tau)(z)
\end{equation}
which holds as an equality of meromorphic functions holomorphic at $z=0$.
 Setting
$$\psi^{{\cal R}}:= \Psi_+^{{\cal
    R}}(0),$$
and applying the stuffle relations (\ref{eq:stufflehol}) at $z=0$  yields 
$$\psi^{{\cal R}}\left(\sigma\,\star_\bullet\,
 \tau\right)= \psi^{{\cal
    R}}(\sigma)\, \psi^{{\cal R}}(\tau).$$
    The tensor products $\sigma=\sigma_1\otimes\cdots\otimes\sigma_k$ with $\sigma_j\in{\cal
  P}^{\alpha_j,0}$ where
 $\mop{Re }(\alpha_1+\cdots +\alpha_m)<-m$ for any $m\in\{1,\ldots ,k\}$ span
 a right co-ideal ${\cal J}$ of ${\cal T}({\cal A})$, namely $\Delta({\cal
   J})\subset {\cal J}\otimes{\cal T}({\cal A})$. The restriction of
 $\Psi^{{\cal R}}$ to this right co-ideal takes values into functions which are
 holomorphic at $z=0$. By construction we then have:
$$\psi^{{\cal R}}\left(\sigma\right)=\Psi_+^{{\cal
    R}}(\sigma)(0)= \Psi^{{\cal
    R}}(\sigma)(0)= \cutoffsum^{\smop{Chen}}_{\leq} \sigma$$
by a similar argument as in the proof of Theorem \ref{thm:renChenint}. The
weak inequality case can be derived similarly  setting $\bullet=-\cdot$.\endsquare\\ \\
On the grounds of this result we set the following definition:
 \begin{defn}\label{defn:Rregmultiple zeta} For any $\sigma\in{\cal T}({\cal A})$,
   the {\rm renormalised nested sums\/} of $\sigma$ (in both weak and strict
   inequality versions) are defined by:
$$\cutoffsum^{\smop{Chen}, {\cal R}}_{<}
 \sigma:= 
 {\psi}^{{\cal R}}(\sigma), \hskip 12mm \cutoffsum^{\smop{Chen}, {\cal R}}_{\leq}
 \sigma:= 
 {\psi'}^{{\cal R}}(\sigma).$$
 \end{defn}
\begin{rk}
There are other multiplicative maps from $\big({\cal T}({\cal
  A}),\star_\bullet\big)$ to $\C$ defined by $\phi_r^{{\cal R}}\circ\log$ with the
notations of Theorem \ref{thm:renChenint}. Those maps are not related to nested sums.
\end{rk}
\vfill \eject  \noindent
\section{RENORMALISED MULTIPLE ZETA FUNCTIONS} \label{sect:entiersnegatifs}
 Recall that multiple Hurwitz zeta functions $\zeta(s_1, \ldots, s_k; v_1,
 \dots, v_k)$ and
 $\overline\zeta(s_1, \ldots, s_k; \, v_1, \ldots, v_k)$ converge whenever $\mop{Re }(s_1+\cdots +
  s_m)>m$ for any $m\in\{1,\ldots ,k\}$, in which case they obey stuffle relations. In this
   section, we implement the  renormalisation procedure described in section 10 to extend them
  to other complex values of $s_i$ while preserving the stuffle relations. We
  then specialise to nonnegative integer arguments for which
we prove rationality of the renormalised multiple zeta values.
\subsection{Stuffle relations for renormalised  multiple zeta functions }\label{sect:stufflezeta}
Let  $\widetilde{\cal A}$ be the subalgebra of ${\cal P}^{*,0}$ generated by
the continuous functions with support inside the interval $]0,1[$ and the set
$$\{f\in {\cal P}^{*,0},\, \exists v\in \R_+, \exists s\in\C,\, f(\xi)=
(\xi+v)^{-s}\hbox{ when }\xi\ge 1\}.$$
Consider the ideal ${\cal N}$ of $\widetilde{\cal A}$ of continuous functions
with support included in the interval $]0,1[$. The quotient algebra ${\cal
  A}=\widetilde {\cal A}/{\cal N}$ is then generated by the elements
$\sigma_{s}^v$, where $\sigma_{s}^v$ is the class of any $f\in\widetilde {\cal A}$
such that $f(\xi)=(\xi+v)^{-s}$ for $\xi\ge 1$. Notice that for any
$v\in\R_+$ the subspace ${\cal A}_v$ of $\cal A$ generated by
$\{\sigma_{s}^v,\,s\in\C\}$ is a subalgebra of ${\cal A}$. We choose the product
$\bullet$ as the opposite of the ordinary product, so that we have:
$$\sigma\bullet\sigma^\prime=-\sigma\sigma^\prime\quad\forall (\sigma,
\sigma^\prime)\in {\cal A}^2; \quad {\rm resp.}\quad
\sigma_s^v\bullet\sigma_{s^\prime}^{v}=-\sigma_{s+s^\prime}^v\quad \forall (\sigma_s^v,
\sigma_{s^\prime}^v)\in {\cal A}_v^2.$$
Let ${\cal W}$ be the $\C$-vector space freely spanned by sequences
$(u_1,\ldots,u_k)$ of real numbers. Let us define the stuffle product on
${\cal W}$
by:
\begin{equation}\label{eq:stuffleR}
(u_1,\ldots ,u_{k})\star(u_{k+1},\ldots ,u_{k+l})=
\sum_{0\le r\le \smop{min}(k,l)}(-1)^r\sum_{\pi\in\smop{mix
    sh}(k,l;r)}(u^\pi_1,\ldots ,u^\pi_{k+l-r}),
\end{equation}
with:
$$u^\pi_j=\sum_{i\in\{1,\ldots,k+l\},\pi(i)=j}u_i$$
(the sum above contains only one or two terms).
Define a map $u\mapsto\sigma_u^v$ from $W$ to ${\cal T}({\cal A}_v)$ by:
$$\sigma_{(u_1,\ldots ,u_k)}^v:=\sigma_{u_1}^v\otimes\cdots\otimes\sigma_{u_k}^v.$$
Then
$$\sigma_u^v\star_\bullet \sigma_{u^\prime}^v= \sigma_{u\star
  u^\prime}^v.$$
The same holds with $\bullet=\cdot$ provided we drop the signs $(-1)^r$ in
equation (\ref{eq:stuffleR}) defining the stuffle product on ${\cal W}$.
\begin{defn}\label{mzrenorm}
Let $W$ be an open neighbourhood of $0$ in $\C$. Let ${\cal R}:\sigma\mapsto
\{\sigma(z)\}_{z\in W}$ be a holomorphic regularisation procedure on $\widetilde{\cal
  A}$ such that the order condition of Proposition \ref{prop:meroChensums2} is satisfied, and which
respects the ideal ${\cal N}$, hence giving rise to a regularisation procedure
on the quotient ${\cal A}$. The renormalised multiple Hurwitz zeta functions (with
respect to ${\cal R}$) are defined by:
\begin{eqnarray*}
 \zeta^{{\cal R}}(s_1, \ldots, s_k; \, v_1, \ldots, v_k)&:=&{\psi}^{{\cal
    R}}(\sigma_{s_1}^{v_1}\otimes \cdots \otimes\sigma_{s_k}^{ v_k}),\\
\hskip 12mm
\overline\zeta^{{\cal R}}(s_1, \ldots, s_k; \, v_1, \ldots, v_k
)&:=&{\psi'}^{{\cal R}}(\sigma_{s_1}^{v_1}\otimes \cdots \otimes\sigma_{s_k}^{v_k} ).
\end{eqnarray*}
For $v_1=\cdots= v_k=v$ we simply set 
$$
 \zeta^{{\cal R}}(s_1, \ldots, s_k;\, v):= {\psi}^{{\cal R}}(\sigma_{s_1, \ldots,
  s_k}^v);\quad
 \overline \zeta^{{\cal R}}(s_1, \ldots, s_k;\, v):= 
{\psi^\prime}^{{\cal R}}(\sigma_{s_1, \ldots, s_k}^v).$$
\end{defn}
By Theorem \ref{thm:renChensums} applied to $\widetilde{\cal A}$  defined at the
beginning of the section, combined with the fact that the nested sums and the regularisation
are well-defined on the quotient ${\cal A}$, we know that ${\psi}^{{\cal
    R}}$ is compatible with the stuffle product:
\begin{equation*}
\hskip -4mm {\psi}^{{\cal
    R}}\left((\sigma_{s_1}^{v_1}\otimes \cdots \otimes\sigma_{s_k}^{v_k})\star
  (\sigma_{s^\prime_1}^{v^\prime_1}\otimes \cdots\otimes \sigma_{s^\prime_{k^\prime}}^
    {v^\prime_{k^\prime}})\right)={\psi}^{{\cal
    R}}\left(\sigma_{s_1}^{v_1}\otimes \cdots\otimes \sigma_{s_k}^{v_k}\right)\, 
{\psi}^{{\cal
    R}}\left(  \sigma_{s^\prime_1}^{v^\prime_1}\otimes \cdots \otimes\sigma_{s^\prime_{k^\prime}}^
    {v^\prime_{k^\prime}}\right)
\end{equation*}
and similarly for ${\psi^\prime}^{{\cal
    R}}$. When setting the $v_i$'s equal to
some fixed $v$,  Theorem  \ref{thm:renChensums} applied to $\widetilde{\cal A}_v$ leads
to the following stuffle properties of families of Hurwitz multiple zeta
functions which contain the ordinary multiple zeta functions obtained by setting $v=0$.
\begin{thm}\label{thm:murenChensum}
\begin{enumerate}
\item Renormalised Hurwitz multiple zeta values have the following properties:
\begin{equation}\label{stuffle-rel}
\overline\zeta^{{\cal R}}(u\star u^\prime; \, v) =
  \overline\zeta^{{\cal R}}(u;v) \, \overline\zeta^{{\cal R}}(u^\prime;v) 
\end{equation}
when the stuffle product $\star$ is defined by (\ref{eq:stuffleR}), and:
\begin{equation}\label{stuffle-relbis}
\zeta^{{\cal R}}(u\star u^\prime; v) =
  \zeta^{{\cal R}}(u;v) \, \zeta^{{\cal R}}(u^\prime; v) 
\end{equation}
when the stuffle product $\star$ is defined by (\ref{eq:stuffleR}) with signs
$(-1)^r$ removed.
\item for any positive
integer $k$, and for $(s_1,\ldots ,s_k)\in\C^k$ such that $\mop{Re }(s_{1}+\cdots +s_m)>m$ for any
$m\in\{1,\ldots ,k\}$, the renormalised values are independent of ${\cal R}$
and can be written as ordinary nested sums:
\begin{eqnarray*}
\overline\zeta^{{\cal R}}(s_1, \ldots, s_k;\, v)&=& \sum_{0< s_1\le\cdots\le 
    n_k} (n_1+v)^{-s_{1}}\cdots (n_k+v)
    ^{-s_k}
=\overline\zeta(s_{1}, \ldots, s_{k}; \,v),\\
\zeta^{{\cal R}}(s_1, \ldots, s_k; \, v)&=& \sum_{0< s_1<\cdots<
    n_k} (n_1+v)^{-s_{1}}\cdots (n_k+v)
    ^{-s_k}
=\zeta(s_{1}, \ldots, s_{k}; \,v).
\end{eqnarray*}
 \end{enumerate}
\end{thm}
This provides the  uniqueness of the
  extension of Riemann multiple zeta functions to regularised multiple zeta
  functions satisfying stuffle relations, once the value  $
   \theta$ at the argument $1$  is imposed; see  \cite{H},
  \cite{W}, \cite{Z}. Indeed, the expressions $\overline\zeta^{{\cal R}} (s_1,
\ldots, s_k)$ and $\zeta^{{\cal R}} (s_1,
\ldots, s_k)$ converge
whenever $s_1>1$  since by assumption, all the $s_i$ are no smaller than $1$. 
When they converge, they obey the stuffle relations (\ref{stuffle-rel}) and
(\ref{stuffle-relbis}) respectively. The   uniqueness of the extension to the case $s_1=1$ then 
follows by induction on the length $k$  from the  stuffle relations
(\ref{stuffle-rel}) which ``push'' the leading term $s_1=1$ whenever it arises, away from the
first position and therefore expresses divergent expressions in terms of
convergent expressions. Given a real number  $\mu$, we can consider
  the holomorphic 
  regularisation (cfr. (\ref{eq:regmu}))
\begin{equation}\label{eq:Rmu}\xi\mapsto {\cal R}_\mu (\tau)(\xi):= \tau(\xi) \, \vert \xi
  \vert^{-\gamma_\mu(z)}\end{equation} which is well defined on ${\cal
  A}$. For example we can consider $\gamma_\mu(z)= z+ \mu z^  2$. In this case
the parameter $\mu$ plays the role of the constant $\theta$.  
\subsection{Rationality of  (Hurwitz) multiple zeta values at nonpositive arguments}
The following theorem  provides refined meromorphicity and  holomorphicity results for multiple
zeta functions which go beyond  the
meromorphicity results of Corollary   \ref{cor:meroChensums}. 
\begin{thm}\label{thm:recurrence}
Let $\gamma$ be a holomorphic function in an open neighbourhood $W$ on $0$ in
$\C$ such that $\gamma(0)=0$ and $\gamma'(0)=1$. Let $l\in\N-\{0\}$. Let
$(\beta_1,\ldots,\beta_l)$ be functions $\beta_j(z)=b_j-c_j\gamma(z)$ such that $\mop{Re }c_j>0$ for
$j=1,\ldots l$. 
\begin{enumerate}
\item  For any $v_1\geq 0, \ldots, v_l \geq 0$ the maps 
$$(z_1, \ldots, z_l)\longmapsto \cutoffsum_{1\le n_l<\cdots< n_1} (n_1+v_1)^{\beta_1(z_1)}\cdots
(n_l+v_l)^{\beta_l(z_l)}$$ are meromorphic with poles on hyperplanes 
$\sum_{i=1}^m c_iz_i \in -\sum_{i=1}^m b_i -m+\N$, $m$ varying in $\{1,
\ldots, l\}$ so that poles at $(z_1, \ldots, z_l)=(0, \ldots,0)$ can arise
whenever  $\sum_{i=1}^m b_i  \in [-m; +\infty[\cap \Z$. When $l=2$ and  $v_1=\cdots=
v_l=v$ for some non negative real number $v$, poles
actually only arise when  $b_1=-1$ or
  $b_1+b_2\in \{-2,-1, 0, 2,4, 6, \ldots\}$.
\item Let us assume that
$b_1,\ldots ,b_{l-1}\in\N$ and $b_l\in\Z$.
\begin{itemize}
\item If $b_l\in\N$, for any non negative real number $v$,  the map
$$z\longmapsto \cutoffsum_{1\le n_l<\cdots< n_1} (n_1+v)^{\beta_1(z)}\cdots
(n_l+v)^{\beta_l(z)}$$
is holomorphic around $z=0$. Its value at $z=0$ is a polynomial expression in $v$  with  coefficients given by
  rational functions in  the $c_j$'s   that 
do not depend on the choice of $\gamma$.
\item If $b_l\le -1$, the  above map is meromorphic with a  possible pole of order at most
  $1$ at $z=0$. The residue at $z=0$ is a polynomial expression in $v$  with  coefficients given by
  rational functions in  the $c_j$'s that do not depend on the choice of
  $\gamma$.
\end{itemize}
\end{enumerate}
\end{thm}
\begin{rk} Our proof, which relies on the algebra structure of the set
  $\{\xi\mapsto (\xi+v)^{a},\quad \xi >0, a\in \N\}$ does not  a priori extend to show  holomorphicity and rationality for all non
  negative values $v_1, \ldots, v_k$ which  therefore remains an open question.
\end{rk}
{\bf Proof:} By Theorem   \ref{thm:meroChensums}  applied to $\sigma_i(\xi)=
(\xi+v_i)^{-s_i}\, \chi(x)$ where $\chi$ is a smooth cut-off function which
vanishes at zero and is identically one outside the unit interval,  the maps 
$$(z_1, \ldots, z_l)\longmapsto \cutoffsum_{1\le n_l<\cdots< n_1} (n_1+v_1)^{\beta_1(z_1)}\cdots
(n_l+v_l)^{\beta_l(z_l)}$$ are meromorphic with poles on hyperplanes 
$\sum_{i=1}^m c_iz_i \in -\sum_{i=1}^m b_i -m+\N$ with $m$ varying in $\{1,
\ldots, l\}$ from which we infer the first part of the proposition concerning the pole
structure for a general $l$ and general $b_i$'s which generalise the results
of  Corollary   \ref{cor:meroChensums}. In the case $l=2$, more can be said on the structure of the
poles; we postpone the proof, leaving  this case for later. Since ${\rm Re}(c_i)<0$ for any $i\in\{1, \ldots, l\}$ we can
set $z_i= z$, from which the meromorphicity  of the map $z\mapsto \cutoffsum_{1\le n_l<\cdots< n_1} (n_1+v_1)^{\beta_1(z)}\cdots
(n_l+v_l)^{\beta_l(z)}$ follows. We actually recover the meromorphicity of this map  by the
inductive proof below. To prove Part 2 of the theorem, we indeed proceed by induction on the depth $l$ and  set in
  order to lighten the notations:
\begin{equation}\label{eq:abbrev}
\sigma_j(z)(\xi):=(\xi+v)^{\beta_j(z)}, \hbox to 15mm{\hfill with }\beta_j(z)=b_j-c_j\gamma(z),
\end{equation}
for $\xi\ge 1$.  
\begin{enumerate}
\item In the case $l=1$ we have $\cutoffsum_{1\le
  n_1}(n_1+v)^{\beta_1(z)}=\zeta(-\beta_1(z);\, v)$ where, for any $v\ge 0$, we have set
$\zeta(s; \,v)=
\sum_{n=1}^\infty (n+v)^{-s}$ which relates to the usual Hurwitz zeta function
$\zeta(s, v):= \sum_{n=0}^\infty (n+v)^{-s}$ by 
$\zeta(s; \, v)= \zeta(s, v+1)$.  The first step of the induction follows from
Proposition \ref{prop:sumlogn} which gives  \begin{equation}
\mopl{Res }_{z=0}\zeta(-\beta_1(z))=\delta_{b_1+1}\, \frac{1}{c_1}; \quad \zeta(-k;v)=
-\frac{B_{k+1}(1+v)}{k+1}\quad{\rm if}\quad b_1\in \N.
\end{equation} 
\item Now suppose $l\ge 2$. It is useful to observe that each
summation raises the order of the corresponding interpolated symbols by $1$,
so that to pick out the finite part as $N\to \infty$, which amounts to
extracting the homogeneous part of the symbol of degree $0$, we need
to take into consideration homogeneous parts of the symbol of negative degree as we go into the depths of the nested sums,
i.e. of degree $\geq
-J$ after $J$ summations. For a symbol $\sigma$ of order $a$ and for any real
$\omega$ we call $\sigma_{\geq \omega}$ the (finite) sum of terms of real order
$\ge \omega$ in  the
asymptotic expansion  
of the symbol $\sigma$. We have for any $z$ in a small neighbourhood of $0$:
\begin{eqnarray}\label{eq:practicalmultiple zeta}
&{}&\left(\sum_{0<n_l<\cdots< n_1< N}
(n_1+v)^{\beta_1(z)}\cdots (n_l+v)^{\beta_l(z)}\right)_{ \geq 0}\nonumber\\
&\hskip -40mm=&\hskip -20mm
\left(\sum_{n_1=1}^{N-1} (n_1+v)^{\beta_1(z)} \left(\sum_{0<n_l<\cdots< n_1}
(n_2+v)^{\beta_2(z)}\cdots (n_l+v)^{\beta_l(z)}\right)_{\geq -b_1-1}\right)_{\geq  0}\\
&\hskip -40mm=&\hskip -20mm
\left(\sum_{n_1=0}^{N-1} (n_1+v)^{\beta_1(z)} \left(\sum_{n_2=1}^{n_1-1}
(n_2+v)^{\beta_2(z)}\cdots\left(\sum_{n_l=1}^{n_{l-1}-1}
(n_l+v)^{\beta_l(z)}\right)_{\geq -b_{l-1}-b_{l-2} \cdots -b_1-l+1}\cdots\right)_{\geq -b_1 -1}\right)_{\geq 0}. \nonumber
\end{eqnarray}
\item We extend the
notation $[a]_j$ to $j=0$ and $j=-1$ by setting:
\begin{equation}
[a]_0:=1,\hskip 20mm [a]_{-1}:=\frac{1}{a+1}.
\end{equation}
Let us use (\ref{eq:cutoffsumN}) to expand the interpolating symbol
  $\overline Q\big(\sigma_l(z)\big)$ for the
  deepest sum $Q\big(\sigma_{l}(z)\big)(N)= \sum_{0<n_l< N}\sigma_{l}(z)(n_l)$, writing:
 \begin{eqnarray}\label{eq:practical1}
\overline Q\big(\sigma_{l}(z)\big)(\xi)
&=&\sum_{j=0}^{2J_l} B_j\,\frac{ [\beta_l(z)]_{j-1}}{j!}\left(
  (\xi+v)^{\beta_l(z)-j+1}-1\right)\nonumber\\
&\hbox to 3mm{}&+\frac{ [\beta_l(z)]_ {2J_l+1}}{(2J_l+1)!} \int_1^\xi \overline{B_{2J_l+1}}
(y)\,(y+v)^{\beta_l(z)-2J_l-1}\, dy,
\end{eqnarray}
choosing $J_l \ge\left[\frac 12\left(\sum_{j=1}^l b_j+l\right)\right]$ we
have $b_l-2J_l-1 < -b_1-\cdots -b_{l-1} -l+1$  so that  we  do
not miss any terms in (\ref{eq:practicalmultiple zeta}).
\item  The $l=2$ case is instructive and worth being treated separately since
  it presents some specificities; in that case
  $J_2\geq \left[\frac 12\,(b_1+b_2)\right]+1$ so that if $b_1+b_2$ is
    even, we take $J_2=\frac 12\,(b_1+b_2)+1$ and  if  $b_1+b_2$ is
    odd, we take $J_2=\frac 12\,(b_1+b_2)+\frac 32.$
 By
(\ref{eq:practical1}), we have 
\begin{eqnarray*}
\hskip -12mm \altcutoffsum_{<}^{\smop{Chen}}\sigma_1(z_1)\otimes\sigma_2(z_2)&=&\cutoffsum_0^\infty
\sigma_1(z_1)\overline Q\big(\sigma_{2}(z_2)\big)\\
&=& 
\sum_{j=0}^{2J_2} B_j\,\frac{ [\beta_2(z_2)]_{j-1}}{j!}\left(
  \cutoffsum_0^\infty (n+v)^{\beta_1(z_1)}\,(n+v)^{\beta_2(z_2)-j+1}-\cutoffsum_0^\infty (n+v)^{\beta_1(z_1)}\right)\nonumber\\
&\hbox to 3mm{}&-
\frac{ [\beta_2(z_2)]_ {2J_2+1}}{(2J_2+1)!}\cutoffsum_0^\infty\left(
(n+v)^{\beta_1(z_1)}\,  \int_1^n \overline{B_{2J_l+1}}
(y)\,(y+v)^{\beta_2(z_2)-2J_2-1} \, dy\right).\nonumber\\
&=& 
\sum_{j=0}^{2J_2} B_j\,\frac{ [\beta_2(z_2)]_{j-1}}{j!}\left(
\zeta(-\beta_1(z_1)-\beta_2(z_2)+j-1;\, v)  -\zeta(-\beta_1(z_1);\, v\right)\nonumber\\
&\hbox to 3mm{}&+
\frac{ [\beta_2(z_2)]_ {2J_2+1}}{(2J_2+1)!}\cutoffsum_0^\infty\left(
(n+v)^{\beta_1(z_1)}\,  \int_1^n \overline{B_{2J_l+1}}
(y)\,(y+v)^{\beta_2(z_2)-2J_2-1} \, dy\right)\nonumber\\
\end{eqnarray*}
Since $J_2\geq\frac 12\,(b_1+b_2)+1$, it follows that $b_1+b_2-2J_2<-1$ so that the
  last term is absolutely convergent and hence holomorphic. Moreover, we
  observe that $[b_2]_{2J_2+1}=0$ whenever $b_2$ is an integer since $b_2-2J_2\leq -b_1-2< 0$, so that the
  last term actually vanishes at $(z_1,z_2)=(0, 0)$ in that case. The map 
$(z_1, z_2)\mapsto
\altcutoffsum_{<}^{\smop{Chen}}\sigma_1(z_1)\otimes\sigma_2(z_2)$ is therefore
meromorphic with poles on a finite number of hyperplanes $b_1-c_1\,z_1=-1$ and $c_1
\,z_1+c_2\,z_2=b_1+b_2-j+2, \quad j\in \{0,1, 2, 4, \ldots, 2J_2\}$ since
Bernoulli numbers
$B_j$ vanish for odd $j$ larger than 2 and  the only pole of
  the zeta function is $1$. Poles at zero therefore only arise if $b_1=-1$ or
  $b_1+b_2\in \{-2,-1, 0, 2,4, 6, \ldots\}$, thus confirming known results
 \cite{AET}, see also \cite{Zh} and \cite{G} Theorem 2.25.\\
The map  is actually holomorphic at zero for any value of $(b_1, b_2)$ along any  hyperplane $z_1=\lambda\, z_2$ with
$\lambda>0$. Indeed, setting  $z_1=\lambda \,z_2$ in the above expression gives
rise to (simple) poles in $z$ which are compensated by terms
$[\beta_2(z)]_{j-1}$   since these involve a factor $z$ as a consquence of the inequality
$b_2-j+2=-b_1\leq 0$ resulting from $b_1+b_2+2=j$. Combining the above results
moreover shows that the evaluation at $z=0$
gives rise to a rational number whenever $b_1$ and $b_2$ are integers as a
result of the rationality of Bernoulli numbers. 
\item Let us now assume $l\geq 2$, thus including the case $l=2$
  even though it was treated separately. 
We write:
$$\altcutoffsum_{<}^{\smop{Chen}}\sigma_1(z)\otimes\cdots\otimes
\sigma_l(z)=\altcutoffsum\sigma_1(z)\otimes \cdots
\otimes\sigma_{l-2}(z_{l-2})\otimes 
  \sigma_{l-1}(z_{l-1})\overline Q(\sigma_l)(z_l)$$
which is a finite linear combination of nested sums 
 $\altcutoffsum_<^{\smop{Chen}}\sigma_1(z)\otimes \cdots
\otimes\sigma_{l-2}(z)\otimes 
  \sigma_{l-1}(z)\sigma_{l,j}(z)$
 of depth $l-1$: with the above choice of $J_l$, $b_l-2J_l\leq -\sum_{i=1}^l
 b_i -l+1<0$ so that  by similar argument as in the
 $l=2$ case, one checks that the remainder term $\sigma_1(z)\otimes \cdots
\otimes\sigma_{l-2}(z)\otimes 
  \sigma_{l-1}(z)\rho_{\beta_l}(z)$ with:
$$ \rho_{\beta_l}(z)(\xi):=-\frac{ [\beta_l(z_l)]_{2J_l+1}}{(2J_l+1)!} \int_1^\xi \overline{B_{2J_l+1}}
(y)\,(y+v)^{\beta_l(z_l)-2J_l-1}\, dy$$
 contributes by a holomorphic term which vanishes at $z=0$. Applying
(\ref{eq:practical1}) we have:
\begin{eqnarray}\label{recurrence2}
&{}&\cutoffsum_{1\le n_l<\cdots< n_l} (n_1+v)^{\beta_1(z)}\cdots
(n_l+v)^{\beta_l(z)}\nonumber\\
&=&\sum_{j=0}^{2J_l}H_l^j(z_l)\cutoffsum_{1\le n_{l-1}<\cdots< n_{1}} (n_1+v)^{\beta_1(z_1)}\cdots
(n_{l-2}+v)^{\beta_{l-2}(z)}(n_{l-1}+v)^{\kappa_{l-1}^j(z_{l-1},z_l)}\nonumber\\
&+& K_l(z_l)\cutoffsum_{1\le n_{l-1}<\cdots< n_{1}} (n_1+v)^{\beta_1(z_1)}\cdots
(n_{l-1}+v)^{\beta_{l-1}(z_{l-1})}
\end{eqnarray}
with
\begin{equation}\label{eq:kappalj}
\kappa_{l-1}^j(z_{l-1}, z_l)=\beta_{l-1}(z_{l-1})+\beta_l(z_l)+1-j,
\end{equation}
\begin{equation}\label{eq:hlj}
H_l^j(z_l)= \frac{B_j}{j!}[\beta_l(z_l)]_{j-1},
\end{equation}
\begin{equation}\label{eq:kl}
K_l(z_l)=
-\sum_{j=0}^{2J_l}\frac{B_j}{j!}[\beta_l(z_l)]_{j-1}+\rho_{\beta_l(z_l)}.
\end{equation}
The function $z\mapsto K_l(z)$ is holomorphic at $z=0$ unless
$b_l=-1$. The functions $z\mapsto H_l^j(z)$ are holomorphic in $z=0$ unless $j=0$ and
$b_l=-1$. If $b_1=-1$, the functions $z\mapsto H_l^0(z)$ and $z\mapsto K_l(z)$
have a simple pole in
$z=0$ with residue $-\frac{1}{c_l}$ and $\frac{1}{c_l}$ respectively.
\item We are now ready to carry out the proof of the second item in the
  proposition  by induction on the depth $l$. Suppose that the proposition
is verified up to depth $l-1$. The left-hand side of \eqref{recurrence2} is
then meromorphic at $z=0$. For any meromorphic function $f$ let us call $\mopl{Res$^j$ }_{z=0}(f)$ the
coefficient of $z^{-j}$ in the Laurent expansion of $f$ at $z=0$. In particular
$\mopl{Res$^0$ }_{z=0}(f)$ stands for the finite part, $\mopl{Res$^1$
}_{z=0}(f)$ stands for the ordinary residue and  $\mopl{Res$^{-1}$ }_{z=0}(f)$
stands for the coefficient of $z$ in the Laurent expansion of $f(z)$. Picking out the coefficient of
$z^{-m}$ in \eqref{recurrence2} (for some $m\in\Z$) we get:
\begin{eqnarray}
&{}&\cutoffsum_{1\le n_l<\cdots<
  n_1}\hskip -5mm (n_1+v)^{\beta_1(z)}\cdots
(n_l+v)^{\beta_l(z)}\nonumber\\
&\hskip -50mm =&\hskip -25mm\mopl{Res$^1$}_{z=0}H_l^0(z)\mopl{Res$^{m-1}$ }_{z=0}\hskip -7mm\cutoffsum_{1\le
  n_{l-1}<\cdots< n_1}\hskip -3mm (n_1+v)^{\beta_1(z)}\cdots
(n_{l-2}+v)^{\beta_{l-2}(z)}(n_{l-1}+v)^{\beta_{l-1}(z)+\beta_l(z)+1}\label{res1}\\
&\hskip -50mm +&\hskip -25mm\mopl{Res$^1$}_{z=0}K_l(z)\mopl{Res$^{m-1}$}_{z=0}\ \cutoffsum_{1\le n_{l-1}<\cdots< n_1} (n_1+v)^{\beta_1(z)}\cdots
(n_{l-1}+v)^{\beta_{l-1}(z)}\label{res2}\\
&\hskip -50mm +&\hskip -25mm \sum_{p\ge
  0}\frac{{H_l^0}^{(p)}(0)}{p!}\mopl{Res$^{m+p}$ }_{z=0}\hskip
-8mm\cutoffsum_{1\le n_{l-1}<\cdots< n_1} \hskip -5mm (n_1+v)^{\beta_1(z)}\cdots
(n_{l-2}+v)^{\beta_{l-2}(z)}(n_{l-1}+v)^{\beta_{l-1}(z)+\beta_l(z)+1}\label{res3}\\
&\hskip -50mm +& \hskip -25mm\sum_{p\ge
  0}\frac{{K_l}^{(p)}(0)}{p!}\mopl{Res$^{m+p}$}_{z=0}\cutoffsum_{1\le n_{l-1}<\cdots< n_1} (n_1+v)^{\beta_1(z)}\cdots
(n_{l-1}+v)^{\beta_{l-1}(z)}\label{res4}\\
&\hskip -50mm +&\hskip -25mm\sum_{j=1}^{2J_l}\sum_{p\ge 0}\frac{{H_l^j}^{(p)}(0)}{p!}\mopl{Res$^{m+p}$ }_{z=0}\hskip -2mm\cutoffsum_{1\le n_{l-1}<\cdots< n_1}\hskip -6mm (n_1+v)^{\beta_1(z)}\cdots
(n_{l-2}+v)^{\beta_{l-2}(z)}(n_{l-1}+v)^{\kappa_{l-1}^j(z)}.
\label{res5}
\end{eqnarray}
\begin{itemize}
\item  Suppose first $b_l\in\N$ and $m\ge 1$ in the equation above. The terms
  \eqref{res1} and \eqref{res2} vanish. The terms \eqref{res3} and 
  \eqref{res4} also vanish according to the induction
  hypothesis. According to the induction assumption, all the terms in the sum
  \eqref{res5} vanish as well, except possibly when $p=0$ and $m=1$, and when $j$ is big
  enough to ensure $\kappa_{l-1}^j(0)\le -1$. But in that case $j\ge
  b_{l-1}+b_l+2\ge b_l+2$, which implies $H_l^j(0)=0$. The left-hand side
  vanishes then for any $m\ge 1$, hence the function $z\mapsto\cutoffsum_{1\le n_l<\cdots<
  n_1}\hskip -5mm (n_1+v)^{\beta_1(z)}\cdots
(n_l+v)^{\beta_l(z)}$ is holomorphic at $z=0$.\\
\vskip -2mm
Now set $m=0$ in the equation above. The terms
  \eqref{res1} and \eqref{res2} still vanish, and only $p=0$ and $p=1$ give
  nonzero contribution in the terms \eqref{res3}, \eqref{res4} and
  \eqref{res5} according to the induction hypothesis. For $p=1$, these three
  terms are rational expressions in $c_1,\ldots
  ,c_{l}$ and in the first derivative of $\gamma$ at $z=0$ (which is $1$ by
  hypothesis).  For $p=0$ this is still the case for \eqref{res3} and
  \eqref{res4} as well as for the terms in the sum \eqref{res5} corresponding
  to $j\le b_l+b_{l-1}+1$, again  due to the induction hypothesis. The other terms in the sum \eqref{res5} do not
  contribute as $H_l^j(0)=0$ for $j\ge b_l+1$.
\item Now suppose $b_l\le -1$ and suppose $m\ge 2$ in the equation above. The
  terms \eqref{res3}, \eqref{res4} and \eqref{res5} vanish according to the
  induction hypothesis. The terms \eqref{res1} and \eqref{res2} also vanish
  unless possibly when $b_l=-1$. But in that case we have
  $b_{l-1}+b_l+1=b_{l-1}\ge 0$, and these two terms also vanish in view of
  the induction hypothesis. Hence the pole is of order at most one.\\
\vskip -2mm
Finally set $m=1$ in the equation above. Only $p=0$ gives a
  nonzero contribution in the terms \eqref{res3}, \eqref{res4} and
  \eqref{res5} according to the induction hypothesis. The residue $\mopl{Res$^1$ }_{z=0}$ is then a rational expression in $c_1,\ldots
  ,c_{l}$ independent of $\gamma$ according to the induction hypothesis and
  the explicit expression of the residues of $H_l^0$ and $K_l$ at $z=0$.
\endsquare
\end{itemize}
\end{enumerate}
\begin{thm}\label{thm:fpRieszsumhom} 
Let $v\in\R^+$ and let ${\cal R}_{0,v}: \tau \mapsto \left(\xi\mapsto
  \tau(z)(\xi)=\tau(\xi)(\xi+v)^{-z}\right)$ be the ordinary Riesz
regularisation shifted by $v$,
and let ${\cal R}_v(\tau)(z):= {\cal R}_{0,v}(\tau) (\gamma(z))$
where $\gamma$  is a holomorphic function of $z$ with $\gamma'(0)=1$.  
Let us set $\sigma_{s_i}^v(\xi)= \chi(\xi)\, (\xi+v)^{-s_i}$ and 
$\sigma^v:=\sigma_{-a_1}^v\otimes\cdots\otimes\sigma_{-a_k}^v$. At  non-positive
integer arguments $s_1=-a_1, \ldots,s_k= -a_k$, the   map 
$$z\mapsto \cutoffsum_<^{\smop{Chen}}\tilde{\cal R}_v^*(\sigma^v)(z)
$$ is holomorphic at $z=0$ and 
the renormalised multiple zeta values 
\begin{equation}
\zeta(-a_1, \ldots, -a_k; \, v)=\mopl{lim}_{z\to
  0}\cutoffsum_<^{\smop{Chen}}\tilde{\cal R}_v^*(\sigma^v)(z)
\end{equation}
is a polynomial of degree $k$ in $v$ with rational coefficients independent of the regularisation function $\gamma$.
\end{thm}
{\bf Proof:} We have $\zeta^{{\cal R}_v}(-a_1, \ldots, -a_k)=\Psi^{{\cal R}_v}_+(\sigma^v)(z)\restr{z=0}$, where $\Psi^{{\cal R}_v}(z)(\sigma^v)=\cutoffsum_<^{\smop{Chen}}\widetilde {\cal R}_v^\star(\sigma^v)(z)$ is a finite linear combination with rational coefficients of terms of the type:
$$\cutoffsum_<^{\smop{Chen}}\sigma_{-b_1+c_1\gamma(z),v}\otimes\cdots\otimes\sigma_{-b_l+c_l\gamma(z),v}$$
where $l\in\{1,\ldots ,k\}$, where $b_1,\ldots , b_l\in\N$ and where the
$c_j$'s are positive integers. This comes from the explicit expression:
\begin{equation}\label{eq:explicitexpr}
{\cal R}_v^\star(\sigma^v)(z)
=\hskip -3mm \sum_{I=(i_1,\ldots ,i_r)\in{\cal P}(k)}\frac{(-1)^{k-r}}{i_1\cdots i_r}
\sum_{J=(c_1,\ldots ,c_l)\in{\cal P}(r)}\frac{1}{c_1!\cdots c_l!}\sigma_{-b_1+c_1\gamma(z)}^{v}\otimes\cdots\otimes\sigma_{-b_l+c_l\gamma(z)}^{v},
\end{equation}
where $b_j$ stands for the sum of the $a_s's$ for $s$ inside the
$j^{\smop{th}}$ packet of the product composition $J\circ I$. The result then follows from  Theorem
\ref{thm:recurrence}, according to which
these maps are holomorphic at $z=0$ and their  values at $z=0$  rational
numbers independent of the choice of $\gamma$. Thus,  the restriction of
$\Psi^{{\cal R}}$ to the Hopf subalgebra ${\cal T}({\cal A}_{v,\N})$ of ${\cal
  T}({\cal A}_v)$, 
spanned by $\{\sigma_{-a}^{v}:\xi\mapsto \chi(\xi)\,(\xi+v)^{a},\,a\in\N\}$, takes its values in holomorphic
functions at $z=0$. This implies that $\Psi^{{\cal R}_v}_+\restr{{\cal T}({\cal
    A}_{v, \N})}=\Psi^{{\cal R}_v}\restr{{\cal T}({\cal A}_{v,\N})}$.\endsquare
\\ \\
Following these steps and using \eqref{eq:Rtildestar}, one obtains the following explicit formula for multiple zeta
values (i.e. $v=0$) at two nonpositive arguments:
\begin{eqnarray}\label{multiple zeta2arg}
\displaystyle
\hskip -30mm\zeta(-a,-b)&=&\mopl{fp }_{z=0}\ \cutoffsum_{1\le n_2< n_1}n_1^{a-z}n_2^{b-z}\nonumber\\
&\hskip -40mm=&\hskip -22mm \frac{1}{b+1}\sum_{s=0}^{b+1}{b+1\choose s}B_s\zeta(-a-b+s-1)+\zeta(-a)\zeta(-b)+(-1)^{a+1}\frac{a!b!}{2(a+b+2)!}B_{a+b+2}.
\end{eqnarray}
Formula \eqref{multiple zeta2arg} yields the following table of values
$\zeta(-a, -b)$ for $a,b\in \{0, \ldots, 6\}$:
\begin{equation*}\hskip -8mm
\begin{disarray}{c|c|c|c|c|c|c|c|}
\zeta(-a,-b)&a=0 &a=1&a=2&a=3&a=4&a=5&a=6\\ \hline
&&&&&&&\\
b=0 &\frac{3}{8}&\frac{1}{12}&\frac{7}{720}&-\frac{1}{120}
&-\frac{11}{2\, 520}&\frac{1}{252}&\frac{1}{224}\\ 
&&&&&&&\\\hline
&&&&&&&\\
b=1 &\frac{1}{24}&\frac{1}{288}&-\frac{1}{240}&-\frac{19}{10\, 080}
&\frac{1}{504}&\frac{41}{20\, 160}&-\frac{1}{480}\\ 
&&&&&&&\\ \hline
&&&&&&&\\
b=2&-\frac{7}{720}&-\frac{1}{240}&0&\frac{1}{504}&\frac{113}{151\, 200}&-\frac{1}{480}&-\frac{307}{166\,320}\\ 
&&&&&&&\\ \hline
&&&&&&&\\
b=3&-\frac{1}{240}&\frac{1}{840}&\frac{1}{504}&\frac{1}{28\,800}&-\frac{1}{480}&-\frac{281}{332\,640}&\frac{1}{264}\\ 
&&&&&&&\\ \hline
&&&&&&&\\
b=4&\frac{11}{2\,520}&\frac{1}{504}&-\frac{113}{151\,200}&-\frac{1}{480}&0&\frac{1}{264}&\frac{117\,977}{75\,675\,600}\\ 
&&&&&&&\\ \hline
&&&&&&&\\
b=5&\frac{1}{504}&-\frac{103}{60\,480}&-\frac{1}{480}&\frac{1}{1232}&\frac{1}{264}&\frac{1}{127\,008}&-\frac{691}{65\,520}\\ 
&&&&&&&\\ \hline
&&&&&&&\\
b=6&-\frac{1}{224}&-\frac{1}{480}&\frac{307}{166\,320}&\frac{1}{264}&-\frac{117\,977}{75\,675\,600}&-\frac{691}{65\,520}&0\\
&&&&&&&\\ \hline
\end{disarray}
\end{equation*}
\begin{rk}
If  $a+b$
is odd and $b\not =0$, all terms in equation (\ref{multiple zeta2arg}) vanish
except the one with $s=1$. This yields:
\begin{equation}
\zeta(-a,-b)=-\frac 12 \zeta(-a-b),
\end{equation}
a fact also established by L. Guo and B. Zhang by different means \cite{GZ}. We also have for odd
$a$~:
\begin{equation}
\zeta(-a,0)=-\zeta(-a).
\end{equation}
This is due to the fact that both versions of the
renormalised value coincide with the analytic continuation of the
double zeta function, which is holomorphic at those points. As a
consequence of the quasi-shuffle relations, the table of
values in \cite{GZ} also matches with ours on the diagonal points $(-a,-a)$. It
however diverges from ours at other singular points.
\end{rk}
\subsection {An alternative renormalization method}\label{sec:AET}
An alternative proposal for multiple zeta values at non-positive arguments can be found in
\cite{AET} (Remark 2 therein) and easily extended to Hurwitz multiple zeta
values at non-positive arguments  defining  as follows\footnote{modulo the reverse convention in the
order of the arguments}:
\begin{equation*}
\zeta^{\rm alt}(-a_1,\ldots,-a_k; v):=\mopl{lim}_{z\to 0}\zeta(-a_1+z,\ldots,-a_k+z;v).
\end{equation*}
This uses the $\widetilde {\cal R}$ regularisation instead of our twisted
$\widetilde{\cal R}^\star$ regularisation,  so that  the stuffle
relations are not expected to hold. Multiple Hurwitz zeta functions at non-positive integer arguments $-a_1,\ldots, -a_k$ are a priori
 sensitive to the
twisting of the regularisation by the Hoffman isomorphism because of the 
dependence on the coefficients $c_1>0, \ldots, c_k>0$ of the limit $\lim_{z\to 0} \zeta(-a_1+c_1\, z,
\ldots, -a_k+c_k\,z;v)$ studied in Theorem \ref{thm:recurrence}. They
  nevertheless coincide for depths $1$ and $2$:
\begin{prop}\label{prop:alt12}For any non positive integers $a_1$ and $a_2$
  and any non negative real number $v$, we have: $$\zeta(-a_1;v) = \zeta^{\rm alt}(-a_1
  ;v)\quad{\rm and}\quad   \zeta(-a_1,  -a_2 ;v)=\zeta^{\rm alt}(-a_1,
  -a_2:v).$$
\end{prop}
{\bf Proof:}  This follows from 
 the holomorphicity of the ordinary Hurwitz zeta function at
non-positive integers.\endsquare\\

However, the two methods are expected to  differ for depth $k>2$;
  this is confirmed 
 in 
the depth $3$ case  investigated in the following example:  
\begin{ex} For $v=0$ and any three symbols $\sigma_1,\sigma_2,\sigma_3$ corresponding to non-positive integer
arguments $s_1=-a_1, s_2=-a_2, s_3=-a_3$, one infers from the explicit
formulae for Hoffman's $\log$ and $\exp$ given in Section \ref{sect:9} that:
\begin{eqnarray*}
&&\tilde{\cal
  R}^*(\sigma_1\otimes\sigma_2\otimes\sigma_3)(z)-\widetilde {\cal R}\left(\sigma_1\otimes\sigma_2\otimes\sigma_3\right)(z)\\
&=&\frac
12\Big(\sigma_1(z)\sigma_2(z)\otimes\sigma_3(z)+\sigma_1(z)\otimes\sigma_2(z)\sigma_3(z)-(\sigma_1\sigma_2)(z)\otimes\sigma_3(z)-\sigma_1(z)\otimes(\sigma_2\sigma_3)(z)\Big)\\
&+&\frac 16\sigma_1(z)\sigma_2(z)\sigma_3(z)-\frac
14(\sigma_1\sigma_2)(z)\sigma_3(z)-\frac 14 \sigma_1(z)(\sigma_2\sigma_3)(z)
+\frac 13\left( \sigma_1\sigma_2\sigma_3\right)(z),\\
\end{eqnarray*}
so that, using the holomorphicity of the zeta function at non-positive integer
arguments, we get  
\begin{eqnarray*}
&&\Delta(-a_1,-a_2, -a_3)= \zeta(-a_1,-a_2, -a_3)-\zeta^{\rm alt}(-a_1,-a_2, -a_3) \\ 
&&\hskip -35mm=\lim_{z\to 0}\frac
12\Big(\zeta(-a_1-a_2+2z, -a_3+z) +\zeta(-a_1+z, -a_2-a_3+2z)-\zeta(-a_1-a_2+z,
-a_3+z) -\zeta(-a_1+z, -a_2-a_3+z)\Big)\\
&&\hskip -30mm+\lim_{z\to 0}\left[\frac 16 \zeta (-a_1-a_2-a_3+3z)-\frac
14\zeta(-a_1-a_2-a-_3+2z) -\frac 14\zeta(-a_1-a_2-a_3+2z) 
+\frac 13\zeta(-a_1-a_2-a_3+z)\right] \\
&&\hskip -35mm= \lim_{z\to 0}\frac
12\Big(\zeta(-a_1-a_2+2z, -a_3+z) +\zeta(-a_1+z, -a_2-a_3+2z)-\zeta(-a_1-a_2+z,
-a_3+z) -\zeta(-a_1+z, -a_2-a_3+z)\Big).\\
\end{eqnarray*}
In particular, we have
\begin{eqnarray*}
&&\Delta(0,-1,-1;v) =\zeta(0,-1,-1;v)-\zeta^{\rm alt}(0,-1,-1;v)\\
&=&\frac 12\mopl{fp}_{z\to
    0}\mopl{fp}_{N\to\infty}\left(\sum_{n_1=1}^N\sum_{n_2=1}^{n_1-1}(n_1+v)^{1-2z}(n_2+v)^{1-z}
+\sum_{n_1=1}^N\sum_{n_2=1}^{n_1-1}(n_1+v)^{-z}(n_2+v)^{2-2z}\right)\\
&&-\frac 12\mopl{fp}_{z\to
    0}\mopl{fp}_{N\to\infty}\left(\sum_{n_1=1}^N\sum_{n_2=1}^{n_1-1}(n_1+v)^{1-z}(n_2+v)^{1-z}-\sum_{n_1=1}^N\sum_{n_2=1}^{n_1-1}(n_1+v)^{-z}(n_2+v)^{2-z}\right)\\
&=&\frac{B_4}{288}(2-8-3+6)\not =0.
\end{eqnarray*}
Remarkably enough, the result does not depend on $v$.
\end{ex}
\subsection{Two identities preserved at nonpositive arguments}
Hurwitz multiple zeta functions verify the following two identities on the
domain of convergence, namely for any $(s_1,\ldots ,s_k)\in\C^k$ with
$\mop{Re}(s_1+\cdots +s_m)>m$ for all $m\in\{1,\ldots ,k\}$:
\begin{equation}\label{eq:hurwitz1bis}
\zeta(s_1,\ldots ,s_k;v+1)=\zeta(s_1,\ldots
,s_k;v)-(v+1)^{-s_k}\zeta(s_1,\ldots ,s_{k-1};v+1)
\end{equation}
and
\begin{equation}\label{eq:hurwitz2bis}
\frac{d}{dv}\zeta(s_1,\ldots ,s_k;v)=\sum_{j=1}^k -s_j\zeta(s_1,\ldots
,s_{j-1},s_j+1,s_{j+1},\ldots ,s_k; v).
\end{equation}
The proof is straightforward and left to the reader. We show here that both
identities are preserved at nonpositive integer arguments, for both
renormalisations $\zeta$ and $\zeta^{\smop{alt}}$.
\begin{prop}\label{prop:deriveHurwitz} For any $a_1,\ldots , a_k\in\N$ we have:
\begin{equation}
\zeta^{\rm alt}(-a_1,\ldots ,-a_k;v+1)=\zeta^{\rm alt}(-a_1,\ldots
,-a_k;v)-(v+1)^{a_k}\zeta^{\rm alt}(-a_1,\ldots ,-a_{k-1};v+1)
\end{equation}
and moreover if $a_1,\ldots , a_k\ge 1$:
\begin{equation}
\frac{d}{dv}\zeta^{\rm alt}(-a_1,\ldots ,-a_k;v)=-\sum_{j=1}^k a_j\zeta^{\rm alt}(-a_1,\ldots
,-a_{j-1},-a_j+1,-a_{j+1},\ldots ,-a_k; v).
\end{equation}
\end{prop}
{\bf Proof:} This follows from the  two identities of meromorphic functions (with respect to the variables $s_j$):
\begin{equation}
\zeta(s_1,\ldots ,s_k;v+1)=\zeta(s_1,\ldots ,s_k;v)-(v+1)^{-s_k}\zeta(s_1,\ldots ,s_{k-1};v+1)
\end{equation}
and:
\begin{equation}
\frac{d}{dv}\zeta(s_1,\ldots ,s_k;v)=-\sum_{j=1}^k s_j\zeta(s_1,\ldots
,s_{j-1},s_j+1,s_{j+1},\ldots ,s_k; v).
\end{equation}
which in turn yield  the following  identities of holomorphic functions:
\begin{equation}
\zeta(-a_1+z,\ldots ,-a_k+z;v+1)=\zeta(-a_1+z,\ldots ,-a_k+z;v)-(v+1)^{a_k+z}\zeta(-a_1+z,\ldots ,-a_{k-1}+z;v+1)
\end{equation}
and:
\begin{equation}
\frac{d}{dv}\zeta(-a_1+z,\ldots ,-a_k+z;v)=\sum_{j=1}^k a_j\zeta(-a_1,\ldots
,-a_{j-1},-a_j+1,-a_{j+1},\ldots ,-a_k; v).
\end{equation}
\endsquare
\begin{thm}\label{prop:deriveHurwitz2} For any $a_1,\ldots , a_k\in\N$ we have:
\begin{equation}\label{eq:hurwitz1ter}
\zeta(-a_1,\ldots ,-a_k;v+1)=\zeta(-a_1,\ldots
,-a_k;v)-(v+1)^{a_k}\zeta(-a_1,\ldots ,-a_{k-1};v+1)
\end{equation}
and moreover when $a_1,\ldots , a_k\ge 1$:
\begin{equation}\label{eq:hurwitz2ter}
\frac{d}{dv}\zeta(-a_1,\ldots ,-a_k;v)=-\sum_{j=1}^k a_j\zeta(-a_1,\ldots
,-a_{j-1},-a_j+1,-a_{j+1},\ldots ,-a_k; v).
\end{equation}
\end{thm}
{\bf Proof:} Equation \eqref{eq:hurwitz1ter} is equivalent to:
\begin{equation}\label{eq:althurwitz1}
\zeta(-a_1,\ldots, -a_k;v+1)=\sum_{r=0}^{k}(-1)^{r}(v+1)^{a_k+\cdots +a_{k-r+1}}\zeta(-a_1,\ldots, -a_{k-r};v),
\end{equation}
where the term corresponding to $r=0$ (resp. $r=k$) is $\zeta(-a_1,\ldots,
-a_k;v)$ (resp. $(-1)^k(v+1)^{a_k+\cdots +a_{1}}$).
We prove
\eqref{eq:althurwitz1} by induction on the depth $k$. We denote by $\delta_1$ any positively supported symbol such that $\delta_1(1)=1$ and $\delta_1(\xi)=0$ for $\xi\ge 2$. We set
$\sigma^v(\xi):=\sigma(\xi+v)$ for any symbol $\sigma$, as well as
$(\sigma\otimes\tau)^v:=\sigma^v\otimes\tau^v$ and so on. Using
the notations of Theorem \ref{thm:fpRieszsumhom} we compute:
\begin{eqnarray*}
&&\zeta(-a_1,\ldots, -a_k;v)-\zeta(-a_1,\ldots, -a_k;v+1)\\
&=&\mopl{lim}_{z\to 0}\cutoffsum^{\smop{Chen}}{\widetilde {\cal R}}^*_v(\sigma_{s_1}^v\otimes\cdots\otimes\sigma_{s_k}^v)(z)-{\widetilde {\cal R}}^*_{v+1}(\sigma_{s_1}^{v+1}\otimes\cdots\otimes\sigma_{s_k}^{v+1})(z)\\
&=&\mopl{lim}_{z\to 0}\cutoffsum^{\smop{Chen}}\big({\widetilde {\cal R}}^*_0(\sigma_{s_1}\otimes\cdots\otimes\sigma_{s_k})(z)^v-{\widetilde {\cal R}}^*_0(\sigma_{s_1}\otimes\cdots\otimes\sigma_{s_k})(z)^{v+1}\big)\\
&=&\mopl{lim}_{z\to 0}\cutoffsum^{\smop{Chen}}{\widetilde {\cal R}}^*_v(\sigma_{s_1}^v\otimes\cdots\otimes\sigma_{s_{k-1}}^v\otimes \sigma_{s_k}^v\delta_1)(z)\\
&=&\mopl{lim}_{z\to 0}\cutoffsum^{\smop{Chen}}{\widetilde {\cal R}}^*_v\big((\sigma_{s_1}^v\otimes\cdots\otimes\sigma_{s_{k-1}}^v)\star \sigma_{s_k}^v\delta_1\big)(z)-\mopl{lim}_{z\to 0}\cutoffsum^{\smop{Chen}}{\widetilde {\cal R}}^*_v\big((\sigma_{s_1}^v\otimes\cdots\otimes\sigma_{s_{k-2}}^v)\otimes\sigma_{s_{k-1}}^v \sigma_{s_k}^v\delta_1\big)(z)\\
&=&\sum_{r=1}^k(-1)^{r-1}\mopl{lim}_{z\to 0}\cutoffsum^{\smop{Chen}}{\widetilde {\cal R}}^*_v\big((\sigma_{s_1}^v\otimes\cdots\otimes\sigma_{s_{k-r}}^v)\star (\sigma_{s_{k-r+1}}^v\cdots\sigma_{s_k}^v\delta_1)\big)(z)\\
&=&\sum_{r=1}^k(-1)^{r-1}(v+1)^{a_k+\cdots +a_{k-r+1}}\mopl{lim}_{z\to 0}\cutoffsum^{\smop{Chen}}{\widetilde {\cal R}}^*_v(\sigma_{s_1}^v\otimes\cdots\otimes\sigma_{s_{k-r}}^v)(z)\\
&=&\sum_{r=1}^k(-1)^{r-1}(v+1)^{a_k+\cdots +a_{k-r+1}}\zeta(-a_1,\ldots, -a_{k-r};v).
\end{eqnarray*}
Here we have used Theorem \ref{thm:fpRieszsumhom} and the expression
\eqref{eq:explicitexpr} implicitly. Equation \eqref{eq:hurwitz2ter} is proved as follows,
still applying Theorem \ref{thm:fpRieszsumhom}. Here $a_1,\ldots,a_k$ are
positive integers (i.e.$\ge 1$) with $s_j=-a_j,\,j=1,\ldots ,k$, and the operator $\frac d{dv}$
will also be denoted by a dot:
\begin{eqnarray*}
&&\frac{d}{dv}\zeta(-a_1,\ldots
,-a_k;v)=\frac{d}{dv}\left(\cutoffsum^{\smop{Chen}}{\widetilde {\cal
      R}}^*_v(\sigma_{s_1}^v\otimes\cdots\otimes\sigma_{s_k}^v)(z)\restr{z=0}\right)\\
&=&\left(\frac{d}{dv}\cutoffsum^{\smop{Chen}}{\widetilde {\cal
      R}}^*_v(\sigma_{s_1}^v\otimes\cdots\otimes\sigma_{s_k}^v)(z)\right)\restr{z=0}
\hbox{ (by holomorphy at }z=0)\\
&=&\left(\cutoffsum^{\smop{Chen}}\frac{d}{dv}{\widetilde {\cal
      R}}^*_v(\sigma_{s_1}^v\otimes\cdots\otimes\sigma_{s_k}^v)(z)\right)\restr{z=0}\\
&=&\left(\cutoffsum^{\smop{Chen}}\frac{d}{dv}\exp\circ{\widetilde {\cal
      R}}_v\circ\log(\sigma_{s_1}^v\otimes\cdots\otimes\sigma_{s_k}^v)(z)\right)\restr{z=0}\\
&=&\left(\cutoffsum^{\smop{Chen}}\Big(\exp\circ\dot{{\widetilde {\cal
      R}}}_v\circ\log(\sigma_{s_1}^v\otimes\cdots\otimes\sigma_{s_k}^v)(z)
+\sum_{j=0}^k\exp\circ{\widetilde {\cal
      R}}_v\circ\log(\sigma_{s_1}^v\otimes\cdots\otimes\dot\sigma_{s_j}^v\otimes\cdots\otimes\sigma_{s_k}^v)(z)\Big)\right)\restr{z=0}\\
&=&-\sum_{j=1}^k(z+s_j)\left(\cutoffsum^{\smop{Chen}}\exp\circ{\widetilde {\cal
      R}}_v\circ\log(\sigma_{s_1}^v\otimes\cdots\otimes\sigma_{s_j+1}^v\otimes\cdots\otimes\sigma_{s_k}^v)(z)\right)\restr{z=0}\\
&=&\sum_{j=1}^ka_j\zeta(-a_1,\ldots,-a_{j-1},-a_j+1,-a_{j+1},\ldots,a_k;v).
\end{eqnarray*}
\endsquare 
\vfill \eject \noindent
\section{A  HIGHER DIMENSIONAL ANALOG OF MULTIPLE  ZETA
  FUNCTIONS}\label{sect:hdim}
Using radial symbols  we extend the constructions carried out previously in
the one-dimensional case to a higher dimensional setup using nested
  cubes. Proposition \ref{prop:HigherdimChensum} and Theorem
  \ref{thm:hmz-explicit} relate the higher dimensional multiple zeta functions to the one-dimensional
ones, thus enabling us to carry out to this higher
dimensional setting the techniques implemented in the one dimensional case.
\subsection{An algebra of radial symbols}\label{sect:mrb}
 Let  $\vert \cdot\vert$ be any  norm on $\R^n$.  Let   $\bullet $ be the
 product  on  complex valued functions on  $\Z^n$ vanishing at the origin,
defined  by
\begin{equation}\label{eq:bullet}\sigma_1 \bullet \sigma_2 (x):=  \frac{1}{A_n(|x|)} 
\sum_{\vert x_1\vert = \vert x_2\vert=\vert x\vert } \sigma_1(x_1)\cdot
 \sigma_2(x_2) \quad \forall x\in \Z^n,
\end{equation}
where for a non-negative number $t=|x|$ for some $x\in\Z^n$, the number of points with integer coordinates in the sphere
centered at zero of radius $t$ is denoted by  $A_n(t)$. We have $A_n(0)=1$
and, when $\vert \cdot\vert$ is the
supremum norm $\vert x\vert={\rm sup}_{i=1}^n\vert x_i\vert$ on $\R^n$, then $A_n(t)= (2t+1)^n- (2t-1)^n$ for any positive integer $t$. The product $\bullet$ is commutative and associative 
 as a result of  the commutativity and associativity of the ordinary product.
\begin{defn} For any complex valued function $\sigma$ on $\Z^n$ and any $x\in \Z^n$ let us set:
 \begin{equation}\label{eq:R}
 R_n(\sigma)(x):=\sum_{0<\vert  y \vert<  
 \vert x\vert, \,  y\in \Z^n} \sigma( y).\end{equation}
\end{defn}
This construction can be immediately extended to complex valued functions on $\R^n$ provided the function
$A_n$ is suitably interpolated to all positive real numbers by a nowhere
vanishing function. We now specialise to radial functions $\sigma=f\circ\vert \cdot\vert$ on
$\R^n$ in which
case the product $\bullet$ reads:
\begin{equation}\label{bulletradial}
\sigma_1 \bullet \sigma_2(x)= A_n(|x|) \sigma_1(x)\sigma_2(x).
\end{equation}
The radial part $\rho_n$ of the operator $R_n$, defined by:
\begin{equation}\label{eq:radial}
R_n(f\circ|\cdot|)=\rho_n(f)\circ|\cdot|
\end{equation}
reads as follows:
\begin{equation}
\rho_n(f)(N):= \sum_{0<p <  
 N} A_n(p) \, f(p),
\end{equation}
where the sum runs over all positive real numbers $p$ such that the sphere of
radius $p$ has non-empty intersection with $\Z^n$. Let ${\cal P}_1={\cal
  P}^{*,0}$ be the algebra of positively supported classical symbols, let
${\cal A}$ be a subalgebra of ${\cal P}_1$,
and let $ M$ the linear operator on ${\cal A}$ with values into radial
functions given by:
\begin{equation}
M(f):=\frac{f}{A_n}\circ \vert \cdot\vert.
\end{equation}
Let ${\cal A}_n$ be the image of ${\cal A}$ by the operator $ M$;
 by (\ref{bulletradial})  we have  $M(f_1\, f_2)=M(f_1)\bullet M(f_2)$ so that  $M$ yields an algebra isomorphism from $({\cal A},.)$ onto $({\cal
  A}_n,\bullet)$. This induces a Hopf algebra isomorphism
$$\cal M:\left({\cal
  H}_1, \star\right)\to \left({\cal H}_n, \star_\bullet\right)$$
between the
two corresponding quasi-shuffle Hopf algebras.
\subsection{Higher-dimensional renormalised multiple zeta values}
We now specialise to the supremum norm $\vert x\vert={\rm sup}_{i=1}^n\vert
x_i\vert$ on $\R^n$. The reason for this is that the image of $\Z^n$ by this
norm is equal to $\N$, so that sums over $\Z^n$ boil down to sums over the
integers. The methods developped previously in the one-dimensional
  situation can then be transposed to the $n$-dimensional 
  one by means of the Hopf algebra isomorphism ${\cal M}$ defined in the
  previous paragraph. We do not know at
  this stage how to build
  meromorphic extensions of multiple zeta functions for radial
  symbols when $\vert\cdot\vert$ is the
  Euclidean norm since we do not know how to count the integer points inside
  the unit sphere.\\
  
We interpolate the number of points with integer coordinates on spheres by setting:
\begin{equation}
A_n(t)=(2t+1)^n-(2t-1)^n
\end{equation}
for any $t\ge 0$. The vanishing of this interpolated $A_n$ at $t=0$ is harmless as we deal with symbols which vanish at the origin.
\begin{defn}\label{def:hdrensum}
Let ${\cal A}$ be a subalgebra of ${\cal P}_1$ and let ${\cal R}$ be a
holomorphic regularisation on ${\cal A}$. Let
$\cutoffsum_{<}^{\smop{Chen},{\cal R}}$ be the
renormalised nested sum character of Definition \ref{defn:Rregmultiple zeta}
(in the strict inequality version) on ${\cal H}_1={\cal
  T}({\cal A})$. The {\rm higher-dimensional renormalised nested sum
  character on radial symbols\/} is defined by~:
\begin{equation}
\cutoffsum_{<}^{\smop{Chen},{\cal R},n}\sigma:=\cutoffsum_{<}^{\smop{Chen},{\cal
    R}}{\cal M}^{-1}(\sigma).
\end{equation}
for any $\sigma\in{\cal H}_n$.
\end{defn}
The fact that this defines a character  for the
$\star_{\bullet}$ quasi-shuffle product on ${\cal H}_n$  follows from 
$$\cutoffsum_{<}^{\smop{Chen},{\cal
    R}}{\cal M}^{-1}(\sigma)=  \psi^{\cal R} \circ {\cal M}^{-1}(\sigma)$$  where
$\psi^{\cal R}: ({\cal H}_1, \star)\to \C$ is 
 the Birkhoff-Hopf character introduced in Theorem \ref{thm:renChensums}. A similar definition
holds for the weak inequality version $\cutoffsum_{\le}^{\smop{Chen},{\cal
    R},n}$: details are left to the reader. We now apply the above constructions to symbols $\sigma^v_{s}\sim
 (\vert\cdot\vert+v)^{-s}$, $s\in \C$, $v\in \R_+$ i.e. to the
 algebra  of radial classical symbols:
  $$\widetilde {\cal A}_n:=\{f\circ\vert
 \cdot\vert, \quad f\in {\cal A}\}$$
 where ${\cal A}$
 was defined in Paragraph \ref{sect:stufflezeta}.
 We set  the following
 definition:
\begin{defn}
Let ${\cal R}$ be a
holomorphic regularisation on ${\cal A}$. Given  any $s_i\in \C, v_i\in \R_+$ and
$\sigma_{s_i}^{v_i}=(\vert \cdot\vert+v_i)^{-s_i}$ for some non negative real number
$v$,  we define the renormalised multiple Hurwitz zeta
 values at $(s_1, \ldots, s_k)$  by
 \begin{equation}
\zeta_n^{\cal R}(s_1, \ldots, s_k; \, v_1,\ldots,v_k):= \cutoffsum^{\smop{Chen}, {\cal
    R}, n}_<\sigma_{s_1}^{v_1}\otimes\cdots\otimes\sigma_{s_k}^{v_k}
\end{equation} where   the renormalised nested sum is taken as in Definition
\ref{def:hdrensum}. As in the one-dimensional case, we denote by 
$\zeta_n^{\cal R}(s_1, \ldots, s_k; \, v)$ the Hurwitz multiple zeta value
with parameters $v_1=\cdots= v_k=v$ and 
 define renormalised
higher-dimensional multiple zeta values by setting the parameters $v_i$ to zero:
\begin{equation}
\zeta_n^{\cal R}(s_1, \ldots, s_k):=\zeta_n^{\cal R}(s_1, \ldots, s_k;0,\ldots,0).
\end{equation}
\end{defn}
\subsection{Explicit formulae for higher-dimensional multiple zeta values}
The linear map  
\begin{eqnarray*}
N: {\cal
   A}&\to& \widetilde {\cal A}_n\\
 f&\mapsto & N(f):= f\circ \vert
 \cdot\vert
\end{eqnarray*} where as before $\vert\cdot\vert$ stands for the supremum norm, induces an isomorphism of tensor
 algebras ${\cal N}: {\cal T}({ \cal A})\to {\cal T}\left(\widetilde {\cal
     A}_n\right)$.    On the other hand, the composition
 $M\circ N^{-1}:\left(\widetilde {\cal A}_n , \cdot\right)\rightarrow
 \left({\cal A}_n, \bullet\right)$ 
 yields an isomorphism of  algebras. Let us observe that
\begin{eqnarray}
\zeta_n^{\cal R}(s_1, \ldots, s_k; \, v_1,\ldots,v_k)
&=&\cutoffsum^{\smop{Chen}, {\cal
    R}}_<{\cal  M}^{-1}\circ {\cal N}\left(f_{s_1}^{v_1}\otimes\cdots\otimes f_{s_k}^{v_k}\right)\nonumber\\
&=& \cutoffsum^{\smop{Chen}, {\cal
    R}}_< A_n \, f_{s_1}^{v_1}\otimes\cdots\otimes
    A_n\, f_{s_k}^{v_k}\label{eq:zetanAn}
\end{eqnarray}
where we have set $f_{s_i}^{v_i}(t):=(t+v_i)^{-s_i}$ and where
 $A_n$ was defined in Paragraph
\ref{sect:mrb}. Equation (\ref{eq:zetanAn})  yields  the following explicit expression for
higher-dimensional renormalised multiple zeta values (here $v=0$):
\begin{thm}\label{thm:hmz-explicit}
\begin{equation}\label{truc}
\hskip -11mm\zeta_n^{{\cal R}}(s_1,\ldots ,s_k)=\sum_{J=(j_1,\ldots
  ,j_k),\,j_r=0,\ldots,[\frac{n}{2}]}2^{nk-2|J]}{n\choose
  2j_1+1}\cdots{n\choose 2j_k+1}\,\zeta^{{\cal
    R}}(s_1-n+2j_1+1,\ldots,s_k-n+2j_k+1).
\end{equation}
\end{thm}
  Examples of
particular instances of \eqref{truc} are:
\begin{eqnarray*}
\zeta^{{\cal R}}_1(s)&=&2\zeta^{{\cal R}}(s)\\
\zeta^{{\cal R}}_2(s)&=&8\zeta^{{\cal R}}(s-1)\\
\zeta^{{\cal R}}_3(s)&=&24\zeta^{{\cal R}}(s-2)+2\zeta^{{\cal R}}(s)\\
\zeta^{{\cal R}}_1(s_1,s_2)&=&4\zeta^{{\cal R}}(s_1,s_2)\\
\zeta^{{\cal R}}_2(s_1,s_2)&=&64\zeta^{{\cal R}}(s_1-1,s_2-1)\\
\zeta^{{\cal R}}_3(s_1,s_2)&=&576\zeta^{{\cal R}}(s_1-2,s_2-2)+48\zeta^{{\cal R}}(s_1-2,s_2)+48\zeta^{{\cal R}}(s_1-2,s_2)+4\zeta^{{\cal R}}(s_1,s_2)\\
&\cdots&
\end{eqnarray*}
\subsection{Higher-dimensional  Hurwitz multiple zeta values
    renormalized at non positive integers}
Since $A_n(t)$ is polynomial in $t$ with integer coefficients:
\begin{equation}
A_n(t)=\sum_{j=0}^{[\frac{n}{2}]}2^{n-2j}{n\choose 2j+1}t^{n-2j-1},
\end{equation}
it can be written as a polynomial in the variable $t+v$ with
   coefficients given by polynomials  in $v$ with integer coefficients. Thus, 
 when $v_1=\cdots=v_k=v$ the tensor product $ A_n \, f_{s_1}^{v_1}\otimes\cdots\otimes
    A_n\, f_{s_k}^{v_k}$ can be written as a finite linear combination  of
    tensor products 
 $  f_{s_1-b_1}^{v_1}\otimes\cdots\otimes
     f_{s_k-b_k}^{v_k}$ with
     coefficients given by polynomials in $v$  with integer  coefficients and  for some non  negative integers $b_1, \ldots,
     b_k$. 
By Theorem \ref{thm:fpRieszsumhom}, setting $v_i=v$ and $s_i=-a_i$ for any $i\in
\{1, \cdots, k \}$, for 
non negative integers $a_i$  the expression $\cutoffsum^{\smop{Chen}, {\cal
    R}}_< f_{s_1-b_1}^{v_1}\otimes\cdots\otimes
     f_{s_k-b_k}^{v_k}$ is a polynomial expression in $v$ with rational
     coefficients, hence following result. 
\begin{prop}  Multiple zeta values
$\zeta_n^{\cal R}(-a_1, \ldots, -a_k; \, v)$ at non positive integers  are polynomials with rational
coefficients in $v$. In particular, higher-dimensional multiple zeta values are rational at nonpositive arguments.
\end{prop}
\vfill\eject \noindent
\section*{Appendix: A continuous analog of polylogarithms}
We describe a continuous analog of
polylogarithms similar to their discrete counterparts in number theory. In
contrast  to the latter, the shuffle
relations  obtained here are of the same type as  the
original ones but they hold on a  different Hopf
algebra. We define a continuous analog of polylogarithms setting for $1\geq r>0$,  $\vert z\vert
  <1$,   $s_1>1$ and $s_i\geq 1$ for $1<i\leq k$:
 \begin{equation}\label{eq:deftildeLi}
\tilde{Li}^r_{\underline s}(z)= \int_{r\leq  t_k\leq \cdots \leq t_1}
 dt_1\cdots  dt_k\, z^{t_1} \prod_{i=1}^k  \, t_i^{-s_i}. 
\end{equation} Note that 
$$\tilde{Li}^r_{\underline s}(1)=\tilde \zeta^{r}(\underline s).$$
\begin{prop}\label{prop:polylog}
 Let $\underline s\in \R^k$. For any $0< z<1$,  
\begin{enumerate}
\item if $s_1=1$ then $$\frac{d}{dz} \tilde{Li}_{\underline s}(z)= -\frac{1}{z\log z} \tilde{Li}^r_{s_2, \ldots, s_k}(z).$$
\item Otherwise 
$$\frac{d}{dz} \tilde{Li}^r_{\underline s}(z)= \frac{1}{ z} \tilde{Li}^r_{(s_1-1, s_2, \ldots, s_k)}(z).$$
\end{enumerate}
\end{prop}
{\bf Proof:}
\begin{enumerate}
\item Let us first assume that $s_1=1$. 
\begin{eqnarray*}
\frac{d}{dz} \tilde{Li}^r_{\underline s}(z)&=& \frac{d}{dz} \int_{r\leq  t_k\leq \cdots \leq t_2\leq t_1} dt_1\cdots  dt_k\, \frac{z^{t_1}}{t_1} \prod_{i=2}^k  \, t_i^{-s_i}\\
&=& \frac{d}{dz}  \int_{r\leq  t_k\leq \cdots \leq t_3\leq t_2} dt_2\cdots  dt_k\,
\cutoffint_{t_2}^\infty \frac{z^{t_1}}{t_1}dt_1\,  \prod_{i=2}^k  \, t_i^{-s_i}\\
&=&  \int_{r\leq  t_k\leq \cdots \leq t_3\leq t_2} dt_2\cdots  dt_k\,
 \left[\frac{z^{t_1-1}}{\log z}\right]_{t_2}^\infty   \prod_{i=2}^k  \, t_i^{-s_i}\\
&=&- \frac{1}{z \log z} \tilde{Li}^r_{( s_2, \ldots, s_k)}(z).
 \end{eqnarray*}
\item Let us now assume that $s_1\neq 1$. 
\begin{eqnarray*}
\frac{d}{dz} \tilde{Li}^r_{\underline s}(z)&=& \frac{d}{dz} \int_{r\leq  t_k\leq \cdots \leq t_1} dt_1\cdots  dt_k\, z^{t_1} \prod_{i=1}^k  \, t_i^{-s_i}\\
&=&\frac{1}{z} \int_{r\leq  t_k\leq \cdots \leq t_1} dt_1\cdots  dt_k\,  \frac{z^{t_1}}{t_1^{s_1-1}}\prod_{i=2}^k  \, t_i^{-s_i}\\
&=& \frac{1}{z} \tilde{Li}^r_{(s_1-1,  s_2, \ldots, s_k)}(z).
 \end{eqnarray*}
\end{enumerate}
\endsquare\\ \\
Following  number theorists (see e.g. \cite{C1}, \cite{ENR}, \cite{W}, \cite{Z}), we now restrict ourselves to positive integer values of the indices $s_i$, in which case  $\tilde \zeta(\underline s)$ converges if $s_1>1$ since all indices $s_i\geq 1$. To keep track of the condition $s_1>1$ it is useful  to introduce a coding of multi-indices $\underline s$ by words over the alphabet $X=\{x_0, x_1\}$ by the rule:
$$\underline  s \mapsto  x_0^{s_1-1}x_1 x_0^{s_2-1}x_1 \cdots x_0^{s_k-1} x_1.$$
This means that an $x_1$ arises every time we have $s_i=1$ for some $i$. In contrast,  $s_i>1$ gives rise to $x_0^{s_i-1} x_1$. Condition $s_1>1$ therefore translates to the requirement  that the word starts with $x_0$. Since the last letter is $x_1$ whether $s_k=1$ or not,   the words arising from a multiindex $(s_1, \ldots, s_k)$ must end with $x_1$. Following number theorists, we call admissible all words starting with $x_0$ and ending with $x_1$ and  set for any admissible word $w=x_{\underline s}$
$$\tilde \zeta_1^r (x_{\underline s})= \tilde \zeta_1^r (\underline s).$$
Clearly, restricting  to admissible words amounts to considering only convergent cases. \\
Let us introduce the $\Q$-algebra of polynomials $\Q\langle X\rangle=\Q\langle x_0,x_1\rangle$  in two non-commutative variables $x_0$ and $x_1$ graded by the degree, with $x_i$ of degree $1$. $\Q\langle X\rangle$ is identified with the graded $\Q$-vector space ${\cal H}$ spanned by the monomials in the variables $x_0$ and $x_1$. We equip ${\cal H}$
with the same shuffle product as in the bulk of the paper, defined inductively by: 
$$1\shu w= w\shu 1= w; \quad x_j u\shu x_kv=x_j(u\shu x_k v)+ x_k(x_j u\shu v).$$
For any word $w=x_{\underline s}$ ending with $x_1$ we set \footnote{Here, as long as  $\vert z\vert <1$,  the restriction on the first letter needed for the convergence of the multiple zeta functions is not necessary anymore.}: 
$$\tilde {Li}^r_{x_{\underline s}}(z):=  \tilde {Li}^r_{\underline s}(z).$$
Then the above result translates to: 
\begin{equation}\label{eq:tildeLi}
\frac{d}{dz} \tilde{Li}^r_{x_j \, u}(z)= \omega_j(z)\,  \tilde{Li}^r_{u}(z),
\end{equation} for any non-empty word $u$ and where 
\begin{enumerate}
\item $\omega_j(z)=-\frac{1}{z\log z}$ if $x_j=x_1$ which corresponds to $s_1=1$,
\item and  $\omega_j(z)=\frac{1}{z}$ if $x_j=x_0$  which corresponds to $s_1>1$.
\end{enumerate}Note that
\begin{eqnarray*}
\frac{d}{dz}\tilde {Li}^r_{x_1}(z) &=&\frac{d}{dz}\int_r^\infty
\frac{z^t}{t}dt=\int_r^\infty z^{t-1}dt\\
&=&\left[\frac{z^{t-1}}{\log z}\right]_{t=r}^\infty
=-\frac{z^{r-1}}{\log z}.
\end{eqnarray*}
These results motivate the following  alternative  description of generalised
polylogarithms extended to all words: let us set  $\tilde {Li}^r_1(z)=1$ where the subscript
$1$ denotes the empty word. On the grounds of equation (\ref{eq:tildeLi}) we can inductively   extend the
polylogarithms from words ending with $x_1$ 
to all words defining them as primitives that vanish at $z=0$ (as does the
original expression (\ref{eq:tildeLi})) for words $w=x_ju$
containing the letter $x_1$:
$$\tilde {Li}^r_{x_1}(z) :=\int_0^z v^r\omega_1(v),\quad {\rm and}\quad 
 \tilde {Li}^r_{x_j u}(z)=\int_0^z\omega_j(v) \tilde {Li}^r_{u}(v)$$ if
 $u\neq 1$. For words of the form $x_0^s$ we set:
$$\tilde {Li}^r_{x_0^s}(z)=\frac{\log^sz}{s!}.$$
These extended polylogarithms still satisfy equations (\ref{eq:tildeLi}) and
therefore  coincide with the original polylogarithms for words ending with $x_1$. 
\begin{thm}The map $w\mapsto \tilde {Li}^r_{w}(z)$ is a homomorphism of the algebra $\left({\cal H}, \shu\right)$ into continuous functions on $]0, 1[$. Equivalently, 
$$\tilde {Li}^r_{w_1\,\shu w_2}(z)= \tilde {Li}^r_{w_1}(z)\, \tilde
{Li}^r_{w_2}(z).$$
\end{thm}
{\bf Proof:} The proof goes as for usual generalised polylogarithms by a straightforward derivation using (\ref{eq:tildeLi}) (see e.g. Lemma 2 in \cite{Zu}). One shows by induction on the  total degree  $\vert w_1\vert +\vert w_2\vert$ that 
$$\frac{d}{dz}\left(\tilde {Li}^r_{w_1}(z)\tilde {Li}^r_{w_2}(z)\right)=\frac{d}{dz}\left(\tilde {Li}^r_{w_1\,\sshu w_2}(z)\right).$$
Hence $$\tilde {Li}^r_{w_1\,\sshu w_2}(z)= \tilde {Li}^r_{w_1}(z)\, \tilde
{Li}^r_{w_2}(z)+ C. $$ The substitution $z=0$ gives $C=0$ for words containing
$x_1$ and the substitution $z=1$ gives $C=0$ for words built from $x_0$ only. 
\endsquare\\ \\
By  (\ref{eq:tildeLi})  we obtain the following nested integral representation for the
words $w=x_{i_1}\cdots x_{i_k}$ ending by $x_1$ in the domain $0<z<1$:
$$\tilde{Li}^r_w(z)= \int_{0\leq z_{k}\leq \cdots\leq
  z_1\leq z}\omega_{i_1}(z_1) \cdots  \omega_{i_{k-1}}(z_{k-1})\, z_k^r \omega_{i_k}(z_k)\, dz_1\cdots dz_k.$$
If $x_{i_1}\neq x_1$ then this expression is also defined for $z=1$ and
provides an alternative nested integral  representation for the continuous
analog of the multiple zeta
functions:
$$\tilde \zeta^r (w)= \int_{0\leq z_{k}\leq \cdots\leq
  z_1\leq 1}\omega_{i_1}(z_1) \cdots\omega_{i_{k-1}}(z_{k-1})\,
z_k^r\omega_{i_k}(z_k)\, dz_1\cdots dz_k,$$
which therefore obey the following shuffle relations:
\begin{equation}\label{eq:tildezetashuffle2}
\tilde \zeta^r(w_1\,\shu w_2)= \tilde \zeta^r (w_1)\, \tilde
\zeta^r (w_2).
\end{equation}
\vfill \eject \noindent                                                                    
\bibliographystyle{plain}

\end{document}